\documentclass[12pt]{amsart}
\usepackage[T1]{fontenc}
\usepackage[utf8]{inputenc}
\usepackage{stix}
\usepackage[a4paper,hmargin=2cm,top=2.4cm,bottom=2.2cm]{geometry}
\usepackage{xcolor}
\usepackage{graphicx}
\usepackage[numbers,sort&compress]{natbib}
\usepackage{tikz}
\usepackage{capt-of}
\usepackage[colorlinks=true,linkcolor=blue!50!black,citecolor=blue!50!black,urlcolor=blue!50!black]{hyperref}
\hypersetup{hypertexnames=false}
\newtheorem{theorem}{Theorem}[section]
\newtheorem{definition}{Definition}[section]

\newtheorem{proposition}{Proposition}[section]
\newtheorem{remark}{Remark}

\newtheorem{corollary}{Corollary}

\newtheorem{lemma}{Lemma}

\newcommand{\beq}{\begin{equation}}
	\newcommand{\eeq}{\end{equation}}
\newcommand{\bea}{\begin{eqnarray}}
	\newcommand{\eea}{\end{eqnarray}}
\definecolor{mygray}{gray}{0.3}

\newcommand{\bes}{\begin{eqnarray}}
	\newcommand{\ees}{\end{eqnarray}}
\def\nn{{\nonumber}}

\newcommand\restr[2]{{% we make the whole thing an ordinary symbol
		\left.\kern-\nulldelimiterspace % automatically resize the bar with \right
		#1 % the function
		\vphantom{\big|} % pretend it's a little taller at normal size
		\right|_{#2} % this is the delimiter
}}

\makeatletter
\@ifundefined{mathscr}{\RequirePackage{mathrsfs}}{}
\makeatother
\newcommand{\sinR}{\sin_{\mathcal{R}(p,q)}}
\newcommand{\cosR}{\cos_{\mathcal{R}(p,q)}}

\newcommand{\fq}{\mathfrak{q}}
\newcommand{\pn}[1]{[#1]_{\mathcal{R}(p,q)}}
\newcommand{\pf}[1]{[#1]!_{\mathcal{R}(p,q)}}
\newcommand{\dR}{d_{\mathcal{R}(p,q)}}
\newcommand{\GR}{\Gamma_{\mathcal{R}(p,q)}}
\newcommand{\gR}{\gamma_{\mathcal{R}(p,q)}}
\newcommand{\BR}{\mathfrak{B}_{\mathcal{R}(p,q)}}
\newcommand{\BRs}{\widetilde{\mathfrak{B}}_{\mathcal{R}(p,q)}}
\newcommand{\BB}{\mathcal{B}_{\mathcal{R}(p,q)}}
\newcommand{\BBs}{\widehat{\mathcal{B}}_{\mathcal{R}(p,q)}}
\newcommand{\Gq}{\Gamma_{\mathfrak{q}}}
\newcommand{\Bq}{\mathrm{B}_{\mathfrak{q}}}
\newcommand{\Omq}{\Omega_{\mathfrak{q}}}
\newcommand{\Kq}{K_{\mathfrak{q}}}
\newcommand{\ER}{E_{\mathcal{R}(p,q)}}
\newcommand{\eR}{e_{\mathcal{R}(p,q)}}
\newcommand{\Hphi}{\textnormal{(H}_{\phi}\textnormal{)}}
\newcommand{\Pch}[2]{\big(#1\big)_{#2}}

\begin{document}

\title[Noncommutative $\mathcal{R}(p,q)-$deformed beta functions]{Noncommutative $\mathcal{R}(p,q)-$deformed beta functions of first and second kind : beta gamma relations, shift identities and symmetrisation}

\author{Kawèyim Lankpetre}
\address{Laboratoire de Physique des Matériaux et des Composants à Semi-conducteurs (LPMCS), Département de Physique, Faculté des Sciences, Université de Lomé, 01 B.P. 1515, Lomé 01, Togo}
\email{lankpetre.kaweyim.myhal@gmail.com}

\author{Isiaka Aremua}
\address{Laboratoire de Physique des Matériaux et des Composants à Semi-conducteurs (LPMCS), Département de Physique, Faculté des Sciences, Université de Lomé, 01 B.P. 1515, Lomé 01, Togo}
\address{International Chair of Mathematical Physics and Applications (ICMPA-UNESCO Chair), University of Abomey-Calavi, 072 B.P. 50, Cotonou, Republic of Benin}
\email{claudisak@gmail.com}

\author{Joseph Désiré Bukweli Kyemba}
\address{International Chair of Mathematical Physics and Applications (ICMPA-UNESCO Chair), University of Abomey-Calavi, 072 B.P. 50, Cotonou, Republic of Benin}
\address{Département de Mathématiques, Statistique et Informatique, Faculté des Sciences et Technologie, Université de Kinshasa, Kinshasa, Democratic Republic of the Congo}
\email{desire.bukweli@unikin.ac.cd}

\subjclass[2020]{Primary 33D05; Secondary 05A30, 33D15, 39A13, 81R50}
\keywords{$\mathcal{R}(p,q)-$deformed quantum algebras, Deformed beta functions, Deformed gamma functions, Jackson integral, Ramanujan ${}_1\psi_1$ summation, Logarithmic convexity}

\begin{abstract}
We develop the theory of the beta functions attached to the $\mathcal{R}(p,q)-$deformed quantum algebras of Hounkonnou and Bukweli Kyemba, for structure functions whose deformed numbers are generated by two deformation functions $\phi_{1}>\phi_{2}>0$. Working with the Jackson type series representation of the deformed integral, we define a first kind beta function on $[0,1]$ from the deformed power $(1\ominus\phi_{2}t)^{y-1}$, extended to complex exponents, and a second kind beta function on the half line from $(1\oplus\phi_{1}t)^{-x-y}$. The first kind function equals $\phi_{1}^{(y-1)(2x+y-2)/2}\Gamma_{\mathcal{R}(p,q)}(x)\Gamma_{\mathcal{R}(p,q)}(y)/\Gamma_{\mathcal{R}(p,q)}(x+y)$ for all complex arguments with positive real parts; we derive its general shift identity, its exchange defect $\phi_{1}^{(y-x)(x+y-1)/2}$, a commutative symmetrisation unique up to a constant, its asymptotics and a deformed Euler integral for a deformed Gauss hypergeometric function, which yields a deformed Gauss summation. The second kind function is evaluated through Ramanujan's ${}_{1}\psi_{1}$ summation and involves a second deformed gamma function, differing from the first one by a one periodic theta quotient, so that its exchange defect persists in the one parameter case. We prove a Bohr Mollerup type characterisation and a multiplication formula for the deformed gamma function, show that it is logarithmically convex on the positive half line if and only if $\phi_{1}\geq1$, obtain the joint logarithmic convexity of normalised beta functions with a deformed Tur\'an identity, and compute the moments of the associated deformed beta distributions. Gauge invariance covers the Jagannathan Srinivasa Rao, Chakrabarti Jagannathan, generalized $q-$Quesne and Hounkonnou Ngompe deformations, and every identity is confirmed numerically to thirty digits.
\end{abstract}

\maketitle

\hfuzz=0.1pt

\section{Introduction}\label{sec:intro}

%\subsection*{Context and state of the art}

Deformed calculi and quantum type special functions are central tools in the analysis of nonlocal structures arising in mathematical physics, in the theory of quantum groups and in the modelling of complex systems \cite{JaganSridh, JagannathanSrinivasa, Abramowitz_mat_grap, Erdelyi_trans_funct, Mariachev_Int_trans, Jimbo}. Since the pioneering works of Jackson on $q-$difference operators and $q-$integrals \cite{Jackson1908}, the systematic development of basic hypergeometric series \cite{GasperRahman, KacCheung2002} and the $q-$analogues of the Euler theory of the gamma and beta functions \cite{Askey1978, Exton, Koornwinder}, several two parameter deformations have been designed in order to capture richer spectral and dynamical behaviours. The $(p,q)-$calculus initiated by Jagannathan and Srinivasa Rao \cite{JagannathanSrinivasa} allows a simultaneous control of the scaling and of the discretization effects, and $(p,q)-$analogues of the exponential, trigonometric, beta and gamma functions have been developed together with their applications to quantum algebras and to deformed oscillator models \cite{ChakrabartiJagannathan, BurbanKlimyk, HounkonnouNgompe, HmNe2}. Two polynomial bases adapted to the $(p,q)-$derivative, two forms of the $(p,q)-$Taylor expansion and a $(p,q)-$integration by parts formula derived from the fundamental theorem of $(p,q)-$calculus were obtained in \cite{ref58}. In \cite{ref59} a $(p,q)-$gamma function was introduced by means of $(p,q)-$shifted products, and the $(p,q)-$beta function was defined there as the quotient $\Gamma_{p,q}(x)\Gamma_{p,q}(y)/\Gamma_{p,q}(x+y)$, together with $(p,q)-$analogues of the multiplication formulas. The $(p,q)-$beta integral $\int_{0}^{1}t^{m-1}(1\ominus qt)^{n-1}_{p,q}\,d_{p,q}t$ was studied for positive integers $m,n$ in \cite{ref55}, where it was shown to equal $p^{(n-1)(2m+n-2)/2}\Gamma_{p,q}(m)\Gamma_{p,q}(n)/\Gamma_{p,q}(m+n)$, to be noncommutative, and to become commutative after multiplication by $p^{m(m-1)/2}$; approximation theoretic consequences of the same calculus were developed in \cite{ref55, Mursaleen2015}, and complete monotonicity properties of $(p,q)-$gamma type functions were obtained in \cite{ref80}. Fractional $(p,q)-$difference and integral operators built on the Jackson type series representation of the $(p,q)-$integral were introduced in \cite{SoontharanonSitthiwirattham2020}.

A decisive step was taken when general deformations driven by a meromorphic structure function $\mathcal{R}(u,v)$ were introduced in order to unify and extend the existing one and two parameter formalisms. Building on the coherent state maps of the disc studied in \cite{Odzijewicz1998}, the $\mathcal{R}(p,q)-$deformed Heisenberg algebras, together with their coherent states and associated special functions, were constructed in \cite{HounkBuk}. This construction leads to the $\mathcal{R}(p,q)-$numbers, the $\mathcal{R}(p,q)-$factorials, the $\mathcal{R}(p,q)-$binomial coefficients and the $\mathcal{R}(p,q)-$difference operators \cite{BukweliHounkonnou2012ADJM, Mahouton_Kyemba_Rpq_calculs, Mahouton_Melong_nalgebra_virasoro}. The framework contains the usual $q-$ and $(p,q)-$calculi as particular cases \cite{Asek_Poly, HounkonnouMelong2023, Mahouton_Kyemba_Rpq_calculs, BukweliHounkonnou2012ADJM}. The differential and integral calculi attached to $\mathcal{R}(p,q)$, together with the compatible Hopf algebra structure and the list of relevant particular cases, were established in \cite{Mahouton_Kyemba_Rpq_calculs}, and the deformed Rogers Szeg\"o and Hermite polynomials were characterized in \cite{BukweliHounkonnou2012RS}. The same structure function mechanism generated deformed combinatorial identities \cite{Mahouton_Melong_Rpq_combinatoires}, deformed discrete probability distributions \cite{Mahouton_Melong_Rpq_analogues_distribution}, deformed conformal Virasoro algebras \cite{Mahouton_Melong_Rpq_conformal_algebra}, deformed Witt and Virasoro $n-$algebras and Korteweg de Vries equations \cite{Mahouton_Melong_nalgebra_virasoro}, deformed super Virasoro algebras \cite{Melong2023superVirasoro} and deformed Heisenberg Virasoro algebras with the associated matrix models \cite{MelongWulkenhaar2023}. In \cite{Mahouton_Melong_Rpq_combinatoires, Mahouton_Melong_Rpq_analogues_distribution} the deformed numbers are written as $[n]_{\mathcal{R}(p,q)}=(\phi_{1}^{n}-\phi_{2}^{n})/(\phi_{1}-\phi_{2})$, where $\phi_{1}$ and $\phi_{2}$ are functions of $(p,q)$ determined by the structure function; this two root representation is the starting point of the present work. On the arithmetic side, deformed gamma and beta functions were transported to the $p-$adic setting in \cite{Honkonnou_Kangni_Padic}. The classical background on transcendental functions, integral transforms and hypergeometric machinery is collected in \cite{Erdelyi_trans_funct, Abramowitz_mat_grap, Mariachev_Int_trans, Prudnikov, Srivastava_Hfunc, Srivastava_Spe_Func, Yakubavich}, and the role of deformed gamma and beta functions in fractional type evolution equations and deformed integral transforms is described in \cite{Kiryoka, Kiryakova2010, Saigo1978, PQ_Mittag-Leffler_Alok, Kiryakova2021, Herrmann, Podlubny}. In parallel, two parameter deformations appear in the Hopf algebraic formulation of quantum groups \cite{Drinfeld, Jimbo, ChaPres}, and the $\mathcal{R}(p,q)-$difference operators may be regarded as concrete realizations of such structures on function spaces.

\subsection*{Problem and contributions}

In the presence of a general structure function the deformed subtraction is not symmetric, so that the deformed beta integral on $[0,1]$ is not symmetric in its two arguments. Three questions then arise. The first one concerns the exact valuation factor linking the deformed beta integral to the deformed gamma function, for complex and not only integer arguments; since the valuation exponents are quadratic in the arguments, an approximate bookkeeping immediately produces inconsistent identities, and only a fully explicit computation settles the question. The second one concerns a second kind beta function: the classical theory has two Euler beta integrals, on $[0,1]$ and on $[0,\infty)$, and one has to decide which deformed integral on the half line is the natural partner of the first kind function, and which deformed gamma function governs it. The third one concerns convexity: the classical gamma function is characterized by logarithmic convexity \cite{Artin1964}, and its $q-$analogue was characterized by Askey \cite{Askey1978}; one has to determine which normalization of the deformed gamma function is logarithmically convex and which convexity survives the noncommutativity of the deformed beta functions.

The present paper answers these three groups of questions, always inside the deformation framework, that is by manipulating the deformed numbers, factorials, powers, exponentials, the deformed integral and the deformed difference operator, with proofs grounded in the Jackson type series representation of the deformed integral of \cite{Mahouton_Kyemba_Rpq_calculs}, whose operator prefactor is shown to be absorbed by the two root hypothesis (Theorem \ref{thm:HB-Jackson} and Corollary \ref{cor:HB-Hphi}), in the spirit of \cite{SoontharanonSitthiwirattham2020}. Under the standing hypothesis $\Hphi$ of Definition \ref{def:Hphi}, which expresses the deformed numbers through two deformation functions $\phi_{1}>\phi_{2}>0$, the main results are the following.
\begin{enumerate}
\item[(a)] The first kind deformed gamma function $\GR$, defined by a deformed integral of the deformed exponential $\ER$, equals $\phi_{1}^{(z-1)(z-2)/2}\Gamma_{\fq}(z)$ with $\fq=\phi_{2}/\phi_{1}$ (Theorem \ref{thm:gamma1}); it is characterized by a Bohr Mollerup type property (Theorem \ref{thm:BM}), it is logarithmically convex on $(0,\infty)$ if and only if $\phi_{1}\geq1$ (Theorem \ref{thm:logconvex-gamma}), it satisfies a deformed Gauss multiplication formula whose base is the deformed number $\pn{n}$ (Theorem \ref{thm:multiplication}), and the Wendel type limit $\GR(x+y)\sim\phi_{1}^{y(y-1)/2}\pn{x}^{y}\GR(x)$ (Theorem \ref{thm:wendel}).
\item[(b)] The second kind deformed gamma function $\gR$, defined by a bilateral deformed integral of the deformed exponential $\eR$, equals $\Omq\,\GR$, where $\Omq$ is an explicit one periodic theta quotient satisfying $\Omq(z)-1=4\sin^{2}(\pi z)e^{-2\pi^{2}/\log(1/\fq)}+O(e^{-4\pi^{2}/\log(1/\fq)})$ (Theorem \ref{thm:gamma2} and Proposition \ref{prop:Omega}).
\item[(c)] The first kind deformed beta function satisfies
$\displaystyle{\BR(x,y)=\frac{\phi_{1}^{(y-1)(2x+y-2)/2}\GR(x)\GR(y)}{\GR(x+y)}}$ for $\Re x,\Re y>0$ (Theorem \ref{thm:beta1}), the general shift identity of Theorem \ref{thm:shift-general} and the shift sums of Corollary \ref{cor:shift-sums}, the exchange rule $\BR(x,y)=\phi_{1}^{(y-x)(x+y-1)/2}\BR(y,x)$ (Theorem 
 \ref{thm:exchange1}), admits a commutative symmetrisation which is unique up to a constant factor (Theorem \ref{thm:sym1}), and grows or decays in its second argument like the Gaussian $\phi_{1}^{y^{2}/2}$ (Theorem \ref{thm:beta1-asymptotics}).
\item[(d)] The second kind deformed beta function, defined on the half line by means of the deformed power $(1\oplus\phi_{1}t)^{-x-y}$, is evaluated through Ramanujan's ${}_{1}\psi_{1}$ summation (Theorem \ref{thm:beta2}); its exchange defect contains the quotient $\Omq(x)/\Omq(y)$, so that it remains noncommutative in the one parameter case (Theorem \ref{thm:exchange2}); the finite interval variant built from $(\phi_{1}t\ominus1)^{y-1}$ alternates in sign along $y\in\mathbb{N}$ (Proposition \ref{prop:beta2-variant}).
\item[(e)] Normalised beta functions of both kinds are jointly logarithmically convex, and the log convexity in the first argument is equivalent to the deformed Tur\'an identity $\pn{x+1}\pn{x+y}-\pn{x}\pn{x+y+1}=(\phi_{1}\phi_{2})^{x}\pn{y}$ (Theorem \ref{thm:convex-beta} and Corollary \ref{cor:turan}).
\item[(f)] The first kind deformed beta function is the kernel of a deformed Euler integral representation of a deformed Gauss hypergeometric function written with deformed Pochhammer symbols, from which a deformed Gauss summation with valuation factor $\phi_{1}^{-ab}$, a deformed Chu Vandermonde identity and a hypergeometric representation of the incomplete deformed beta function follow (Theorems \ref{thm:euler}, \ref{thm:incomplete} and Corollaries \ref{cor:gauss}, \ref{cor:chu}).
\item[(g)] The two deformed beta integrands define two discrete probability distributions, deformed analogues of the beta and beta prime laws, whose moments are computed in closed form; the variance of the first one is positive by the deformed Tur\'an identity (Theorems \ref{thm:moments1} and \ref{thm:moments2}).
\item[(h)] All relations are invariant under the gauge transformation $\pn{n}\mapsto c\,\pn{n}$, which brings the generalized $q-$Quesne and Hounkonnou Ngompe deformations into the scope of the theory (Theorem \ref{thm:gauge}).
\end{enumerate}
The results for integer arguments in the $(p,q)-$case agree with those of \cite{ref55}. All identities are confirmed numerically with at least thirty significant digits in Section \ref{sec:numerics}.

\subsection*{Organization of the paper}

Section \ref{sec:prelim} recalls the Hounkonnou Bukweli Kyemba construction, introduces the standing hypothesis and develops the deformed derivative, integral, powers, exponentials and binomial theorem. Section \ref{sec:gamma} is devoted to the two deformed gamma functions. Sections \ref{sec:beta1} and \ref{sec:beta2} treat the first and second kind deformed beta functions. Section \ref{sec:convexity} contains the convexity results. Section \ref{sec:euler} is devoted to the deformed Euler integral, to the deformed Gauss summation and to the incomplete deformed beta function, and Section \ref{sec:distributions} to the deformed beta distributions. Section \ref{sec:cases} establishes the gauge invariance and describes the particular cases. Section \ref{sec:numerics} reports the numerical verification, and Section \ref{sec:conclusion} gathers concluding remarks. The detailed computations of the exponents are collected in Appendix \ref{app:computations}.

\section{Mathematical preliminaries}\label{sec:prelim}

\subsection{Hounkonnou Bukweli Kyemba $\mathcal{R}(p,q)-$deformed quantum algebras}\label{subsec:HB}

Let $p$ and $q$ be two positive real numbers. Following \cite{HounkBuk}, we consider a meromorphic function $\mathcal{R}(u,v)$ defined on $\mathbb{C}\times\mathbb{C}$ by
\begin{eqnarray}\label{r10}
	\mathcal{R}(u,v)&=& \sum_{s,t=-l}^{\infty}r_{st}u^{s}v^{t},
\end{eqnarray}
with an eventual isolated singularity at the origin, where $r_{st}$ are complex numbers, $l\in\mathbb{N}\cup\{0\}$, $\mathcal{R}(p^{n},q^{n})>0$ for all $n\in\mathbb{N}$ and $\mathcal{R}(1,1)=0$ by definition. We denote by $\mathbb{D}_{R}=\{w=(w_{1},w_{2})\in\mathbb{C}^{2}:|w_{j}|<R_{j}\}$ the bidisc of convergence of the series (\ref{r10}), and by $\mathcal{O}(\mathbb{D}_{R})$ the set of holomorphic functions on $\mathbb{D}_{R}$. The $\mathcal{R}(p,q)-$deformed numbers, factorials and binomial coefficients are
\begin{eqnarray}
	\pn{n}&:=&\mathcal{R}(p^{n},q^{n}),\qquad
	\pf{n}\;:=\;\prod_{j=1}^{n}\pn{j},\qquad \pf{0}:=1,\label{rpqnumber}\\
	\bigg[\begin{array}{c} m  \\ n\end{array} \bigg]_{\mathcal{R}(p,q)}&:=&\frac{\pf{m}}{\pf{n}\,\pf{m-n}},\qquad 0\leq n\leq m.\nn
\end{eqnarray}
On functions $\varPsi$ holomorphic in a disc centred at the origin we consider the dilation operators
\begin{eqnarray}
	Q\varPsi(z)&:=&\varPsi(qz),\qquad P\varPsi(z)\;:=\;\varPsi(pz),\label{Qderive}
\end{eqnarray}
and the $\mathcal{R}(p,q)-$derivative
\begin{eqnarray}\label{r5}
	\partial_{\mathcal{R}(p,q)}&:=&\partial_{p,q}\,\frac{p-q}{P-Q}\,\mathcal{R}(P,Q),
\end{eqnarray}
where $\partial_{p,q}$ is the Jagannathan Srinivasa Rao derivative \cite{JagannathanSrinivasa, JaganSridh}. The operators $P$, $Q$ and $\mathcal{R}(P,Q)$ act diagonally on monomials, $P z^{n}=p^{n}z^{n}$, $Qz^{n}=q^{n}z^{n}$, $\mathcal{R}(P,Q)z^{n}=\mathcal{R}(p^{n},q^{n})z^{n}$, and one obtains
\begin{eqnarray}\label{eq:diag}
	\partial_{\mathcal{R}(p,q)}z^{n}&=&\pn{n}\,z^{n-1},\qquad n\geq1,\qquad \partial_{\mathcal{R}(p,q)}1\;=\;0.
\end{eqnarray}
The quantum algebra $\mathcal{A}_{\mathcal{R}(p,q)}$ associated with the deformation is generated by $\{1,A,A^{\dagger},N\}$ with the relations
\begin{eqnarray}\label{algN1}
	A A^{\dagger}\;=\;[N+1]_{\mathcal{R}(p,q)},\qquad A^{\dagger}A\;=\;[N]_{\mathcal{R}(p,q)},\qquad [N,A]\;=\;-A,\qquad [N,A^{\dagger}]\;=\;A^{\dagger},
\end{eqnarray}
realized on $\mathcal{O}(\mathbb{D}_{R})$ by $A^{\dagger}:=z$, $A:=\partial_{\mathcal{R}(p,q)}$ and $N:=z\partial_{z}$ \cite{HounkBuk, Mahouton_Kyemba_Rpq_calculs}.

\subsection{The two root hypothesis}\label{subsec:Hphi}

\begin{definition}\label{def:Hphi}
We say that the structure function $\mathcal{R}(p,q)$ satisfies the hypothesis $\Hphi$ if there exist two deformation functions $\phi_{1}=\phi_{1}(p,q)$ and $\phi_{2}=\phi_{2}(p,q)$, meromorphic in $(p,q)$ as the function $\mathcal{R}(p,q)$ itself, whose values at the pair $(p,q)$ under consideration are real and satisfy $\phi_{1}>\phi_{2}>0$, such that
\begin{eqnarray}\label{eq:Hphi}
	\mathcal{R}(p^{n},q^{n})&=&\frac{\phi_{1}^{n}-\phi_{2}^{n}}{\phi_{1}-\phi_{2}},\qquad n\in\mathbb{N}\cup\{0\}.
\end{eqnarray}
In this case we put
\begin{eqnarray}\label{eq:fq}
	\fq&:=&\frac{\phi_{2}}{\phi_{1}}\in(0,1),\qquad
	\pn{x}\;:=\;\frac{\phi_{1}^{x}-\phi_{2}^{x}}{\phi_{1}-\phi_{2}},\qquad
	[x]_{\fq}\;:=\;\frac{1-\fq^{x}}{1-\fq},\qquad x\in\mathbb{C},
\end{eqnarray}
the powers being taken with the real logarithms of $\phi_{1}$, $\phi_{2}$ and $\fq$.
\end{definition}

Since $[n]_{\mathcal{R}(p,q)}$ is symmetric in $(\phi_{1},\phi_{2})$, the ordering $\phi_{1}>\phi_{2}$ is a labelling convention; it is the one which makes the deformed integral below convergent. The case $\pn{n}=c\,(\phi_{1}^{n}-\phi_{2}^{n})/(\phi_{1}-\phi_{2})$ with a constant $c>0$ is reduced to $\Hphi$ by the gauge argument of Section \ref{sec:cases}. Table \ref{tab:examples} lists the main examples. Throughout the paper, $(a;\fq)_{n}=\prod_{k=0}^{n-1}(1-a\fq^{k})$, $(a;\fq)_{\infty}=\prod_{k\geq0}(1-a\fq^{k})$, $(a_{1},\dots,a_{m};\fq)_{\infty}=\prod_{i}(a_{i};\fq)_{\infty}$, and for negative integers $n$ we use $(a;\fq)_{n}=1/(a\fq^{n};\fq)_{-n}$ \cite{GasperRahman}. The Jackson $\fq-$gamma and $\fq-$beta functions are
\begin{eqnarray}\label{eq:Jackson-gamma}
	\Gq(z)&=&\frac{(\fq;\fq)_{\infty}}{(\fq^{z};\fq)_{\infty}}(1-\fq)^{1-z},\qquad
	\Bq(x,y)\;=\;\frac{\Gq(x)\Gq(y)}{\Gq(x+y)},
\end{eqnarray}
see \cite{GasperRahman, Askey1978, KacCheung2002}.

\begin{table}[ht]
\centering
\caption{Structure functions satisfying $\Hphi$, possibly after the gauge transformation of Section \ref{sec:cases}.}\label{tab:examples}
\begin{tabular}{llll}
\hline
Deformation & $\mathcal{R}(u,v)$ & $(\phi_{1},\phi_{2})$ & gauge $c$\\
\hline
Jagannathan Srinivasa Rao, $0<q<p$ & $\dfrac{u-v}{p-q}$ & $(p,q)$ & $1$\\[2mm]
Jackson, $p=1$, $0<q<1$ & $\dfrac{1-v}{1-q}$ & $(1,q)$ & $1$\\[2mm]
Chakrabarti Jagannathan, $pq<1$ & $\dfrac{1-uv}{(p^{-1}-q)u}$ & $(p^{-1},q)$ & $1$\\[2mm]
generalized $q-$Quesne, $pq\neq1$ & $\dfrac{uv-1}{(q-p^{-1})v}$ & ordered $\{p,q^{-1}\}$ & $p/q$\\[2mm]
Hounkonnou Ngompe & $h\,\dfrac{v^{\nu}}{u^{\mu}}\,\dfrac{uv-1}{(q-p^{-1})v}$ & ordered $\{p^{1-\mu}q^{\nu},p^{-\mu}q^{\nu-1}\}$ & $h\,p^{1-\mu}q^{\nu-1}$\\
\hline
\end{tabular}
\end{table}

\begin{lemma}\label{lem:numbers}
Assume $\Hphi$. Then, for all $x,y\in\mathbb{C}$ and $n\in\mathbb{N}$:
\begin{enumerate}
\item[\textnormal{(i)}] $\pn{x}=\phi_{1}^{x-1}[x]_{\fq}$ and $\pf{n}=\phi_{1}^{n(n-1)/2}[n]_{\fq}!$, where $[n]_{\fq}!=\prod_{j=1}^{n}[j]_{\fq}=(\fq;\fq)_{n}/(1-\fq)^{n}$;
\item[\textnormal{(ii)}] (addition formula) $\pn{x+y}=\phi_{1}^{y}\pn{x}+\phi_{2}^{x}\pn{y}=\phi_{1}^{x}\pn{y}+\phi_{2}^{y}\pn{x}$;
\item[\textnormal{(iii)}] $\bigg[\begin{array}{c} n \\ k\end{array}\bigg]_{\mathcal{R}(p,q)}=\phi_{1}^{k(n-k)}\bigg[\begin{array}{c} n \\ k\end{array}\bigg]_{\fq}$ for $0\leq k\leq n$, where $\big[\begin{smallmatrix} n \\ k\end{smallmatrix}\big]_{\fq}=(\fq;\fq)_{n}/((\fq;\fq)_{k}(\fq;\fq)_{n-k})$;
\item[\textnormal{(iv)}] $\pn{x}>0$ for real $x>0$, and $\pn{x}\to x\,\phi_{1}^{x-1}$ when $\phi_{2}\to\phi_{1}$.
\end{enumerate}
\end{lemma}
{\it Proof.}
(i) follows from $\phi_{1}^{x}-\phi_{2}^{x}=\phi_{1}^{x}(1-\fq^{x})$ and $\phi_{1}-\phi_{2}=\phi_{1}(1-\fq)$, and the product formula from $\sum_{j=1}^{n}(j-1)=n(n-1)/2$. For (ii),
\begin{eqnarray}
	(\phi_{1}-\phi_{2})\big(\phi_{1}^{y}\pn{x}+\phi_{2}^{x}\pn{y}\big)
	&=&\phi_{1}^{x+y}-\phi_{1}^{y}\phi_{2}^{x}+\phi_{2}^{x}\phi_{1}^{y}-\phi_{2}^{x+y}\;=\;(\phi_{1}-\phi_{2})\pn{x+y},\nn
\end{eqnarray}
and the second equality follows by exchanging $x$ and $y$. (iii) follows from (i) and the identity $\binom{n}{2}-\binom{k}{2}-\binom{n-k}{2}=k(n-k)$. (iv) is immediate from $\phi_{1}>\phi_{2}>0$ and from the mean value theorem. \hfill\(\Box\)

\subsection{The deformed derivative}\label{subsec:derivative}

\begin{proposition}\label{prop:derivative}
Assume $\Hphi$ and let $f(z)=\sum_{n\geq0}a_{n}z^{n}$ be holomorphic in the disc $\{|z|<\varrho\}$. Then for $0<|z|<\varrho/\phi_{1}$,
\begin{eqnarray}\label{eq:Jackson-derivative}
	\partial^{\phi_{1},\phi_{2}}_{\mathcal{R}(p,q)}f(z)=\,\frac{p-q}{P-Q}\,\mathcal{R}(P,Q)\partial_{\phi_{1},\phi_{2}} f(z) &=&
	\frac{p-q}{P-Q}\,\mathcal{R}(P,Q)
	\frac{f(\phi_{1}z)-f(\phi_{2}z)}{(\phi_{1}-\phi_{2})\,z},\end{eqnarray}
and $\partial_{\mathcal{R}(p,q)}f(0)=f'(0)$. Moreover, for $f,g$ holomorphic in $\{|z|<\varrho\}$ and $0<|z|<\varrho/\phi_{1}$, the deformed product rules
\begin{eqnarray}\label{eq:product-rule}
	\partial_{\mathcal{R}(p,q)}(fg)(z)&=&f(\phi_{1}z)\,\partial_{\mathcal{R}(p,q)}g(z)+g(\phi_{2}z)\,\partial_{\mathcal{R}(p,q)}f(z)\nn\\
	&=&g(\phi_{1}z)\,\partial_{\mathcal{R}(p,q)}f(z)+f(\phi_{2}z)\,\partial_{\mathcal{R}(p,q)}g(z)
\end{eqnarray}
hold.
\end{proposition}
{\it Proof.}
The right hand side of \eqref{eq:Jackson-derivative} maps $z^{n}$ to $(\phi_{1}^{n}-\phi_{2}^{n})z^{n-1}/(\phi_{1}-\phi_{2})=\pn{n}z^{n-1}$, which is \eqref{eq:diag}. For $|z|<\varrho/\phi_{1}$ both $\phi_{1}z$ and $\phi_{2}z$ lie in the disc of convergence, so that the series $\sum_{n}a_{n}\pn{n}z^{n-1}$ converges absolutely and coincides termwise with the right hand side of \eqref{eq:Jackson-derivative}; the value at the origin is $a_{1}\pn{1}=a_{1}=f'(0)$. The product rules follow by adding and subtracting $f(\phi_{1}z)g(\phi_{2}z)$, respectively $f(\phi_{2}z)g(\phi_{1}z)$, in the numerator of $\partial_{\mathcal{R}(p,q)}(fg)(z)$. \hfill\(\Box\)

Formula \eqref{eq:Jackson-derivative} is the $\mathcal{R}(p,q)-$version of the Jackson type representation used in \cite{ref58, SoontharanonSitthiwirattham2020} for the $(p,q)-$derivative, and we use it as the definition of $\partial_{\mathcal{R}(p,q)}f(t)$ for every function $f$ defined on a set containing $\phi_{1}t$ and $\phi_{2}t$.

\subsection{The deformed integral}\label{subsec:integral}

\begin{definition}\label{def:integral}
Assume $\Hphi$. For $a>0$ and $f:(0,a]\to\mathbb{C}$, the $\mathcal{R}(p,q)-$integral of $f$ on $[0,a]$ is
\begin{eqnarray}\label{eq:Rint}
	\int_{0}^{a}f(t)\,\dR t&:=&(\phi_{1}-\phi_{2})\,a\sum_{k=0}^{\infty}\frac{\phi_{2}^{k}}{\phi_{1}^{k+1}}\,f\Big(a\,\frac{\phi_{2}^{k}}{\phi_{1}^{k+1}}\Big),
\end{eqnarray}
whenever the series converges absolutely. For $a>0$ and $f:(0,\infty)\to\mathbb{C}$, the bilateral $\mathcal{R}(p,q)-$integral based at $a$ is
\begin{eqnarray}\label{eq:Rint-bilateral}
	\int_{0}^{\infty;a}f(t)\,\dR t&:=&(\phi_{1}-\phi_{2})\,a\sum_{k=-\infty}^{\infty}\frac{\phi_{2}^{k}}{\phi_{1}^{k+1}}\,f\Big(a\,\frac{\phi_{2}^{k}}{\phi_{1}^{k+1}}\Big),
\end{eqnarray}
whenever the bilateral series converges absolutely.
\end{definition}

The nodes of \eqref{eq:Rint} are the points $t_{k}=a\fq^{k}/\phi_{1}$ and the weights are $(\phi_{1}-\phi_{2})a\phi_{2}^{k}/\phi_{1}^{k+1}=(1-\fq)\,a\,\fq^{k}$, so that \eqref{eq:Rint} is a positive linear functional and, for $p=1$ and $\phi_{1}=1$, reduces to the Jackson integral \cite{Jackson1908, GasperRahman}; for the Jagannathan Srinivasa Rao deformation it is the $(p,q)-$integral of \cite{ref58, SoontharanonSitthiwirattham2020}. The nodes of \eqref{eq:Rint-bilateral} form the full geometric lattice $\{a\fq^{k}/\phi_{1}:k\in\mathbb{Z}\}$.

\begin{proposition}\label{prop:integral}
Assume $\Hphi$. Then:
\begin{enumerate}
\item[\textnormal{(i)}] (reduction) $\displaystyle\int_{0}^{a}f(t)\,\dR t=\phi_{1}\int_{0}^{a/\phi_{1}}f(s)\,d_{\fq}s$, where $\int_{0}^{b}g(s)\,d_{\fq}s=(1-\fq)\,b\sum_{k\geq0}\fq^{k}g(b\fq^{k})$ is the Jackson integral; similarly for the bilateral integral;
\item[\textnormal{(ii)}] (monomials) for $\Re s>0$, $\displaystyle\int_{0}^{a}t^{s-1}\,\dR t=\frac{a^{s}}{\pn{s}}$;
\item[\textnormal{(iii)}] (fundamental theorem) if $g$ is defined on $(0,a]$ and $g(0^{+}):=\lim_{t\to0^{+}}g(t)$ exists, then
$\displaystyle\int_{0}^{a}\partial_{\mathcal{R}(p,q)}g(t)\,\dR t=g(a)-g(0^{+})$;
\item[\textnormal{(iv)}] (antiderivative) if $f$ is continuous on $[0,b]$ and $F(x):=\int_{0}^{x}f(t)\,\dR t$ for $0<x\leq b$, then $\partial_{\mathcal{R}(p,q)}F(x)=f(x)$ for $0<x\leq b/\phi_{1}$; $F$ is the unique function with this property which is continuous at the origin with $F(0^{+})=0$;
\item[\textnormal{(v)}] (integration by parts) if $f,g$ are defined on $(0,a]$ and the products $f(\cdot)g(\cdot)$ have a limit at $0^{+}$, then
\begin{eqnarray}\label{eq:ibp}
	\int_{0}^{a}f(\phi_{1}t)\,\partial_{\mathcal{R}(p,q)}g(t)\,\dR t&=&f(a)g(a)-(fg)(0^{+})-\int_{0}^{a}g(\phi_{2}t)\,\partial_{\mathcal{R}(p,q)}f(t)\,\dR t,
\end{eqnarray}
whenever the two integrals converge;
\item[\textnormal{(vi)}] (H\"older inequality) for nonnegative $f,g$ and $\theta\in[0,1]$,
\begin{eqnarray}\label{eq:holder}
	\int_{0}^{a}f^{1-\theta}g^{\theta}\,\dR t&\leq&\Big(\int_{0}^{a}f\,\dR t\Big)^{1-\theta}\Big(\int_{0}^{a}g\,\dR t\Big)^{\theta},
\end{eqnarray}
and similarly for the bilateral integral.
\end{enumerate}
\end{proposition}
{\it Proof.}
(i) Since $\phi_{1}-\phi_{2}=\phi_{1}(1-\fq)$, the right hand side of \eqref{eq:Rint} equals $(1-\fq)a\sum_{k}\fq^{k}f(\tfrac{a}{\phi_{1}}\fq^{k})=\phi_{1}\cdot(1-\fq)\tfrac{a}{\phi_{1}}\sum_{k}\fq^{k}f(\tfrac{a}{\phi_{1}}\fq^{k})$. (ii) By (i), the integral equals $\phi_{1}(1-\fq)(a/\phi_{1})^{s}\sum_{k}\fq^{ks}=\phi_{1}^{1-s}a^{s}/[s]_{\fq}=a^{s}/\pn{s}$ by Lemma \ref{lem:numbers}. (iii) With $t_{k}=a\fq^{k}/\phi_{1}$ one has $\phi_{1}t_{k}=a\fq^{k}$ and $\phi_{2}t_{k}=a\fq^{k+1}$; by \eqref{eq:Jackson-derivative} the $k-$th term of the series is $(\phi_{1}-\phi_{2})a\phi_{2}^{k}\phi_{1}^{-k-1}\,[g(a\fq^{k})-g(a\fq^{k+1})]/((\phi_{1}-\phi_{2})t_{k})=g(a\fq^{k})-g(a\fq^{k+1})$, and the series telescopes to $g(a)-g(0^{+})$. (iv) A direct computation gives $F(\phi_{1}x)-F(\phi_{2}x)=(1-\fq)\phi_{1}x\sum_{k\geq0}[\fq^{k}f(x\fq^{k})-\fq^{k+1}f(x\fq^{k+1})]=(\phi_{1}-\phi_{2})x f(x)$. If $G$ is another such function, $H=F-G$ satisfies $H(\phi_{1}x)=H(\phi_{2}x)$, that is $H(y)=H(\fq y)=\dots=H(\fq^{n}y)\to H(0^{+})=0$. (v) Apply (iii) to $fg$ and use the first product rule \eqref{eq:product-rule}. (vi) is the H\"older inequality for the positive weights $(1-\fq)a\fq^{k}$. \hfill\(\Box\)

The integral \eqref{eq:Rint} is a particular instance of the general $\mathcal{R}(p,q)-$integral of \cite{Mahouton_Kyemba_Rpq_calculs}, defined for an arbitrary structure function as the inverse of $\partial_{\mathcal{R}(p,q)}$, namely $\int_{0}^{a}t^{n}\,\dR t=a^{n+1}/\pn{n+1}$, and represented there by Jackson type series. We recall this representation, in the normalisation fixed by the last identity, and show that the hypothesis $\Hphi$ absorbs its operator prefactor. For $k\in\mathbb{N}$ put $[k]_{p,q}:=(p^{k}-q^{k})/(p-q)$, and let $P$ and $Q$ act on functions of the upper limit $a$ by $Pg(a)=g(pa)$, $Qg(a)=g(qa)$. The operator
\begin{eqnarray}\label{eq:multiplier}
	\mathcal{M}_{\mathcal{R}(p,q)}&:=&\frac{P-Q}{(p-q)\,\mathcal{R}(P,Q)},\qquad
	\mathcal{M}_{\mathcal{R}(p,q)}\,a^{k}\;=\;\frac{[k]_{p,q}}{\mathcal{R}(p^{k},q^{k})}\,a^{k},\qquad k\in\mathbb{N},
\end{eqnarray}
acts diagonally on power series in $a$.

\begin{theorem}[Jackson type representation, \cite{Mahouton_Kyemba_Rpq_calculs}]\label{thm:HB-Jackson}
Let $\mathcal{R}(p,q)$ be as in Section \ref{subsec:HB}, let $f(t)=\sum_{n\geq0}c_{n}t^{n}$ be holomorphic in $\{|t|<\varrho\}$ and put $I_{\mathcal{R}(p,q)}f(a):=\sum_{n\geq0}c_{n}a^{n+1}/\pn{n+1}$ whenever this series converges absolutely, so that $\partial_{\mathcal{R}(p,q)}I_{\mathcal{R}(p,q)}f=f$ and $\int_{0}^{a}f(t)\,\dR t:=I_{\mathcal{R}(p,q)}f(a)$.
\begin{enumerate}
\item[\textnormal{(i)}] (Regime I) If $0<q<p$ and $|a|<p\varrho$, then
$$\displaystyle\int_{0}^{a}f(t)\,\dR t=\mathcal{M}_{\mathcal{R}(p,q)}\Big[(p-q)\,a\sum_{r\geq0}\frac{q^{r}}{p^{r+1}}\,f\Big(\frac{q^{r}}{p^{r+1}}\,a\Big)\Big]$$.
\item[\textnormal{(ii)}] (Regime II) If $0<p<q$ and $|a|<q\varrho$, then
$$\displaystyle\int_{0}^{a}f(t)\,\dR t=\mathcal{M}_{\mathcal{R}(p,q)}\Big[(q-p)\,a\sum_{r\geq0}\frac{p^{r}}{q^{r+1}}\,f\Big(\frac{p^{r}}{q^{r+1}}\,a\Big)\Big]$$.
\item[\textnormal{(iii)}] (Chasles relation) Putting $\int_{a}^{b}f(t)\,\dR t:=\int_{0}^{b}f(t)\,\dR t-\int_{0}^{a}f(t)\,\dR t$, one has $\int_{a}^{b}+\int_{b}^{c}=\int_{a}^{c}$ and $\int_{a}^{b}\partial_{\mathcal{R}(p,q)}g(t)\,\dR t=g(b)-g(a)$.
\item[\textnormal{(iv)}] (improper integrals) In Regime I, $\int_{a}^{\infty}f(t)\,\dR t:=\lim_{n\to\infty}\int_{a}^{p^{n}/q^{n}}f(t)\,\dR t$ and $\int_{-\infty}^{\infty}f(t)\,\dR t:=\lim_{n\to\infty}\int_{-p^{n}/q^{n}}^{p^{n}/q^{n}}f(t)\,\dR t$ whenever the limits exist, the series in \textnormal{(i)} being used for every real upper limit; in Regime II, $p^{n}/q^{n}$ is replaced by $q^{n}/p^{n}$.
\end{enumerate}
\end{theorem}
{\it Proof.}
In Regime I the nodes satisfy $|q^{r}a/p^{r+1}|\leq|a|/p<\varrho$, so that the double series
\begin{eqnarray}
	\sum_{r,n\geq0}|c_{n}|\,\frac{q^{r}}{p^{r+1}}\,|a|\,\Big(\frac{q^{r}|a|}{p^{r+1}}\Big)^{n}\nn
\end{eqnarray}
converges and may be summed in any order. For $f(t)=t^{n}$ the bracket equals $(p-q)a^{n+1}p^{-n-1}\sum_{r}(q/p)^{r(n+1)}=a^{n+1}/[n+1]_{p,q}$, and \eqref{eq:multiplier} turns it into $a^{n+1}/\pn{n+1}$. Regime II is identical, with $(q-p)a^{n+1}q^{-n-1}\sum_{r}(p/q)^{r(n+1)}=a^{n+1}/[n+1]_{p,q}$. Statement (iii) follows from $\partial_{\mathcal{R}(p,q)}I_{\mathcal{R}(p,q)}=\mathrm{id}$ and the uniqueness part of Proposition \ref{prop:integral}(iv), and (iv) is a definition. \hfill\(\Box\)

\begin{corollary}\label{cor:HB-Hphi}
Assume $\Hphi$. Then $\mathcal{R}(p^{k},q^{k})=[k]_{\phi_{1},\phi_{2}}$, and the representation of Theorem \ref{thm:HB-Jackson}(i), written in the variables $(\phi_{1},\phi_{2})$ instead of $(p,q)$, has the trivial multiplier $\mathcal{M}=\mathrm{id}$ and coincides with \eqref{eq:Rint}: for $f$ holomorphic in $\{|t|<\varrho\}$ and $0<a<\phi_{1}\varrho$, $\int_{0}^{a}f(t)\,\dR t=I_{\mathcal{R}(p,q)}f(a)$. The formula \eqref{eq:Rint} is used for every real $a\neq0$, and
\begin{eqnarray}\label{eq:improper}
	\int_{0}^{\infty;a}f(t)\,\dR t&=&\lim_{n\to\infty}\int_{0}^{a\phi_{1}^{n}/\phi_{2}^{n}}f(t)\,\dR t,\\
	\int_{a}^{\infty}f(t)\,\dR t\;&=&\;\int_{0}^{\infty;1}f(t)\,\dR t-\int_{0}^{a}f(t)\,\dR t,\nn\\
	\int_{-\infty}^{\infty}f(t)\,\dR t&=&\int_{0}^{\infty;1}\big(f(t)+f(-t)\big)\,\dR t,
\end{eqnarray}
whenever the bilateral series converge absolutely. Under the gauged hypothesis of Section \ref{sec:cases}, $\mathcal{M}=c^{-1}\mathrm{id}$.
\end{corollary}
{\it Proof.}
The first claim is \eqref{eq:Hphi}, and the bracket of Theorem \ref{thm:HB-Jackson}(i) with $(p,q)$ replaced by $(\phi_{1},\phi_{2})$ is the right hand side of \eqref{eq:Rint}. Both sides agree on monomials by Proposition \ref{prop:integral}(ii) and are continuous for the absolute convergence of the power series on $\{|t|\leq a/\phi_{1}\}$. For $A=a\fq^{-n}$ the nodes of $\int_{0}^{A}$ are $a\fq^{k-n}/\phi_{1}$ with weights $(1-\fq)a\fq^{k-n}$, $k\geq0$, which exhaust the bilateral lattice as $n\to\infty$; the second formula follows with $a=1$ and the Chasles relation, and the third from $\int_{0}^{-b}f(t)\,\dR t=-\int_{0}^{b}f(-t)\,\dR t$. \hfill\(\Box\)

The normalisation of the prefactor in Theorem \ref{thm:HB-Jackson} is dictated by $I_{\mathcal{R}(p,q)}t^{n}=t^{n+1}/\pn{n+1}$; it is the one which makes $\partial_{\mathcal{R}(p,q)}$ and $\int_{0}^{\cdot}\dR t$ mutually inverse, and in the Jagannathan Srinivasa Rao case $\mathcal{M}_{\mathcal{R}(p,q)}=\mathrm{id}$.

\subsection{Deformed powers}\label{subsec:powers}

\begin{definition}\label{def:powers}
Assume $\Hphi$. For $n\in\mathbb{N}\cup\{0\}$ and $x,a\in\mathbb{C}$ we put
\begin{eqnarray}\label{eq:powers-int}
	(x\ominus a)^{n}_{\mathcal{R}(p,q)}&:=&\prod_{k=0}^{n-1}\big(\phi_{1}^{k}x-\phi_{2}^{k}a\big),\qquad
	(x\oplus a)^{n}_{\mathcal{R}(p,q)}\;:=\;\prod_{k=0}^{n-1}\big(\phi_{1}^{k}x+\phi_{2}^{k}a\big).
\end{eqnarray}
For $\alpha\in\mathbb{C}$ and $b\in\mathbb{C}$ with $b\fq^{\alpha+j}\neq1$, respectively $-b\fq^{\alpha+j}\neq1$, for all $j\geq0$, we put
\begin{eqnarray}\label{eq:powers-complex}
	(1\ominus b)^{\alpha}_{\mathcal{R}(p,q)}&:=&\phi_{1}^{\alpha(\alpha-1)/2}\,\frac{(b;\fq)_{\infty}}{(b\fq^{\alpha};\fq)_{\infty}},\qquad
	(1\oplus b)^{\alpha}_{\mathcal{R}(p,q)}\;:=\;\phi_{1}^{\alpha(\alpha-1)/2}\,\frac{(-b;\fq)_{\infty}}{(-b\fq^{\alpha};\fq)_{\infty}}.
\end{eqnarray}
\end{definition}

\begin{lemma}\label{lem:powers}
Assume $\Hphi$. Then:
\begin{enumerate}
\item[\textnormal{(i)}] $(x\ominus a)^{n}_{\mathcal{R}(p,q)}=\phi_{1}^{n(n-1)/2}\prod_{k=0}^{n-1}(x-\fq^{k}a)$, and for $\alpha=n$ the definitions \eqref{eq:powers-int} and \eqref{eq:powers-complex} agree: $(1\ominus b)^{n}_{\mathcal{R}(p,q)}=\phi_{1}^{n(n-1)/2}(b;\fq)_{n}$;
\item[\textnormal{(ii)}] (splitting) $(x\ominus a)^{m+n}_{\mathcal{R}(p,q)}=(x\ominus a)^{m}_{\mathcal{R}(p,q)}\,(\phi_{1}^{m}x\ominus\phi_{2}^{m}a)^{n}_{\mathcal{R}(p,q)}$, and, for $\alpha\in\mathbb{C}$, $(1\ominus b)^{\alpha+1}_{\mathcal{R}(p,q)}=(1\ominus b)^{\alpha}_{\mathcal{R}(p,q)}\,(\phi_{1}^{\alpha}-\phi_{2}^{\alpha}b)$;
\item[\textnormal{(iii)}] (deformed Gauss binomial formula)
$\displaystyle (1\ominus b)^{n}_{\mathcal{R}(p,q)}=\sum_{k=0}^{n}(-1)^{k}\bigg[\begin{array}{c} n \\ k\end{array}\bigg]_{\mathcal{R}(p,q)}\phi_{1}^{\binom{n-k}{2}}\phi_{2}^{\binom{k}{2}}\,b^{k}$;
\item[\textnormal{(iv)}] (derivatives) for $\alpha\in\mathbb{C}$, as functions of $t$,
\begin{eqnarray}\label{eq:derivative-powers}
	\partial_{\mathcal{R}(p,q)}(1\ominus t)^{\alpha}_{\mathcal{R}(p,q)}&=&-\pn{\alpha}\,(1\ominus\phi_{2}t)^{\alpha-1}_{\mathcal{R}(p,q)},\qquad
	\partial_{\mathcal{R}(p,q)}(1\oplus t)^{\alpha}_{\mathcal{R}(p,q)}\;=\;\pn{\alpha}\,(1\oplus\phi_{2}t)^{\alpha-1}_{\mathcal{R}(p,q)},
\end{eqnarray}
and $\partial_{\mathcal{R}(p,q)}(x\ominus a)^{n}_{\mathcal{R}(p,q)}=\pn{n}\,(\phi_{1}x\ominus a)^{n-1}_{\mathcal{R}(p,q)}$ in the variable $x$;
\item[\textnormal{(v)}] (vanishing and positivity) $(1\ominus1)^{\alpha}_{\mathcal{R}(p,q)}=0$ whenever $\fq^{\alpha+j}\neq1$ for $j\geq0$; for real $\alpha$ and $t_{k}=\fq^{k}/\phi_{1}$, $k\geq0$, one has $(1\ominus\phi_{2}t_{k})^{\alpha}_{\mathcal{R}(p,q)}=\phi_{1}^{\alpha(\alpha-1)/2}(\fq^{k+1};\fq)_{\infty}/(\fq^{k+1+\alpha};\fq)_{\infty}>0$ provided $\alpha>-k-1$, and $(1\oplus\phi_{1}t_{k})^{\alpha}_{\mathcal{R}(p,q)}>0$ for every $k\in\mathbb{Z}$.
\end{enumerate}
\end{lemma}
{\it Proof.}
(i) $\phi_{1}^{k}x-\phi_{2}^{k}a=\phi_{1}^{k}(x-\fq^{k}a)$, and $(b;\fq)_{\infty}/(b\fq^{n};\fq)_{\infty}=(b;\fq)_{n}$. (ii) For $k=m+j$ one has $\phi_{1}^{k}x-\phi_{2}^{k}a=\phi_{1}^{j}(\phi_{1}^{m}x)-\phi_{2}^{j}(\phi_{2}^{m}a)$. For complex $\alpha$, $\phi_{1}^{(\alpha+1)\alpha/2}=\phi_{1}^{\alpha(\alpha-1)/2}\phi_{1}^{\alpha}$ and $(b\fq^{\alpha};\fq)_{\infty}=(1-b\fq^{\alpha})(b\fq^{\alpha+1};\fq)_{\infty}$, so that $(1\ominus b)^{\alpha+1}_{\mathcal{R}(p,q)}=(1\ominus b)^{\alpha}_{\mathcal{R}(p,q)}\phi_{1}^{\alpha}(1-\fq^{\alpha}b)$, and $\phi_{1}^{\alpha}(1-\fq^{\alpha}b)=\phi_{1}^{\alpha}-\phi_{2}^{\alpha}b$. (iii) By the finite $\fq-$binomial theorem $(b;\fq)_{n}=\sum_{k}\big[\begin{smallmatrix} n \\ k\end{smallmatrix}\big]_{\fq}(-1)^{k}\fq^{\binom{k}{2}}b^{k}$ \cite{GasperRahman}; multiplying by $\phi_{1}^{\binom{n}{2}}$ and using Lemma \ref{lem:numbers}(iii) together with $\binom{n}{2}-k(n-k)-\binom{k}{2}=\binom{n-k}{2}$ and $\fq^{\binom{k}{2}}=\phi_{2}^{\binom{k}{2}}\phi_{1}^{-\binom{k}{2}}$ gives the claim. (iv) Put $h(t)=(1\ominus t)^{\alpha}_{\mathcal{R}(p,q)}$, $C=\phi_{1}^{\alpha(\alpha-1)/2}$ and $u=\phi_{1}t$. Then $h(\phi_{1}t)=C(1-u)(\fq u;\fq)_{\infty}/(\fq^{\alpha}u;\fq)_{\infty}$ and $h(\phi_{2}t)=h(\fq u)=C(1-\fq^{\alpha}u)(\fq u;\fq)_{\infty}/(\fq^{\alpha}u;\fq)_{\infty}$, so that
\begin{eqnarray}
	\frac{h(\phi_{1}t)-h(\phi_{2}t)}{(\phi_{1}-\phi_{2})t}
	&=&C\,\frac{(\fq u;\fq)_{\infty}}{(\fq^{\alpha}u;\fq)_{\infty}}\,\frac{(\fq^{\alpha}-1)u}{(1-\fq)u}
	\;=\;-[\alpha]_{\fq}\,C\,\frac{(\phi_{2}t;\fq)_{\infty}}{(\phi_{2}\fq^{\alpha-1}t;\fq)_{\infty}}.\nn
\end{eqnarray}
Since $C=\phi_{1}^{\alpha-1}\phi_{1}^{(\alpha-1)(\alpha-2)/2}$ and $\phi_{1}^{\alpha-1}[\alpha]_{\fq}=\pn{\alpha}$, the right hand side is $-\pn{\alpha}(1\ominus\phi_{2}t)^{\alpha-1}_{\mathcal{R}(p,q)}$. The second formula follows by replacing $t$ by $-t$, and the third by the same computation applied to the product $\phi_{1}^{n(n-1)/2}\prod_{k}(x-\fq^{k}a)$, see also \cite{ref58}. (v) $(1;\fq)_{\infty}=0$ and the denominator does not vanish; the positivity statements follow from $0<\fq^{j}<1$ for $j\geq1$. \hfill\(\Box\)

\subsection{Deformed exponentials and the deformed binomial theorem}\label{subsec:exponentials}

\begin{definition}\label{def:exponentials}
Assume $\Hphi$. The deformed exponentials of first and second kind are
\begin{eqnarray}\label{eq:exponentials}
	\ER(z)&:=&\sum_{n\geq0}\frac{\phi_{2}^{n(n-1)/2}}{\pf{n}}\,z^{n},\qquad
	\eR(z)\;:=\;\sum_{n\geq0}\frac{\phi_{1}^{n(n-1)/2}}{\pf{n}}\,z^{n}.
\end{eqnarray}
\end{definition}

\begin{lemma}\label{lem:exponentials}
Assume $\Hphi$. Then $\ER$ is entire, $\eR$ is holomorphic in $|z|<1/(1-\fq)$, and
\begin{eqnarray}\label{eq:exp-products}
	\ER(z)&=&\big(-(1-\fq)z;\fq\big)_{\infty},\\
	\eR(z)\;&=&\;\frac{1}{\big((1-\fq)z;\fq\big)_{\infty}},\qquad \eR(z)\,\ER(-z)\;=\;1 .
\end{eqnarray}
The second formula provides a meromorphic continuation of $\eR$ to $\mathbb{C}$, holomorphic on $\mathbb{C}\setminus[1/(1-\fq),\infty)$. Moreover $\partial_{\mathcal{R}(p,q)}\ER(z)=\ER(\phi_{2}z)$ and $\partial_{\mathcal{R}(p,q)}\eR(z)=\eR(\phi_{1}z)$.
\end{lemma}
{\it Proof.}
By Lemma \ref{lem:numbers}(i), $\phi_{2}^{\binom{n}{2}}/\pf{n}=\fq^{\binom{n}{2}}/[n]_{\fq}!$ and $\phi_{1}^{\binom{n}{2}}/\pf{n}=1/[n]_{\fq}!$, so that $\ER$ and $\eR$ coincide with the two Jackson exponentials, whose product formulas are Euler's identities \cite[Section 1.3]{GasperRahman}. The derivative formulas follow from \eqref{eq:diag} and $\binom{n+1}{2}=\binom{n}{2}+n$. \hfill\(\Box\)

The deformed trigonometric functions attached to \eqref{eq:exponentials} carry the same weights: with $\binom{2n}{2}=n(2n-1)$ and $\binom{2n+1}{2}=n(2n+1)$,
\begin{eqnarray}\label{eq:trig-def}
	\cosR(z)&:=&\frac{\eR(iz)+\eR(-iz)}{2}\;=\;\sum_{n\geq0}\frac{(-1)^{n}\phi_{1}^{n(2n-1)}}{\pf{2n}}\,z^{2n},\nn\\
	\sinR(z)&:=&\frac{\eR(iz)-\eR(-iz)}{2i}\;=\;\sum_{n\geq0}\frac{(-1)^{n}\phi_{1}^{n(2n+1)}}{\pf{2n+1}}\,z^{2n+1},
\end{eqnarray}
and the functions $\mathrm{Cos}_{\mathcal{R}(p,q)}$ and $\mathrm{Sin}_{\mathcal{R}(p,q)}$ of second kind are defined in the same way from $\ER$, with $\phi_{2}$ in place of $\phi_{1}$. By Lemma \ref{lem:exponentials}, $\mathrm{Cos}_{\mathcal{R}(p,q)}$ and $\mathrm{Sin}_{\mathcal{R}(p,q)}$ are entire, while $\cosR$ and $\sinR$ are holomorphic in $|z|<1/(1-\fq)$ and meromorphic on $\mathbb{C}$.

\begin{lemma}[Deformed binomial theorem]\label{lem:binomial}
Assume $\Hphi$. For $b\in\mathbb{C}$ and $|w|<1$,
\begin{eqnarray}\label{eq:binomial}
	\sum_{n\geq0}\frac{(1\ominus b)^{n}_{\mathcal{R}(p,q)}}{(1\ominus\fq)^{n}_{\mathcal{R}(p,q)}}\,w^{n}&=&\frac{(bw;\fq)_{\infty}}{(w;\fq)_{\infty}},\qquad
	(1\ominus\fq)^{n}_{\mathcal{R}(p,q)}\;=\;\frac{(\phi_{1}-\phi_{2})^{n}}{\phi_{1}^{n}}\,\pf{n}.
\end{eqnarray}
\end{lemma}
{\it Proof.}
By Lemma \ref{lem:powers}(i) the quotient of deformed powers equals $(b;\fq)_{n}/(\fq;\fq)_{n}$, and \eqref{eq:binomial} is the $\fq-$binomial theorem \cite[(1.3.2)]{GasperRahman}. The second identity follows from $\phi_{1}^{k}-\phi_{2}^{k}\fq=\phi_{1}^{-1}(\phi_{1}^{k+1}-\phi_{2}^{k+1})$. \hfill\(\Box\)

\section{The two deformed gamma functions}\label{sec:gamma}

\subsection{The deformed gamma function of first kind}\label{subsec:gamma1}

\begin{definition}\label{def:gamma1}
Assume $\Hphi$ and put $a_{\mathcal{R}(p,q)}:=\phi_{1}/(\phi_{1}-\phi_{2})$. For $\Re z>0$ the first kind $\mathcal{R}(p,q)-$gamma function is
\begin{eqnarray}\label{eq:gamma1-def}
	\GR(z)&:=&\phi_{1}^{\frac{z(z-1)}{2}}\int_{0}^{a_{\mathcal{R}(p,q)}}t^{z-1}\,\ER(-\phi_{2}t)\,\dR t .
\end{eqnarray}
\end{definition}

The nodes of the integral \eqref{eq:gamma1-def} are $t_{k}=\fq^{k}/(\phi_{1}-\phi_{2})$, $k\geq0$, and, by Lemma \ref{lem:exponentials}, $\ER(-\phi_{2}t_{k})=(\fq^{k+1};\fq)_{\infty}$. Since $(\fq^{k+1};\fq)_{\infty}=0$ for every negative integer $k$, the integral \eqref{eq:gamma1-def} coincides with the bilateral integral
\begin{eqnarray}
	\int_{0}^{\infty;a_{\mathcal{R}(p,q)}}t^{z-1}\ER(-\phi_{2}t)\,\dR t;\nn
\end{eqnarray}
in this sense \eqref{eq:gamma1-def} is the deformed analogue of the Euler integral on $(0,\infty)$.

\begin{theorem}\label{thm:gamma1}
Assume $\Hphi$. Then for $\Re z>0$,
\begin{eqnarray}\label{eq:gamma1-closed}
	\GR(z)&=&\phi_{1}^{\frac{(z-1)(z-2)}{2}}\,\Gq(z)
	\;=\;\phi_{1}^{\frac{(z-1)(z-2)}{2}}\,\frac{(\fq;\fq)_{\infty}}{(\fq^{z};\fq)_{\infty}}\,(1-\fq)^{1-z}.
\end{eqnarray}
Consequently $\GR$ extends to a meromorphic function on $\mathbb{C}$ without zeros, whose poles are simple and located at $z=-n+2\pi i m/\log\fq$, $n\in\mathbb{N}\cup\{0\}$, $m\in\mathbb{Z}$, and
\begin{eqnarray}\label{eq:gamma1-functional}
	\GR(z+1)\;=\;\pn{z}\,\GR(z),\qquad \GR(1)\;=\;1,\qquad \GR(n+1)\;=\;\pf{n},\quad n\in\mathbb{N}.
\end{eqnarray}
\end{theorem}
{\it Proof.}
The weights of \eqref{eq:gamma1-def} are $(\phi_{1}-\phi_{2})a_{\mathcal{R}(p,q)}\phi_{2}^{k}/\phi_{1}^{k+1}=\fq^{k}$, so that
\begin{eqnarray}
	\int_{0}^{a_{\mathcal{R}(p,q)}}t^{z-1}\ER(-\phi_{2}t)\,\dR t
	&=&\sum_{k\geq0}\fq^{k}\Big(\frac{\fq^{k}}{\phi_{1}-\phi_{2}}\Big)^{z-1}(\fq^{k+1};\fq)_{\infty}
	\;=\;(\phi_{1}-\phi_{2})^{1-z}(\fq;\fq)_{\infty}\sum_{k\geq0}\frac{\fq^{kz}}{(\fq;\fq)_{k}}.\nn
\end{eqnarray}
By Euler's identity $\sum_{k\geq0}w^{k}/(\fq;\fq)_{k}=1/(w;\fq)_{\infty}$ with $w=\fq^{z}$, $|w|<1$ \cite[(1.3.15)]{GasperRahman}, the integral equals $(\phi_{1}-\phi_{2})^{1-z}(\fq;\fq)_{\infty}/(\fq^{z};\fq)_{\infty}=\phi_{1}^{1-z}\Gq(z)$. Multiplying by $\phi_{1}^{z(z-1)/2}$ and using $z(z-1)/2+1-z=(z-1)(z-2)/2$ gives \eqref{eq:gamma1-closed}. The meromorphic continuation and the location of the poles follow from the zeros of $(\fq^{z};\fq)_{\infty}$, namely $\fq^{z+j}=1$. Finally, $\Gq(z+1)=[z]_{\fq}\Gq(z)$ and $\phi_{1}^{z(z-1)/2-(z-1)(z-2)/2}=\phi_{1}^{z-1}$, so that $\GR(z+1)=\phi_{1}^{z-1}[z]_{\fq}\GR(z)=\pn{z}\GR(z)$ by Lemma \ref{lem:numbers}; the value $\GR(1)=1$ is read off from \eqref{eq:gamma1-closed}. \hfill\(\Box\)

\begin{remark}\label{rem:sadjang}
For the Jagannathan Srinivasa Rao deformation, \eqref{eq:gamma1-closed} reads $\Gamma_{p,q}(z)=p^{(z-1)(z-2)/2}\Gamma_{q/p}(z)$, which interpolates the $(p,q)-$factorials. The $(p,q)-$gamma function of \cite{ref59} is written as the quotient $(p\ominus q)^{\infty}_{p,q}(p^{z}\ominus q^{z})^{-\infty}_{p,q}(p-q)^{1-z}$ of infinite $(p,q)-$products. At level $N$, the partial quotient $(\phi_{1}\ominus\phi_{2})^{N}_{\mathcal{R}(p,q)}(\phi_{1}-\phi_{2})^{1-z}/(\phi_{1}^{z}\ominus\phi_{2}^{z})^{N}_{\mathcal{R}(p,q)}$ equals $\phi_{1}^{(N+1)(1-z)}(\fq;\fq)_{N}(1-\fq)^{1-z}/(\fq^{z};\fq)_{N}$, which converges as written when $\phi_{1}=1$; for $\phi_{1}\neq1$ the factor $\phi_{1}^{(N+1)(1-z)}$ has to be removed, and \eqref{eq:gamma1-closed} is the resulting normalisation, the only one of the form $\phi_{1}^{\textnormal{quadratic}}\Gq$ which satisfies \eqref{eq:gamma1-functional}.
\end{remark}

\begin{theorem}[Bohr Mollerup type characterisation]\label{thm:BM}
Assume $\Hphi$ and let $f:(0,\infty)\to(0,\infty)$ satisfy
\begin{enumerate}
\item[\textnormal{(i)}] $f(1)=1$ and $f(x+1)=\pn{x}\,f(x)$ for $x>0$;
\item[\textnormal{(ii)}] the function $x\mapsto\phi_{1}^{-(x-1)(x-2)/2}f(x)$ is logarithmically convex on $(0,\infty)$.
\end{enumerate}
Then $f=\GR$ on $(0,\infty)$. Conversely, $\GR$ satisfies \textnormal{(i)} and \textnormal{(ii)}.
\end{theorem}
{\it Proof.}
Put $g(x):=\phi_{1}^{-(x-1)(x-2)/2}f(x)$. Then $g(1)=1$ and, since $(x-1)(x-2)/2-x(x-1)/2=1-x$,
$g(x+1)=\phi_{1}^{1-x}\pn{x}g(x)=[x]_{\fq}\,g(x)$. Being positive and logarithmically convex, $g$ coincides with $\Gq$ by Askey's $\fq-$analogue of the Bohr Mollerup theorem \cite[Theorem 3.1]{Askey1978}, hence $f=\phi_{1}^{(x-1)(x-2)/2}\Gq=\GR$ by Theorem \ref{thm:gamma1}. Conversely, (i) is \eqref{eq:gamma1-functional}, and $\phi_{1}^{-(x-1)(x-2)/2}\GR(x)=\Gq(x)$ is logarithmically convex by \eqref{eq:logGq} below. \hfill\(\Box\)

For $x>0$, logarithmic differentiation of \eqref{eq:Jackson-gamma} gives
\begin{eqnarray}\label{eq:logGq}
	\frac{d^{2}}{dx^{2}}\log\Gq(x)&=&(\log\fq)^{2}\sum_{k\geq0}\frac{\fq^{x+k}}{(1-\fq^{x+k})^{2}}\;>\;0 ,
\end{eqnarray}
and the right hand side is a strictly decreasing function of $x$ which tends to $+\infty$ as $x\to0^{+}$ and to $0$ as $x\to+\infty$, because $u\mapsto u/(1-u)^{2}$ is increasing on $(0,1)$.

\begin{theorem}\label{thm:logconvex-gamma}
Assume $\Hphi$. For $x>0$,
\begin{eqnarray}\label{eq:logGR}
	\frac{d^{2}}{dx^{2}}\log\GR(x)&=&\log\phi_{1}+(\log\fq)^{2}\sum_{k\geq0}\frac{\fq^{x+k}}{(1-\fq^{x+k})^{2}} .
\end{eqnarray}
Consequently:
\begin{enumerate}
\item[\textnormal{(i)}] if $\phi_{1}\geq1$, then $\GR$ is strictly logarithmically convex on $(0,\infty)$;
\item[\textnormal{(ii)}] if $\phi_{1}<1$, there exists a unique $x_{*}=x_{*}(\phi_{1},\fq)>0$ such that $\GR$ is strictly logarithmically convex on $(0,x_{*}]$ and strictly logarithmically concave on $[x_{*},\infty)$; in particular $\GR$ is not logarithmically convex on $(0,\infty)$;
\item[\textnormal{(iii)}] in all cases $\frac{d^{2}}{dx^{2}}\log\GR(x)\geq\log\phi_{1}$, with $\frac{d^{2}}{dx^{2}}\log\GR(x)\to\log\phi_{1}$ as $x\to+\infty$.
\end{enumerate}
\end{theorem}
{\it Proof.}
By Theorem \ref{thm:gamma1}, $\log\GR(x)=\tfrac{(x-1)(x-2)}{2}\log\phi_{1}+\log\Gq(x)$, whose second derivative is \eqref{eq:logGR}. Assertions (i), (ii) and (iii) follow from the monotonicity and limits of \eqref{eq:logGq} recalled above: if $\phi_{1}<1$, the function \eqref{eq:logGR} decreases strictly from $+\infty$ to $\log\phi_{1}<0$ and has therefore a unique zero $x_{*}$. \hfill\(\Box\)

\begin{remark}\label{rem:logconvex}
Theorem \ref{thm:logconvex-gamma} shows that the logarithmic convexity of the deformed gamma function is not inherited from the undeformed theory when the dominant deformation function satisfies $\phi_{1}<1$: the quadratic valuation factor $\phi_{1}^{(x-1)(x-2)/2}$ is then logarithmically concave and eventually dominates. For the Jagannathan Srinivasa Rao deformation with $(p,q)=(0.9,0.5)$ one finds $x_{*}=3.6722\ldots$ (Section \ref{sec:numerics}). The correct convexity statement is the normalised one of Theorem \ref{thm:BM}, which holds for every admissible pair $(\phi_{1},\phi_{2})$.
\end{remark}

\subsection{The deformed gamma function of second kind}\label{subsec:gamma2}

\begin{definition}\label{def:gamma2}
Assume $\Hphi$ and let $a_{\mathcal{R}(p,q)}=\phi_{1}/(\phi_{1}-\phi_{2})$. For $\Re z>0$ the second kind $\mathcal{R}(p,q)-$gamma function is
\begin{eqnarray}\label{eq:gamma2-def}
	\gR(z)&:=&\phi_{2}^{\frac{z(z-1)}{2}}\int_{0}^{\infty;a_{\mathcal{R}(p,q)}}t^{z-1}\,\eR(-\phi_{1}t)\,\dR t ,
\end{eqnarray}
where $\eR$ is continued to the negative half line by \eqref{eq:exp-products}.
\end{definition}

At the nodes $t_{k}=\fq^{k}/(\phi_{1}-\phi_{2})$, $k\in\mathbb{Z}$, one has $\eR(-\phi_{1}t_{k})=1/(-\fq^{k};\fq)_{\infty}$. For $k=-m\to-\infty$, $(-\fq^{-m};\fq)_{\infty}=(-1;\fq)_{\infty}\prod_{i=1}^{m}(1+\fq^{-i})\geq2\fq^{-m(m+1)/2}$, so that the bilateral series converges for every $z$ with $\Re z>0$, the condition $\Re z>0$ being needed for $k\to+\infty$.

We introduce the theta quotient
\begin{eqnarray}\label{eq:Omega}
	\Omq(z)&:=&\fq^{\frac{z(z-1)}{2}}\,\Kq(z),\qquad
	\Kq(z)\;:=\;\frac{(-\fq^{z},-\fq^{1-z};\fq)_{\infty}}{(-1,-\fq;\fq)_{\infty}},\qquad z\in\mathbb{C}.
\end{eqnarray}

\begin{proposition}\label{prop:Omega}
The functions $\Kq$ and $\Omq$ are entire, and:
\begin{enumerate}
\item[\textnormal{(i)}] $\Kq(z+1)=\fq^{-z}\Kq(z)$, $\Omq(z+1)=\Omq(z)$, $\Omq(1-z)=\Omq(z)$ and $\Omq(n)=1$ for $n\in\mathbb{Z}$;
\item[\textnormal{(ii)}] $\Omq(z)=\vartheta_{\fq}(z)/\vartheta_{\fq}(0)$ with $\vartheta_{\fq}(z):=\sum_{n\in\mathbb{Z}}\fq^{(n+z)(n+z-1)/2}$;
\item[\textnormal{(iii)}] with $\lambda:=\log(1/\fq)$ and $\varepsilon:=e^{-2\pi^{2}/\lambda}$,
\begin{eqnarray}\label{eq:Omega-Poisson}
	\Omq(z)&=&\frac{\sum_{k\in\mathbb{Z}}(-1)^{k}\varepsilon^{k^{2}}\cos(2\pi kz)}{\sum_{k\in\mathbb{Z}}(-1)^{k}\varepsilon^{k^{2}}},
\end{eqnarray}
so that, for real $z$, $\Omq(z)>0$ and $\Omq(z)-1=4\sin^{2}(\pi z)\,\varepsilon+O(\varepsilon^{2})$ uniformly in $z$; in particular $\Omq$ is not constant.
\end{enumerate}
\end{proposition}
{\it Proof.}
(i) From $(-\fq^{z+1};\fq)_{\infty}=(-\fq^{z};\fq)_{\infty}/(1+\fq^{z})$ and $(-\fq^{-z};\fq)_{\infty}=(1+\fq^{-z})(-\fq^{1-z};\fq)_{\infty}$ we get $\Kq(z+1)/\Kq(z)=(1+\fq^{-z})/(1+\fq^{z})=\fq^{-z}$, and $\fq^{(z+1)z/2}=\fq^{z(z-1)/2}\fq^{z}$ gives the periodicity of $\Omq$. The symmetry follows from the symmetry of $\Kq$ under $z\mapsto1-z$, and $\Omq(1)=\Kq(1)=1$. (ii) By the Jacobi triple product identity \cite[(1.6.1)]{GasperRahman}, $(\fq,-w,-\fq/w;\fq)_{\infty}=\sum_{n\in\mathbb{Z}}\fq^{n(n-1)/2}w^{n}$. With $w=\fq^{z}$ this gives $\fq^{z(z-1)/2}(\fq,-\fq^{z},-\fq^{1-z};\fq)_{\infty}=\sum_{n}\fq^{(n+z)(n+z-1)/2}=\vartheta_{\fq}(z)$, and the normalisation follows from $\Omq(0)=\Omq(1)=1$. (iii) Writing $(n+z)(n+z-1)/2=(n+u)^{2}/2-1/8$ with $u=z-1/2$, we get $\vartheta_{\fq}(z)=e^{\lambda/8}\sum_{n}e^{-\lambda(n+u)^{2}/2}$, and the Poisson summation formula $\sum_{n}e^{-\lambda(n+u)^{2}/2}=\sqrt{2\pi/\lambda}\sum_{k}e^{-2\pi^{2}k^{2}/\lambda}e^{2\pi iku}$ yields \eqref{eq:Omega-Poisson} after grouping $\pm k$ and using $e^{2\pi ik(z-1/2)}=(-1)^{k}e^{2\pi ikz}$. Expanding numerator and denominator in powers of $\varepsilon$ gives $\Omq(z)=1+2\varepsilon(1-\cos2\pi z)+O(\varepsilon^{2})$. Positivity follows from (ii), and $\Omq(1/2)-1=4\varepsilon+O(\varepsilon^{2})\neq0$ for $\varepsilon$ small, while $\Omq(1/2)\neq1$ for all $\fq\in(0,1)$ because the numerator minus the denominator in \eqref{eq:Omega-Poisson} equals $\sum_{k\ \mathrm{odd}}4\varepsilon^{k^{2}}>0$ at $z=1/2$. \hfill\(\Box\)

\begin{theorem}\label{thm:gamma2}
Assume $\Hphi$. For $\Re z>0$,
\begin{eqnarray}\label{eq:gamma2-closed}
	\gR(z)&=&\Omq(z)\,\GR(z).
\end{eqnarray}
Consequently $\gR$ extends meromorphically to $\mathbb{C}$, $\gR(z+1)=\pn{z}\gR(z)$, $\gR(n+1)=\pf{n}$ for $n\in\mathbb{N}\cup\{0\}$, but $\gR\neq\GR$, and the function $x\mapsto\phi_{1}^{-(x-1)(x-2)/2}\gR(x)$ is not logarithmically convex on $(0,\infty)$.
\end{theorem}
{\it Proof.}
The weights of \eqref{eq:gamma2-def} are $\fq^{k}$ as in the proof of Theorem \ref{thm:gamma1}, hence
\begin{eqnarray}
	\int_{0}^{\infty;a_{\mathcal{R}(p,q)}}t^{z-1}\eR(-\phi_{1}t)\,\dR t
	&=&(\phi_{1}-\phi_{2})^{1-z}\sum_{k\in\mathbb{Z}}\frac{\fq^{kz}}{(-\fq^{k};\fq)_{\infty}}
	\;=\;\frac{(\phi_{1}-\phi_{2})^{1-z}}{(-1;\fq)_{\infty}}\sum_{k\in\mathbb{Z}}(-1;\fq)_{k}\,\fq^{kz},\nn
\end{eqnarray}
because $(-\fq^{k};\fq)_{\infty}=(-1;\fq)_{\infty}/(-1;\fq)_{k}$ for all $k\in\mathbb{Z}$. Ramanujan's ${}_{1}\psi_{1}$ summation \cite[(5.2.1)]{GasperRahman}, \cite{AndrewsAskey1978},
\begin{eqnarray}\label{eq:psi11}
	\sum_{k\in\mathbb{Z}}\frac{(a;\fq)_{k}}{(b;\fq)_{k}}\,w^{k}&=&\frac{(\fq,b/a,aw,\fq/(aw);\fq)_{\infty}}{(b,\fq/a,w,b/(aw);\fq)_{\infty}},\qquad |b/a|<|w|<1,
\end{eqnarray}
with $a=-1$, $b=0$ and $w=\fq^{z}$, gives $\sum_{k}(-1;\fq)_{k}\fq^{kz}=(\fq,-\fq^{z},-\fq^{1-z};\fq)_{\infty}/(-\fq,\fq^{z};\fq)_{\infty}$. Therefore the integral equals $(\phi_{1}-\phi_{2})^{1-z}(\fq;\fq)_{\infty}\Kq(z)/(\fq^{z};\fq)_{\infty}=\phi_{1}^{1-z}\Kq(z)\Gq(z)$, and
\begin{eqnarray}
	\gR(z)&=&\phi_{2}^{\frac{z(z-1)}{2}}\phi_{1}^{1-z}\Kq(z)\Gq(z)
	\;=\;\fq^{\frac{z(z-1)}{2}}\Kq(z)\,\phi_{1}^{\frac{z(z-1)}{2}+1-z}\Gq(z)
	\;=\;\Omq(z)\,\GR(z).\nn
\end{eqnarray}
The functional equation and the interpolation property follow from Theorem \ref{thm:gamma1} and Proposition \ref{prop:Omega}(i). Since $\Omq$ is not constant, $\gR\neq\GR$, and Theorem \ref{thm:BM} shows that $\phi_{1}^{-(x-1)(x-2)/2}\gR(x)=\Omq(x)\Gq(x)$ cannot be logarithmically convex, since otherwise it would coincide with $\Gq$. \hfill\(\Box\)

\begin{remark}\label{rem:gamma2}
The two deformed gamma functions coincide on the integers and satisfy the same functional equation, and their quotient is the one periodic function $\Omq$, which is extremely close to $1$: for $\fq=0.5$ one has $\Omq(1/2)-1=1.7\times10^{-12}$. This is the deformed counterpart of the classical fact that the Jackson integral and the bilateral Ramanujan type integral define two $q-$gamma functions which differ by a quasi constant \cite{Askey1978, GasperRahman}. The deformation enters \eqref{eq:gamma2-def} through the prefactor $\phi_{2}^{z(z-1)/2}$ and through the deformed exponential $\eR$ evaluated at $-\phi_{1}t$, which are exchanged with respect to \eqref{eq:gamma1-def}.
\end{remark}

\subsection{Multiplication formula, reflection and Wendel type estimates}\label{subsec:gamma-further}

For $n\in\mathbb{N}$, the pair $(\phi_{1}^{n},\phi_{2}^{n})$ satisfies the hypothesis $\Hphi$ with ratio $\fq^{n}$, and its deformed numbers are $(\phi_{1}^{nm}-\phi_{2}^{nm})/(\phi_{1}^{n}-\phi_{2}^{n})=\pn{nm}/\pn{n}$. We denote by
\begin{eqnarray}\label{eq:gamma-n}
	\Gamma_{\mathcal{R}(p,q)}^{[n]}(w)&:=&\phi_{1}^{\frac{n(w-1)(w-2)}{2}}\,\Gamma_{\fq^{n}}(w)
\end{eqnarray}
the first kind deformed gamma function attached to $(\phi_{1}^{n},\phi_{2}^{n})$ by Theorem \ref{thm:gamma1}; for $n=1$ it is $\GR$.

\begin{theorem}[Deformed Gauss multiplication formula]\label{thm:multiplication}
Assume $\Hphi$ and let $n\in\mathbb{N}$, $n\geq2$. Then, as an identity between meromorphic functions of $z$,
\begin{eqnarray}\label{eq:multiplication}
	\GR(nz)\prod_{k=1}^{n-1}\Gamma_{\mathcal{R}(p,q)}^{[n]}\Big(\frac{k}{n}\Big)
	&=&\pn{n}^{\,nz-1}\prod_{k=0}^{n-1}\Gamma_{\mathcal{R}(p,q)}^{[n]}\Big(z+\frac{k}{n}\Big).
\end{eqnarray}
In particular, for $n=2$ one obtains the deformed Legendre duplication formula
\begin{eqnarray}\label{eq:duplication}
	\GR(2z)\,\Gamma_{\mathcal{R}(p,q)}^{[2]}\big(\tfrac12\big)&=&(\phi_{1}+\phi_{2})^{2z-1}\,\Gamma_{\mathcal{R}(p,q)}^{[2]}(z)\,\Gamma_{\mathcal{R}(p,q)}^{[2]}\big(z+\tfrac12\big).
\end{eqnarray}
\end{theorem}
{\it Proof.}
The $\fq-$analogue of the Gauss multiplication formula \cite[Section 1.10]{GasperRahman} reads
$\Gq(nz)\prod_{k=1}^{n-1}\Gamma_{\fq^{n}}(k/n)=[n]_{\fq}^{\,nz-1}\prod_{k=0}^{n-1}\Gamma_{\fq^{n}}(z+k/n)$.
By Theorem \ref{thm:gamma1}, \eqref{eq:gamma-n} and $\pn{n}=\phi_{1}^{n-1}[n]_{\fq}$, the identity \eqref{eq:multiplication} is equivalent to the exponent identity $A_{n}(z)=(n-1)(nz-1)+B_{n}(z)$, where
\begin{eqnarray}
	A_{n}(z)&=&\frac{(nz-1)(nz-2)}{2}+\frac{n}{2}\sum_{k=1}^{n-1}\Big(\frac{k}{n}-1\Big)\Big(\frac{k}{n}-2\Big),\qquad
	B_{n}(z)\;=\;\frac{n}{2}\sum_{k=0}^{n-1}\Big(z+\frac{k}{n}-1\Big)\Big(z+\frac{k}{n}-2\Big).\nn
\end{eqnarray}
Put $C_{n}=\sum_{k=0}^{n-1}(k/n-1)(k/n-2)$. Since $\sum_{k=0}^{n-1}k/n=(n-1)/2$, one finds $B_{n}(z)=\tfrac{n}{2}\big[nz^{2}-(2n+1)z+C_{n}\big]$ and $A_{n}(z)=\tfrac{1}{2}(n^{2}z^{2}-3nz+2)+\tfrac{n}{2}(C_{n}-2)$, whence
$A_{n}(z)-B_{n}(z)=\tfrac12\big[(2n^{2}-2n)z+2-2n\big]=(n-1)(nz-1)$. For $n=2$, $\pn{2}=\phi_{1}+\phi_{2}$. \hfill\(\Box\)

The duplication formula \eqref{eq:duplication} shows that the base $(\phi_{1}+\phi_{2})$ of the classical factor $2^{2z-1}$ is the deformed number $\pn{2}$, while the gamma functions on the right hand side are those of the squared deformation $(\phi_{1}^{2},\phi_{2}^{2})$; no extra valuation factor appears, the quadratic exponents of the three deformed gamma functions compensating exactly.

\begin{proposition}[Deformed reflection formula]\label{prop:reflection}
Assume $\Hphi$ and let
\begin{eqnarray}\label{eq:theta}
	\theta_{\fq}(w)&:=&(w,\fq/w,\fq;\fq)_{\infty}\;=\;\sum_{n\in\mathbb{Z}}(-1)^{n}\fq^{n(n-1)/2}w^{n}.
\end{eqnarray}
For $0<x<1$,
\begin{eqnarray}\label{eq:reflection}
	\GR(x)\,\GR(1-x)&=&\phi_{1}^{\,x^{2}-x+1}\,\frac{(1-\fq)\,(\fq;\fq)_{\infty}^{3}}{\theta_{\fq}(\fq^{x})},\qquad
	\BR(x,1-x)\;=\;\phi_{1}^{\frac{x^{2}-x+2}{2}}\,\frac{(1-\fq)\,(\fq;\fq)_{\infty}^{3}}{\theta_{\fq}(\fq^{x})}.
\end{eqnarray}
As $\phi_{1}\to1$ and $\fq\to1^{-}$, both right hand sides tend to $\pi/\sin(\pi x)$.
\end{proposition}
{\it Proof.}
By \eqref{eq:Jackson-gamma}, $\Gq(x)\Gq(1-x)=(\fq;\fq)_{\infty}^{2}(1-\fq)/((\fq^{x};\fq)_{\infty}(\fq^{1-x};\fq)_{\infty})$, and the series expression of $\theta_{\fq}$ is the Jacobi triple product identity \cite[Section 1.6]{GasperRahman}. The valuation exponents are $\tfrac{(x-1)(x-2)}{2}+\tfrac{x(x+1)}{2}=x^{2}-x+1$ for the gamma product and, by Theorem \ref{thm:beta1}, $1-x+\tfrac{(-x)(-x-1)}{2}=\tfrac{x^{2}-x+2}{2}$ for the beta function, since $\GR(1)=1$. The limit follows from $\Gq\to\Gamma$ \cite{Askey1978} and from the classical reflection formula. \hfill\(\Box\)

The next result is the deformed counterpart of the two sided estimate of Wendel \cite{Wendel1948}; it is obtained, as in the classical case, from the logarithmic convexity of the normalised deformed gamma function.

\begin{theorem}[Wendel type estimate]\label{thm:wendel}
Assume $\Hphi$. For $x>0$ and $0\leq s\leq1$,
\begin{eqnarray}\label{eq:wendel}
	\phi_{1}^{\frac{s(1-s)}{2}}\,\pn{x}\,\pn{x+s}^{\,s-1}
	\;\leq\;\frac{\GR(x+s)}{\GR(x)}\;\leq\;
	\phi_{1}^{-\frac{s(1-s)}{2}}\,\pn{x}^{\,s}.
\end{eqnarray}
Consequently, for every real $y>0$,
\begin{eqnarray}\label{eq:wendel-limit}
	\lim_{x\to+\infty}\frac{\GR(x+y)}{\phi_{1}^{\frac{y(y-1)}{2}}\,\pn{x}^{\,y}\,\GR(x)}&=&1 .
\end{eqnarray}
\end{theorem}
{\it Proof.}
By \eqref{eq:logGq}, $\Gq$ is logarithmically convex on $(0,\infty)$. Since $x+s=(1-s)x+s(x+1)$ and $x+1=s(x+s)+(1-s)(x+s+1)$, the convexity of $\log\Gq$ and the functional equation $\Gq(u+1)=[u]_{\fq}\Gq(u)$ give
\begin{eqnarray}
	\Gq(x+s)&\leq&\Gq(x)^{1-s}\Gq(x+1)^{s}\;=\;[x]_{\fq}^{\,s}\,\Gq(x),\nn\\
	\Gq(x+1)&\leq&\Gq(x+s)^{s}\,\Gq(x+s+1)^{1-s}\;=\;[x+s]_{\fq}^{\,1-s}\,\Gq(x+s).\nn
\end{eqnarray}
Hence $[x]_{\fq}[x+s]_{\fq}^{s-1}\leq\Gq(x+s)/\Gq(x)\leq[x]_{\fq}^{s}$. By Theorem \ref{thm:gamma1}, $\GR(x+s)/\GR(x)=\phi_{1}^{s(2x+s-3)/2}\Gq(x+s)/\Gq(x)$, and by Lemma \ref{lem:numbers} $[x]_{\fq}=\phi_{1}^{1-x}\pn{x}$. The upper exponent is $s(2x+s-3)/2+s(1-x)=-s(1-s)/2$, and the lower one is $s(2x+s-3)/2+(1-x)-(1-x-s)(1-s)=s(1-s)/2$, which gives \eqref{eq:wendel}. As $x\to+\infty$, $\pn{x+s}/\pn{x}\to\phi_{1}^{s}$, so that the ratio of the two bounds tends to $\phi_{1}^{s(1-s)}\phi_{1}^{s(s-1)}=1$, which is \eqref{eq:wendel-limit} for $y=s\in[0,1]$, because $\phi_{1}^{-s(1-s)/2}=\phi_{1}^{s(s-1)/2}$. For $y=n+s$ with $n\in\mathbb{N}$ and $s\in[0,1)$, $\GR(x+y)=\Pch{\pn{x+s}}{n}\GR(x+s)$, and $\Pch{\pn{x+s}}{n}\sim\phi_{1}^{ns+n(n-1)/2}\pn{x}^{n}$, since $\pn{x+a}/\pn{x}\to\phi_{1}^{a}$; the exponents add up to $ns+n(n-1)/2+s(s-1)/2=y(y-1)/2$. \hfill\(\Box\)

The limit \eqref{eq:wendel-limit} replaces the classical asymptotic relation $\Gamma(x+y)\sim x^{y}\Gamma(x)$: the power $x^{y}$ becomes the power $\pn{x}^{y}$ of the deformed number, corrected by the valuation factor $\phi_{1}^{y(y-1)/2}$, which is the value at $y$ of the exponent of the deformed factorial $\pf{n}=\phi_{1}^{n(n-1)/2}[n]_{\fq}!$.

\section{The first kind noncommutative deformed beta function}\label{sec:beta1}

\subsection{Definition and closed form}\label{subsec:beta1-def}

\begin{definition}\label{def:beta1}
Assume $\Hphi$ and let $\Omega_{+}:=\{(x,y)\in\mathbb{C}^{2}:\Re x>0,\ \Re y>0\}$. The first kind $\mathcal{R}(p,q)-$beta function is
\begin{eqnarray}\label{eq:beta1-def}
	\BR(x,y)&:=&\int_{0}^{1}t^{x-1}\,(1\ominus\phi_{2}t)^{y-1}_{\mathcal{R}(p,q)}\,\dR t,\qquad (x,y)\in\Omega_{+},
\end{eqnarray}
where the deformed power is given by \eqref{eq:powers-complex}, that is
$(1\ominus\phi_{2}t)^{y-1}_{\mathcal{R}(p,q)}=\phi_{1}^{(y-1)(y-2)/2}(\phi_{2}t;\fq)_{\infty}/(\phi_{2}\fq^{y-1}t;\fq)_{\infty}$.
\end{definition}

For $y=n+1$, $n\in\mathbb{N}\cup\{0\}$, the integrand is the polynomial $t^{x-1}\prod_{k=0}^{n-1}(\phi_{1}^{k}-\phi_{2}^{k+1}t)$, and for the Jagannathan Srinivasa Rao deformation and integer arguments \eqref{eq:beta1-def} is the $(p,q)-$beta integral of \cite{ref55}.

\begin{lemma}\label{lem:beta1-series}
Assume $\Hphi$. For $(x,y)\in\Omega_{+}$,
\begin{eqnarray}\label{eq:beta1-series}
	\BR(x,y)&=&(1-\fq)\,\phi_{1}^{1-x+\frac{(y-1)(y-2)}{2}}\sum_{j\geq0}\fq^{jx}\,\frac{(\fq^{j+1};\fq)_{\infty}}{(\fq^{j+y};\fq)_{\infty}},
\end{eqnarray}
the series converges absolutely and uniformly on compact subsets of $\Omega_{+}$, $\BR$ is holomorphic on $\Omega_{+}$, and $\BR(x,y)>0$ for real $x,y>0$.
\end{lemma}
{\it Proof.}
The nodes of \eqref{eq:beta1-def} are $t_{j}=\fq^{j}/\phi_{1}$ with weights $(1-\fq)\fq^{j}$, and $\phi_{2}t_{j}=\fq^{j+1}$, which gives \eqref{eq:beta1-series}. For $\Re y\geq\delta>0$ one has $|1-\fq^{j+y+k}|\geq1-\fq^{\delta+k}$, hence $|(\fq^{j+y};\fq)_{\infty}|\geq(\fq^{\delta};\fq)_{\infty}>0$, while $0<(\fq^{j+1};\fq)_{\infty}\leq1$; the $j-$th term is therefore bounded by $\fq^{j\Re x}/(\fq^{\delta};\fq)_{\infty}$, which gives the uniform convergence on $\{\Re x\geq\delta,\Re y\geq\delta\}$. Each term being holomorphic, $\BR$ is holomorphic on $\Omega_{+}$ by the Weierstrass theorem. Positivity follows from Lemma \ref{lem:powers}(v). \hfill\(\Box\)

The holomorphy of $\BR$ on $\Omega_{+}$ thus holds for every admissible pair $\phi_{1}>\phi_{2}>0$; no normalisation of the form $\phi_{1}\leq1$ is needed, since only the ratio $\fq$ enters the convergence of \eqref{eq:beta1-series}.

\begin{theorem}\label{thm:beta1}
Assume $\Hphi$. For $(x,y)\in\Omega_{+}$,
\begin{eqnarray}\label{eq:beta1-gamma}
	\BR(x,y)&=&\phi_{1}^{\frac{(y-1)(2x+y-2)}{2}}\,\frac{\GR(x)\,\GR(y)}{\GR(x+y)}
	\;=\;\phi_{1}^{1-x+\frac{(y-1)(y-2)}{2}}\,\Bq(x,y).
\end{eqnarray}
In particular $\BR$ extends to a meromorphic function on $\mathbb{C}^{2}$.
\end{theorem}
{\it Proof.}
Since $(\fq^{j+1};\fq)_{\infty}=(\fq;\fq)_{\infty}/(\fq;\fq)_{j}$ and $(\fq^{j+y};\fq)_{\infty}=(\fq^{y};\fq)_{\infty}/(\fq^{y};\fq)_{j}$, Lemma \ref{lem:binomial} with $b=\fq^{y}$ and $w=\fq^{x}$ gives
\begin{eqnarray}
	\sum_{j\geq0}\fq^{jx}\,\frac{(\fq^{j+1};\fq)_{\infty}}{(\fq^{j+y};\fq)_{\infty}}
	&=&\frac{(\fq;\fq)_{\infty}}{(\fq^{y};\fq)_{\infty}}\sum_{j\geq0}\frac{(\fq^{y};\fq)_{j}}{(\fq;\fq)_{j}}\,\fq^{jx}
	\;=\;\frac{(\fq;\fq)_{\infty}\,(\fq^{x+y};\fq)_{\infty}}{(\fq^{x};\fq)_{\infty}\,(\fq^{y};\fq)_{\infty}} .\nn
\end{eqnarray}
Multiplying by $1-\fq$ and comparing with \eqref{eq:Jackson-gamma} gives the second equality in \eqref{eq:beta1-gamma}. By Theorem \ref{thm:gamma1},
\begin{eqnarray}
	\frac{\GR(x)\GR(y)}{\GR(x+y)}&=&\phi_{1}^{\frac{(x-1)(x-2)+(y-1)(y-2)-(x+y-1)(x+y-2)}{2}}\,\Bq(x,y)\;=\;\phi_{1}^{1-xy}\,\Bq(x,y),\nn
\end{eqnarray}
and the first equality follows from the exponent identity $1-x+\tfrac{(y-1)(y-2)}{2}-(1-xy)=\tfrac{(y-1)(2x+y-2)}{2}$, verified in Appendix \ref{app:exponents}. The meromorphic continuation is inherited from $\GR$. \hfill\(\Box\)

\begin{corollary}\label{cor:special-values}
Assume $\Hphi$. For $x,y\in\mathbb{C}$ with positive real parts and $m,n\in\mathbb{N}$,
\begin{eqnarray}
	\BR(x,1)&=&\frac{1}{\pn{x}},\qquad
	\BR(1,y)\;=\;\frac{\phi_{1}^{\frac{y(y-1)}{2}}}{\pn{y}},\nn\\
	\BR(m,n)&=&\phi_{1}^{\frac{(n-1)(2m+n-2)}{2}}\,\frac{\pf{m-1}\,\pf{n-1}}{\pf{m+n-1}},\label{eq:special1}\\
	\BR\big(\tfrac12,\tfrac12\big)&=&\phi_{1}^{\frac{7}{8}}\,\Gq\big(\tfrac12\big)^{2}.\label{eq:special2}
\end{eqnarray}
\end{corollary}
{\it Proof.}
The first value is also Proposition \ref{prop:integral}(ii). The others follow from \eqref{eq:beta1-gamma} with $\GR(1)=1$, $\GR(y+1)=\pn{y}\GR(y)$, $\GR(n)=\pf{n-1}$, and, for \eqref{eq:special2}, $(y-1)(2x+y-2)/2=1/8$ at $x=y=1/2$ together with $\GR(1/2)=\phi_{1}^{3/8}\Gq(1/2)$. \hfill\(\Box\)

Table \ref{tab:integer-values} lists the values of $\BR(m,n)$ for $1\leq m,n\leq3$ obtained from \eqref{eq:special1}. The asymmetry between the entries $(m,n)$ and $(n,m)$ is the exchange defect of Theorem \ref{thm:exchange1} below; for instance $\BR(1,2)/\BR(2,1)=\phi_{1}$ and $\BR(2,3)/\BR(3,2)=\phi_{1}^{2}$.

\begin{table}[ht]
\centering
\caption{Values of the first kind deformed beta function at integer points, $\pn{k}=(\phi_{1}^{k}-\phi_{2}^{k})/(\phi_{1}-\phi_{2})$.}\label{tab:integer-values}
\begin{tabular}{llll}
\hline
$m\backslash n$ & $1$ & $2$ & $3$\\
\hline
$1$ & $1$ & $\dfrac{\phi_{1}}{\pn{2}}$ & $\dfrac{\phi_{1}^{3}}{\pn{3}}$\\[3mm]
$2$ & $\dfrac{1}{\pn{2}}$ & $\dfrac{\phi_{1}^{2}}{\pn{2}\pn{3}}$ & $\dfrac{\phi_{1}^{5}}{\pn{3}\pn{4}}$\\[3mm]
$3$ & $\dfrac{1}{\pn{3}}$ & $\dfrac{\phi_{1}^{3}}{\pn{3}\pn{4}}$ & $\dfrac{\phi_{1}^{7}\pn{2}}{\pn{3}\pn{4}\pn{5}}$\\[2mm]
\hline
\end{tabular}
\end{table}

\begin{corollary}[Deformed partial fraction identity]\label{cor:partial-fractions}
Assume $\Hphi$. For $n\in\mathbb{N}\cup\{0\}$ and $x\in\mathbb{C}\setminus\{-k+2\pi im/\log\fq:\ k\in\mathbb{N}\cup\{0\},\ m\in\mathbb{Z}\}$,
\begin{eqnarray}\label{eq:partial-fractions}
	\sum_{k=0}^{n}(-1)^{k}\bigg[\begin{array}{c} n \\ k\end{array}\bigg]_{\mathcal{R}(p,q)}
	\frac{\phi_{1}^{\binom{n-k}{2}}\phi_{2}^{\binom{k+1}{2}}}{\pn{x+k}}
	&=&\phi_{1}^{\frac{n(2x+n-1)}{2}}\,\frac{\pf{n}}{\prod_{j=0}^{n}\pn{x+j}} .
\end{eqnarray}
\end{corollary}
{\it Proof.}
For $\Re x>0$, integrate the deformed Gauss formula of Lemma \ref{lem:powers}(iii) with $b=\phi_{2}t$ termwise against $t^{x-1}$, using Proposition \ref{prop:integral}(ii) and $\phi_{2}^{\binom{k}{2}}\phi_{2}^{k}=\phi_{2}^{\binom{k+1}{2}}$: the left hand side of \eqref{eq:partial-fractions} is $\BR(x,n+1)$, which equals the right hand side by Theorem \ref{thm:beta1}, since $(y-1)(2x+y-2)/2=n(2x+n-1)/2$ for $y=n+1$ and $\GR(x)/\GR(x+n+1)=1/\prod_{j=0}^{n}\pn{x+j}$. Both sides are meromorphic in $x$, and the identity extends to the stated set. \hfill\(\Box\)

\subsection{Noncommutativity and exchange rule}\label{subsec:beta1-exchange}

\begin{theorem}\label{thm:exchange1}
Assume $\Hphi$. For $(x,y)\in\Omega_{+}$,
\begin{eqnarray}\label{eq:exchange1}
	\BR(x,y)&=&\phi_{1}^{\frac{(y-x)(x+y-1)}{2}}\,\BR(y,x).
\end{eqnarray}
In particular, for real $x,y>0$, $\BR(x,y)=\BR(y,x)$ if and only if $\phi_{1}=1$, or $x=y$, or $x+y=1$.
\end{theorem}
{\it Proof.}
The gamma quotient in \eqref{eq:beta1-gamma} is symmetric, so that $\BR(x,y)/\BR(y,x)=\phi_{1}^{E(x,y)-E(y,x)}$ with $E(x,y)=(y-1)(2x+y-2)/2$, and $E(x,y)-E(y,x)=(y-x)(x+y-1)/2$ (Appendix \ref{app:exponents}). For real arguments, $\phi_{1}^{w}=1$ with $w$ real holds if and only if $\phi_{1}=1$ or $w=0$. \hfill\(\Box\)

The noncommutativity of $\BR$ is thus entirely carried by the dominant deformation function $\phi_{1}$; it disappears in the Jackson case $\phi_{1}=1$, where $\BR=\Bq$ is the commutative Jackson beta function \cite{GasperRahman}.

\subsection{Integration by parts, Pascal rule and unit shifts}\label{subsec:beta1-shifts}

\begin{lemma}\label{lem:ibp-pascal}
Assume $\Hphi$ and let $(x,y)\in\Omega_{+}$. Then
\begin{eqnarray}
	\BR(x+1,y)&=&\phi_{1}^{-x}\,\frac{\pn{x}}{\pn{y}}\,\BR(x,y+1),\label{eq:ibp-beta}\\
	\BR(x,y+1)&=&\phi_{1}^{y-1}\,\BR(x,y)-\phi_{2}^{y}\,\BR(x+1,y).\label{eq:pascal-beta}
\end{eqnarray}
\end{lemma}
{\it Proof.}
For \eqref{eq:ibp-beta} apply the integration by parts formula \eqref{eq:ibp} on $[0,1]$ with
$f(t)=(t/\phi_{1})^{x}$ and $g(t)=-(1\ominus t)^{y}_{\mathcal{R}(p,q)}/\pn{y}$. By Lemma \ref{lem:powers}(iv), $\partial_{\mathcal{R}(p,q)}g(t)=(1\ominus\phi_{2}t)^{y-1}_{\mathcal{R}(p,q)}$, while $f(\phi_{1}t)=t^{x}$ and $\partial_{\mathcal{R}(p,q)}f(t)=\phi_{1}^{-x}\pn{x}t^{x-1}$. The boundary term vanishes: $g(1)=0$ by Lemma \ref{lem:powers}(v), and $(fg)(0^{+})=0$ since $\Re x>0$ and $g$ is bounded near the origin. Since $g(\phi_{2}t)=-(1\ominus\phi_{2}t)^{y}_{\mathcal{R}(p,q)}/\pn{y}$, formula \eqref{eq:ibp} gives
\begin{eqnarray}
	\BR(x+1,y)&=&\int_{0}^{1}f(\phi_{1}t)\,\partial_{\mathcal{R}(p,q)}g(t)\,\dR t
	\;=\;\frac{\pn{x}}{\phi_{1}^{x}\pn{y}}\int_{0}^{1}t^{x-1}(1\ominus\phi_{2}t)^{y}_{\mathcal{R}(p,q)}\,\dR t,\nn
\end{eqnarray}
which is \eqref{eq:ibp-beta}. For \eqref{eq:pascal-beta}, Lemma \ref{lem:powers}(ii) with $b=\phi_{2}t$ and $\alpha=y-1$ gives $(1\ominus\phi_{2}t)^{y}_{\mathcal{R}(p,q)}=(1\ominus\phi_{2}t)^{y-1}_{\mathcal{R}(p,q)}(\phi_{1}^{y-1}-\phi_{2}^{y}t)$, and one integrates against $t^{x-1}$. \hfill\(\Box\)

\begin{theorem}[Unit and double shifts]\label{thm:unit-shifts}
Assume $\Hphi$. For $(x,y)\in\Omega_{+}$,
\begin{eqnarray}
	\BR(x+1,y)&=&\phi_{1}^{\,y-1}\,\frac{\pn{x}}{\pn{x+y}}\,\BR(x,y),\label{eq:shift-x}\\
	\BR(x,y+1)&=&\phi_{1}^{\,x+y-1}\,\frac{\pn{y}}{\pn{x+y}}\,\BR(x,y),\label{eq:shift-y}\\
	\BR(x+1,y+1)&=&\phi_{1}^{\,x+2y-1}\,\frac{\pn{x}\,\pn{y}}{\pn{x+y}\,\pn{x+y+1}}\,\BR(x,y).\label{eq:shift-xy}
\end{eqnarray}
\end{theorem}
{\it Proof.}
We give a proof internal to the deformed calculus, independent of Theorem \ref{thm:beta1}. Inserting \eqref{eq:ibp-beta} into \eqref{eq:pascal-beta},
\begin{eqnarray}
	\BR(x,y+1)\Big(1+\frac{\phi_{2}^{y}\pn{x}}{\phi_{1}^{x}\pn{y}}\Big)&=&\phi_{1}^{y-1}\BR(x,y).\nn
\end{eqnarray}
By the addition formula of Lemma \ref{lem:numbers}(ii), $\phi_{1}^{x}\pn{y}+\phi_{2}^{y}\pn{x}=\pn{x+y}$, so that the bracket equals $\pn{x+y}/(\phi_{1}^{x}\pn{y})$, which is \eqref{eq:shift-y}. Combining \eqref{eq:shift-y} with \eqref{eq:ibp-beta} gives
\begin{eqnarray}
	\BR(x+1,y)&=&\phi_{1}^{-x}\tfrac{\pn{x}}{\pn{y}}\phi_{1}^{x+y-1}\tfrac{\pn{y}}{\pn{x+y}}\BR(x,y),\nn
\end{eqnarray}
which is \eqref{eq:shift-x}. Finally \eqref{eq:shift-xy} follows by applying \eqref{eq:shift-y} at $(x,y)$ and \eqref{eq:shift-x} at $(x,y+1)$, the exponents adding up to $(x+y-1)+y$. \hfill\(\Box\)

Theorem \ref{thm:unit-shifts} also follows from \eqref{eq:beta1-gamma} and the functional equation \eqref{eq:gamma1-functional}; the coincidence of the two derivations is a consistency check of the exponent $E(x,y)$, since $E(x+1,y)-E(x,y)=y-1$ and $E(x,y+1)-E(x,y)=x+y-1$.

\subsection{Finite shifts}\label{subsec:beta1-finite}

For $a\in\mathbb{C}$ and $n\in\mathbb{N}\cup\{0\}$ we use the deformed Pochhammer symbol
\begin{eqnarray}\label{eq:pochhammer}
	\Pch{\pn{a}}{n}&:=&\prod_{j=0}^{n-1}\pn{a+j}\;=\;\frac{\GR(a+n)}{\GR(a)} .
\end{eqnarray}

\begin{theorem}[General shift identity]\label{thm:shift-general}
Assume $\Hphi$. For $(x,y)\in\Omega_{+}$ and $m,n\in\mathbb{N}\cup\{0\}$,
\begin{eqnarray}\label{eq:shift-general}
	\BR(x+m,y+n)&=&\phi_{1}^{\,m(y-1)+n(x+y)+mn+\frac{n(n-3)}{2}}\,
	\frac{\Pch{\pn{x}}{m}\,\Pch{\pn{y}}{n}}{\Pch{\pn{x+y}}{m+n}}\,\BR(x,y).
\end{eqnarray}
In particular,
\begin{eqnarray}
	\BR(x+n,y)&=&\phi_{1}^{\,n(y-1)}\,\frac{\Pch{\pn{x}}{n}}{\Pch{\pn{x+y}}{n}}\,\BR(x,y),\label{eq:finite-x}\\
	\BR(x,y+n)&=&\phi_{1}^{\,\frac{n(2x+2y+n-3)}{2}}\,\frac{\Pch{\pn{y}}{n}}{\Pch{\pn{x+y}}{n}}\,\BR(x,y),\label{eq:finite-y}\\
	\BR(x+n,y)&=&\phi_{1}^{-\frac{n(2x+n-1)}{2}}\,\frac{\Pch{\pn{x}}{n}}{\Pch{\pn{y}}{n}}\,\BR(x,y+n).\label{eq:finite-mixed}
\end{eqnarray}
\end{theorem}
{\it Proof.}
By \eqref{eq:beta1-gamma} and \eqref{eq:pochhammer},
\begin{eqnarray}
	\frac{\BR(x+m,y+n)}{\BR(x,y)}&=&\phi_{1}^{E(x+m,y+n)-E(x,y)}\,
	\frac{\GR(x+m)\GR(y+n)\GR(x+y)}{\GR(x)\GR(y)\GR(x+y+m+n)},\nn
\end{eqnarray}
and $E(x+m,y+n)-E(x,y)=m(y-1)+n(x+y)+mn+n(n-3)/2$ by Appendix \ref{app:exponents}. The particular cases correspond to $(m,n)$ replaced by $(n,0)$ and $(0,n)$, and \eqref{eq:finite-mixed} is the quotient of \eqref{eq:finite-x} by \eqref{eq:finite-y}, the exponent being $n(y-1)-n(2x+2y+n-3)/2=-n(2x+n-1)/2$. \hfill\(\Box\)

\begin{corollary}[Shift sums]\label{cor:shift-sums}
Assume $\Hphi$. For $(x,y)\in\Omega_{+}$ and $n\in\mathbb{N}$,
\begin{eqnarray}
	\BR(x+1,y)+\BR(x,y+1)&=&\frac{\phi_{1}^{\,y-1}\pn{x}+\phi_{1}^{\,x+y-1}\pn{y}}{\pn{x+y}}\,\BR(x,y),\label{eq:shift-sum1}\\
	\BR(x+n,y)+\BR(x,y+n)&=&\frac{\phi_{1}^{\,n(y-1)}(\phi_{1}^{x}\ominus\phi_{2}^{x})^{n}_{\mathcal{R}(p,q)}+\phi_{1}^{\,\frac{n(2x+2y+n-3)}{2}}(\phi_{1}^{y}\ominus\phi_{2}^{y})^{n}_{\mathcal{R}(p,q)}}{(\phi_{1}^{x+y}\ominus\phi_{2}^{x+y})^{n}_{\mathcal{R}(p,q)}}\,\BR(x,y),\label{eq:shift-sumn}
\end{eqnarray}
and, more generally, $(\phi_{1}^{a}\ominus\phi_{2}^{a})^{n}_{\mathcal{R}(p,q)}=(\phi_{1}-\phi_{2})^{n}\Pch{\pn{a}}{n}$ for every $a\in\mathbb{C}$.
\end{corollary}
{\it Proof.}
By \eqref{eq:powers-int}, $(\phi_{1}^{a}\ominus\phi_{2}^{a})^{n}_{\mathcal{R}(p,q)}=\prod_{k=0}^{n-1}(\phi_{1}^{a+k}-\phi_{2}^{a+k})=(\phi_{1}-\phi_{2})^{n}\prod_{k=0}^{n-1}\pn{a+k}$. Formula \eqref{eq:shift-sumn} is the sum of \eqref{eq:finite-x} and \eqref{eq:finite-y}, the factors $(\phi_{1}-\phi_{2})^{n}$ cancelling, and \eqref{eq:shift-sum1} is its case $n=1$. \hfill\(\Box\)

\begin{remark}\label{rem:finite-exponents}
The exponent in \eqref{eq:finite-y} is quadratic in $n$, whereas the exponent in \eqref{eq:finite-x} is linear in $n$: the iteration in the second argument accumulates the exponents $x+y-1,x+y,\dots,x+y+n-2$ of \eqref{eq:shift-y}. Only the deformed subtraction $\ominus$ occurs in the shift quotients, since $\Pch{\pn{a}}{n}$ is a product of the factors $\phi_{1}^{a+k}-\phi_{2}^{a+k}$. In particular, the iterated forms with linear exponent $n(x+y-1)$ in the second argument, with the deformed addition $\oplus$ in the denominators, or with the exponent $x+2y+1$ in the double shift, which can be found in preliminary versions of these identities, are not correct for $n\geq2$, respectively for $\phi_{1}\neq1$; the exponents of \eqref{eq:shift-xy}, \eqref{eq:finite-y}, \eqref{eq:finite-mixed} and \eqref{eq:shift-sumn} are checked numerically in Section \ref{sec:numerics}.
\end{remark}

\subsection{Symmetrisation}\label{subsec:beta1-sym}

\begin{definition}\label{def:beta1-sym}
Assume $\Hphi$. The symmetrised first kind deformed beta function is
\begin{eqnarray}\label{eq:beta1-sym-def}
	\BRs(x,y)&:=&\phi_{1}^{\frac{x(x-1)}{2}}\,\BR(x,y)\\&=&\;\int_{0}^{1}\phi_{1}^{\frac{x(x-1)}{2}}\,t^{x-1}(1\ominus\phi_{2}t)^{y-1}_{\mathcal{R}(p,q)}\,\dR t,\qquad (x,y)\in\Omega_{+}.
\end{eqnarray}
\end{definition}

\begin{theorem}\label{thm:sym1}
Assume $\Hphi$. For $(x,y)\in\Omega_{+}$,
\begin{eqnarray}\label{eq:beta1-sym}
	\BRs(x,y)&=&\phi_{1}^{\frac{(x+y-1)(x+y-2)}{2}}\,\frac{\GR(x)\,\GR(y)}{\GR(x+y)}
	\;=\;\frac{\GR(x)\,\GR(y)}{\Gq(x+y)}\\
	\;&=&\;\phi_{1}^{\frac{(x-1)(x-2)+(y-1)(y-2)}{2}}\,\Bq(x,y),
\end{eqnarray}
so that $\BRs(x,y)=\BRs(y,x)$. Moreover, if $h:\{\Re x>0\}\to\mathbb{C}$ is such that $(x,y)\mapsto\phi_{1}^{h(x)}\BR(x,y)$ is symmetric on $\Omega_{+}$ and $\phi_{1}\neq1$, then $h(x)-x(x-1)/2$ is constant modulo $2\pi i/\log\phi_{1}$, that is the symmetrising factor is unique up to a constant factor.
\end{theorem}
{\it Proof.}
By \eqref{eq:beta1-gamma}, $\BRs(x,y)=\phi_{1}^{x(x-1)/2+E(x,y)}\GR(x)\GR(y)/\GR(x+y)$, and $x(x-1)/2+E(x,y)=(x+y-1)(x+y-2)/2$ (Appendix \ref{app:exponents}). The second expression follows from $\GR(x+y)=\phi_{1}^{(x+y-1)(x+y-2)/2}\Gq(x+y)$, and the third one from $\GR(x)\GR(y)=\phi_{1}^{[(x-1)(x-2)+(y-1)(y-2)]/2}\Gq(x)\Gq(y)$. For the uniqueness, by \eqref{eq:exchange1} the symmetry of $\phi_{1}^{h(x)}\BR(x,y)$ is equivalent to $\phi_{1}^{h(x)-h(y)+(y-x)(x+y-1)/2}=1$, and $(y-x)(x+y-1)/2=[y(y-1)-x(x-1)]/2$. Hence $\phi_{1}^{u(x)-u(y)}=1$ with $u(x)=h(x)-x(x-1)/2$, for all $x,y$. \hfill\(\Box\)

For integer arguments and the Jagannathan Srinivasa Rao deformation, \eqref{eq:beta1-sym} reduces to the commutative $(p,q)-$beta function of \cite{ref55}, whose exponent $(2mn+m^{2}+n^{2}-3m-3n+2)/2$ is precisely $(m+n-1)(m+n-2)/2$.

\subsection{Products of deformed beta functions}\label{subsec:beta1-products}

In the classical theory, $\mathrm{B}(x,y)\mathrm{B}(x+y,z)=\Gamma(x)\Gamma(y)\Gamma(z)/\Gamma(x+y+z)$ is invariant under all permutations of $(x,y,z)$, which expresses the associativity of the Dirichlet integral. In the deformed setting this invariance is broken by explicit valuation factors.

\begin{proposition}\label{prop:products}
Assume $\Hphi$ and let $x,y,z\in\mathbb{C}$ with positive real parts. Then
\begin{eqnarray}\label{eq:product-gamma}
	\BR(x,y)\,\BR(x+y,z)&=&\phi_{1}^{\frac{(y-1)(2x+y-2)+(z-1)(2x+2y+z-2)}{2}}\,\frac{\GR(x)\,\GR(y)\,\GR(z)}{\GR(x+y+z)},
\end{eqnarray}
and consequently
\begin{eqnarray}
	\BR(x,y)\,\BR(x+y,z)&=&\phi_{1}^{\,z-y}\,\BR(x,z)\,\BR(x+z,y),\label{eq:product-swap}\\
	\BR(x,y)\,\BR(x+y,z)&=&\phi_{1}^{-\frac{2x+2yz+z^{2}-3z}{2}}\,\BR(y,z)\,\BR(x,y+z).\label{eq:product-assoc}
\end{eqnarray}
For the symmetrised function,
\begin{eqnarray}\label{eq:product-sym}
	\BRs(x,y)\,\BRs(x+y,z)&=&\phi_{1}^{\frac{(x+y-1)(x+y-2)+(x+y+z-1)(x+y+z-2)}{2}}\,\frac{\GR(x)\,\GR(y)\,\GR(z)}{\GR(x+y+z)}.
\end{eqnarray}
\end{proposition}
{\it Proof.}
Formula \eqref{eq:product-gamma} follows from \eqref{eq:beta1-gamma} applied twice, the factor $\GR(x+y)$ cancelling. The gamma quotient on the right hand side is symmetric in $(x,y,z)$, so that the products in \eqref{eq:product-swap} and \eqref{eq:product-assoc} differ only through the exponents. With $E(a,b)=(b-1)(2a+b-2)/2$ a direct expansion gives
$E(x,y)+E(x+y,z)-E(x,z)-E(x+z,y)=z-y$ and $E(x,y)+E(x+y,z)-E(y,z)-E(x,y+z)=-(2x+2yz+z^{2}-3z)/2$. The last statement follows from \eqref{eq:beta1-sym}. \hfill\(\Box\)

For $x=y=z=1$ one finds $\BR(1,1)\BR(2,1)=1/\pn{2}$ and $\BR(1,1)\BR(1,2)=\phi_{1}/\pn{2}$ by Corollary \ref{cor:special-values}, in agreement with \eqref{eq:product-assoc}, whose exponent equals $-1$ at this point. The defects \eqref{eq:product-swap} and \eqref{eq:product-assoc} vanish in the Jackson case, where the associativity of the Dirichlet type products is restored.

\subsection{Asymptotic behaviour}\label{subsec:beta1-asymptotics}

\begin{theorem}\label{thm:beta1-asymptotics}
Assume $\Hphi$.
\begin{enumerate}
\item[\textnormal{(i)}] For fixed $y$ with $\Re y>0$, as $x\to+\infty$ along the real axis,
\begin{eqnarray}\label{eq:asym-x}
	\BR(x,y)&=&\phi_{1}^{\,1-x+\frac{(y-1)(y-2)}{2}}\,(1-\fq)^{y}\,\Gq(y)\,\big(1+O(\fq^{x})\big)\\
	\;&=&\;\phi_{1}^{(x-1)(y-1)}\,\frac{\GR(y)}{\pn{x}^{\,y}}\,\big(1+O(\fq^{x})\big).
\end{eqnarray}
\item[\textnormal{(ii)}] For fixed $x$ with $\Re x>0$, as $y\to+\infty$ along the real axis,
\begin{eqnarray}\label{eq:asym-y}
	\BR(x,y)&=&\phi_{1}^{\,1-x+\frac{(y-1)(y-2)}{2}}\,(1-\fq)^{x}\,\Gq(x)\,\big(1+O(\fq^{y})\big).
\end{eqnarray}
In particular, for real $x>0$, $\BR(x,y)\to+\infty$ if $\phi_{1}>1$, $\BR(x,y)\to(1-\fq)^{x}\Gq(x)$ if $\phi_{1}=1$, and $\BR(x,y)\to0$ if $\phi_{1}<1$, with $\log\BR(x,y)\sim\tfrac{y^{2}}{2}\log\phi_{1}$ when $\phi_{1}\neq1$.
\end{enumerate}
\end{theorem}
{\it Proof.}
By Theorem \ref{thm:beta1}, $\BR=\phi_{1}^{1-x+(y-1)(y-2)/2}\Bq$, and by \eqref{eq:Jackson-gamma},
\begin{eqnarray}
	\Bq(x,y)&=&\Gq(y)\,(1-\fq)^{y}\,\frac{(\fq^{x+y};\fq)_{\infty}}{(\fq^{x};\fq)_{\infty}}
	\;=\;\Gq(x)\,(1-\fq)^{x}\,\frac{(\fq^{x+y};\fq)_{\infty}}{(\fq^{y};\fq)_{\infty}} .\nn
\end{eqnarray}
For $u\to+\infty$ and bounded $\Re v>0$, $\log(\fq^{u+v};\fq)_{\infty}-\log(\fq^{u};\fq)_{\infty}=O(\fq^{u})$, which gives the first expressions in \eqref{eq:asym-x} and \eqref{eq:asym-y}. For the second expression in \eqref{eq:asym-x}, Lemma \ref{lem:numbers} gives $\pn{x}^{-y}=\phi_{1}^{-y(x-1)}(1-\fq)^{y}(1-\fq^{x})^{-y}$, and $\GR(y)=\phi_{1}^{(y-1)(y-2)/2}\Gq(y)$; the exponents satisfy $(x-1)(y-1)-y(x-1)+\tfrac{(y-1)(y-2)}{2}=1-x+\tfrac{(y-1)(y-2)}{2}$. The last assertions follow from the sign of $\log\phi_{1}$. \hfill\(\Box\)

The behaviour \eqref{eq:asym-y} departs sharply from the classical one, $\mathrm{B}(x,y)\sim\Gamma(x)y^{-x}$: in the deformed setting the decay in the second argument is replaced by a Gaussian growth or decay governed by the valuation factor, and only in the Jackson case $\phi_{1}=1$ does the deformed beta function tend to a finite nonzero limit. On the contrary, the behaviour in the first argument \eqref{eq:asym-x} is the exact deformed analogue of $\mathrm{B}(x,y)\sim\Gamma(y)x^{-y}$, the power $x^{-y}$ being replaced by $\pn{x}^{-y}$ and the valuation factor by $\phi_{1}^{(x-1)(y-1)}$, in accordance with the Wendel type limit \eqref{eq:wendel-limit}.

\begin{proposition}[Monotonicity in the first argument]\label{prop:monotonicity}
Assume $\Hphi$ and let $y>0$. For $x>0$,
\begin{eqnarray}\label{eq:dlogB}
	\frac{\partial}{\partial x}\log\BR(x,y)&=&-\log\phi_{1}+\log\fq\sum_{k\geq0}\Big(\frac{\fq^{x+k}}{1-\fq^{x+k}}-\frac{\fq^{x+y+k}}{1-\fq^{x+y+k}}\Big),
\end{eqnarray}
and the sum is positive. Consequently:
\begin{enumerate}
\item[\textnormal{(i)}] if $\phi_{1}\geq1$, the function $x\mapsto\BR(x,y)$ is strictly decreasing on $(0,\infty)$;
\item[\textnormal{(ii)}] if $\phi_{1}<1$, there exists a unique $x_{\min}(y)>0$ such that $x\mapsto\BR(x,y)$ is strictly decreasing on $(0,x_{\min}(y)]$ and strictly increasing on $[x_{\min}(y),\infty)$, and $\BR(x,y)\to+\infty$ as $x\to+\infty$.
\end{enumerate}
\end{proposition}
{\it Proof.}
By Theorem \ref{thm:beta1}, $\log\BR(x,y)=(1-x)\log\phi_{1}+\tfrac{(y-1)(y-2)}{2}\log\phi_{1}+\log\Gq(x)+\log\Gq(y)-\log\Gq(x+y)$, and logarithmic differentiation of \eqref{eq:Jackson-gamma} gives $\frac{d}{dx}\log\Gq(x)=-\log(1-\fq)+\log\fq\sum_{k\geq0}\fq^{x+k}/(1-\fq^{x+k})$, which yields \eqref{eq:dlogB}. Since $u\mapsto u/(1-u)$ is increasing on $(0,1)$ and $\fq^{x+k}>\fq^{x+y+k}$, each term of the sum is positive, and $\log\fq<0$. The second term of \eqref{eq:dlogB} is therefore negative; it tends to $-\infty$ as $x\to0^{+}$ and to $0$ as $x\to+\infty$, and it is increasing in $x$ by the logarithmic convexity of Theorem \ref{thm:convex-beta}(i). If $\phi_{1}\geq1$ the derivative is negative. If $\phi_{1}<1$, the derivative increases from $-\infty$ to $-\log\phi_{1}>0$ and has a unique zero, and $\BR(x,y)\to+\infty$ by \eqref{eq:asym-x}. \hfill\(\Box\)

\begin{figure}[ht]
\centering
\includegraphics[width=\linewidth]{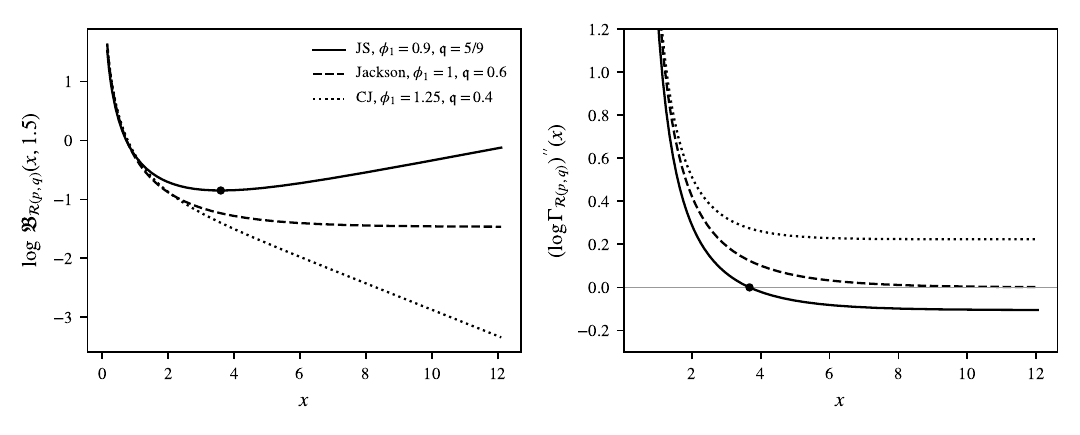}
\caption{Left: $x\mapsto\log\BR(x,1.5)$ for the Jagannathan Srinivasa Rao deformation $(p,q)=(0.9,0.5)$, the Jackson case $q=0.6$ and the Chakrabarti Jagannathan deformation $(p,q)=(0.8,0.5)$; the dot marks the minimum $x_{\min}(1.5)=3.5936\ldots$ of Proposition \ref{prop:monotonicity}(ii). Right: $x\mapsto\frac{d^{2}}{dx^{2}}\log\GR(x)$ for the same deformations; the dot marks the zero $x_{*}=3.6722\ldots$ of Theorem \ref{thm:logconvex-gamma}(ii).}\label{fig:beta-gamma}
\end{figure}

Figure \ref{fig:beta-gamma} illustrates the three regimes $\phi_{1}<1$, $\phi_{1}=1$ and $\phi_{1}>1$. In the first one, the deformed beta function has a minimum in its first argument and the deformed gamma function loses its logarithmic convexity beyond $x_{*}$; in the last one both functions behave as in the classical theory, the deformed beta function being decreasing and the deformed gamma function strictly logarithmically convex.

\section{The second kind noncommutative deformed beta function}\label{sec:beta2}

The classical beta function admits, besides the Euler integral on $[0,1]$, the representation $\int_{0}^{\infty}t^{x-1}(1+t)^{-x-y}\,dt$ on the half line. Its deformed counterpart is built from the deformed addition $\oplus$ and from the bilateral deformed integral.

\subsection{Definition and closed form}\label{subsec:beta2-def}

\begin{definition}\label{def:beta2}
Assume $\Hphi$. For $(x,y)\in\Omega_{+}$ the second kind $\mathcal{R}(p,q)-$beta function is
\begin{eqnarray}\label{eq:beta2-def}
	\BB(x,y)&:=&\int_{0}^{\infty;1}\frac{t^{x-1}}{(1\oplus\phi_{1}t)^{x+y}_{\mathcal{R}(p,q)}}\,\dR t,
\end{eqnarray}
where, by \eqref{eq:powers-complex}, $1/(1\oplus\phi_{1}t)^{s}_{\mathcal{R}(p,q)}=\phi_{1}^{-s(s-1)/2}(-\phi_{1}\fq^{s}t;\fq)_{\infty}/(-\phi_{1}t;\fq)_{\infty}$, an entire function of $s$ for every $t>0$.
\end{definition}

\begin{lemma}\label{lem:beta2-series}
Assume $\Hphi$ and put $s=x+y$. For $(x,y)\in\Omega_{+}$,
\begin{eqnarray}\label{eq:beta2-series}
	\BB(x,y)&=&(1-\fq)\,\phi_{1}^{1-x-\frac{s(s-1)}{2}}\sum_{k\in\mathbb{Z}}\fq^{kx}\,\frac{(-\fq^{k+s};\fq)_{\infty}}{(-\fq^{k};\fq)_{\infty}},
\end{eqnarray}
the bilateral series converges absolutely and uniformly on compact subsets of $\Omega_{+}$, $\BB$ is holomorphic on $\Omega_{+}$, and $\BB(x,y)>0$ for real $x,y>0$.
\end{lemma}
{\it Proof.}
The nodes are $t_{k}=\fq^{k}/\phi_{1}$, $k\in\mathbb{Z}$, the weights are $(1-\fq)\fq^{k}$ and $\phi_{1}t_{k}=\fq^{k}$, which gives \eqref{eq:beta2-series}. For $k\geq0$ the quotient of products is bounded by $(-\fq^{\Re s};\fq)_{\infty}$, so that the terms are $O(\fq^{k\Re x})$. For $k=-m\leq0$, the identity $(-\fq^{-m}a;\fq)_{\infty}=(-a;\fq)_{\infty}\prod_{i=1}^{m}(1+a\fq^{-i})$ gives
\begin{eqnarray}
	\Big|\frac{(-\fq^{s-m};\fq)_{\infty}}{(-\fq^{-m};\fq)_{\infty}}\Big|
	&=&\frac{|(-\fq^{s};\fq)_{\infty}|}{(-1;\fq)_{\infty}}\prod_{i=1}^{m}\frac{|\fq^{i}+\fq^{s}|}{\fq^{i}+1}
	\;\leq\;\frac{(-\fq^{\Re s};\fq)_{\infty}\,(-\fq^{1-\Re s};\fq)_{\infty}}{(-1;\fq)_{\infty}}\,\fq^{m\Re s},\nn
\end{eqnarray}
since $|\fq^{i}+\fq^{s}|\leq\fq^{\Re s}(1+\fq^{i-\Re s})$. The $k-$th term is therefore $O(\fq^{m\Re y})$ for $k=-m$, uniformly on compact subsets of $\Omega_{+}$, and holomorphy follows from the Weierstrass theorem. Positivity follows from Lemma \ref{lem:powers}(v). \hfill\(\Box\)

\begin{theorem}\label{thm:beta2}
Assume $\Hphi$. For $(x,y)\in\Omega_{+}$,
\begin{eqnarray}
	\BB(x,y)&=&\phi_{1}^{\,1-x-\frac{(x+y)(x+y-1)}{2}}\,\Kq(x)\,\Bq(x,y)\label{eq:beta2-q}\\
	&=&\phi_{1}^{-\frac{x^{2}+x+y^{2}-y}{2}}\,\Kq(x)\,\frac{\GR(x)\,\GR(y)}{\GR(x+y)}\label{eq:beta2-gamma1}\\
	&=&\phi_{1}^{-x-\frac{y(y-1)}{2}}\,\phi_{2}^{-\frac{x(x-1)}{2}}\,\frac{\gR(x)\,\GR(y)}{\GR(x+y)} .\label{eq:beta2-gamma2}
\end{eqnarray}
In particular $\BB$ extends to a meromorphic function on $\mathbb{C}^{2}$, and
\begin{eqnarray}\label{eq:beta2-vs-beta1}
	\BB(x,y)&=&\phi_{1}^{-\frac{(x+y)(x+y-1)}{2}-\frac{(y-1)(y-2)}{2}}\,\Kq(x)\,\BR(x,y).
\end{eqnarray}
\end{theorem}
{\it Proof.}
Since $(-\fq^{k+s};\fq)_{\infty}=(-\fq^{s};\fq)_{\infty}/(-\fq^{s};\fq)_{k}$ and $(-\fq^{k};\fq)_{\infty}=(-1;\fq)_{\infty}/(-1;\fq)_{k}$ for $k\in\mathbb{Z}$,
\begin{eqnarray}
	\sum_{k\in\mathbb{Z}}\fq^{kx}\frac{(-\fq^{k+s};\fq)_{\infty}}{(-\fq^{k};\fq)_{\infty}}
	&=&\frac{(-\fq^{s};\fq)_{\infty}}{(-1;\fq)_{\infty}}\sum_{k\in\mathbb{Z}}\frac{(-1;\fq)_{k}}{(-\fq^{s};\fq)_{k}}\,\fq^{kx}.\nn
\end{eqnarray}
Ramanujan's summation \eqref{eq:psi11} with $a=-1$, $b=-\fq^{s}$ and $w=\fq^{x}$, whose condition $|b/a|=|\fq^{x+y}|<|\fq^{x}|<1$ is equivalent to $(x,y)\in\Omega_{+}$, gives
\begin{eqnarray}
	\sum_{k\in\mathbb{Z}}\frac{(-1;\fq)_{k}}{(-\fq^{s};\fq)_{k}}\,\fq^{kx}
	&=&\frac{(\fq,\fq^{x+y},-\fq^{x},-\fq^{1-x};\fq)_{\infty}}{(-\fq^{x+y},-\fq,\fq^{x},\fq^{y};\fq)_{\infty}} .\nn
\end{eqnarray}
Hence the bilateral sum equals $\Kq(x)(\fq,\fq^{x+y};\fq)_{\infty}/(\fq^{x},\fq^{y};\fq)_{\infty}=\Kq(x)\Bq(x,y)/(1-\fq)$, and \eqref{eq:beta2-q} follows from \eqref{eq:beta2-series}. Formula \eqref{eq:beta2-gamma1} follows from $\Bq(x,y)=\phi_{1}^{xy-1}\GR(x)\GR(y)/\GR(x+y)$ (proof of Theorem \ref{thm:beta1}) and from the exponent identity of Appendix \ref{app:exponents}. Since $\Kq(x)=\fq^{-x(x-1)/2}\Omq(x)$ and $\gR=\Omq\GR$ (Theorem \ref{thm:gamma2}), we have $\Kq(x)\GR(x)=\phi_{1}^{x(x-1)/2}\phi_{2}^{-x(x-1)/2}\gR(x)$, which turns \eqref{eq:beta2-gamma1} into \eqref{eq:beta2-gamma2}. Finally \eqref{eq:beta2-vs-beta1} is the quotient of \eqref{eq:beta2-q} by the second expression in \eqref{eq:beta1-gamma}. \hfill\(\Box\)

Formula \eqref{eq:beta2-gamma2} shows that the second kind deformed beta function pairs the second kind deformed gamma function in its first argument with the first kind deformed gamma function in its second argument and in the denominator; both deformation functions $\phi_{1}$ and $\phi_{2}$ enter the valuation factor.

\begin{corollary}\label{cor:beta2-values}
Assume $\Hphi$. For $y\in\mathbb{C}$ with $\Re y>0$ and $m,n\in\mathbb{N}$,
\begin{eqnarray}
	\BB(1,y)&=&\frac{\phi_{1}^{-1-\frac{y(y-1)}{2}}}{\pn{y}},\qquad
	\BB(m,n)\;=\;\phi_{1}^{-m-\frac{n(n-1)}{2}}\,\phi_{2}^{-\frac{m(m-1)}{2}}\,\frac{\pf{m-1}\,\pf{n-1}}{\pf{m+n-1}} .\nn
\end{eqnarray}
In particular $\BB(1,1)=1/\phi_{1}$.
\end{corollary}
{\it Proof.}
Use \eqref{eq:beta2-gamma2} with $\gR(m)=\GR(m)=\pf{m-1}$ and $\gR(1)=1$. The value $\BB(1,1)=1/\phi_{1}$ can also be obtained directly: $(1\oplus\phi_{1}t_{k})^{2}_{\mathcal{R}(p,q)}=\phi_{1}(1+\fq^{k})(1+\fq^{k+1})$, and the bilateral sum $\sum_{k\in\mathbb{Z}}(1-\fq)\fq^{k}/((1+\fq^{k})(1+\fq^{k+1}))=\sum_{k}[(1+\fq^{k+1})^{-1}-(1+\fq^{k})^{-1}]$ telescopes to $1$. \hfill\(\Box\)

\subsection{Exchange rule, shifts and symmetrisation}\label{subsec:beta2-props}

\begin{theorem}\label{thm:exchange2}
Assume $\Hphi$. For $(x,y)\in\Omega_{+}$,
\begin{eqnarray}\label{eq:exchange2}
	\BB(x,y)&=&\phi_{1}^{\frac{(y-x)(3-x-y)}{2}}\,\phi_{2}^{\frac{(y-x)(x+y-1)}{2}}\,\frac{\Omq(x)}{\Omq(y)}\,\BB(y,x).
\end{eqnarray}
In particular, in the Jackson case $\phi_{1}=1$, $\BB(x,y)=\fq^{(y-x)(x+y-1)/2}\,\Omq(x)\Omq(y)^{-1}\,\BB(y,x)$, and $\BB$ is not commutative.
\end{theorem}
{\it Proof.}
By \eqref{eq:beta2-gamma2} and $\gR=\Omq\GR$,
\begin{eqnarray}
	\frac{\BB(x,y)}{\BB(y,x)}&=&\phi_{1}^{\,y-x+\frac{x(x-1)-y(y-1)}{2}}\,\phi_{2}^{\frac{y(y-1)-x(x-1)}{2}}\,\frac{\Omq(x)}{\Omq(y)},\nn
\end{eqnarray}
and $x(x-1)-y(y-1)=(x-y)(x+y-1)$. For $\phi_{1}=1$, $\phi_{2}=\fq$. Since $\Omq$ is one periodic and, by Proposition \ref{prop:Omega}, takes values close to but different from $1$, the factor in front of $\BB(y,x)$ is not identically $1$. \hfill\(\Box\)

The comparison of Theorems \ref{thm:exchange1} and \ref{thm:exchange2} is instructive: the first kind exchange defect is carried by $\phi_{1}$ alone and vanishes in the Jackson case, while the second kind exchange defect is carried by both deformation functions and by the theta quotient $\Omq$, and survives in the Jackson case.

\begin{theorem}\label{thm:beta2-shifts}
Assume $\Hphi$. For $(x,y)\in\Omega_{+}$ and $m,n\in\mathbb{N}\cup\{0\}$,
\begin{eqnarray}\label{eq:beta2-shift-general}
	\BB(x+m,y+n)&=&\phi_{1}^{-m-ny-\frac{n(n-1)}{2}}\,\phi_{2}^{-mx-\frac{m(m-1)}{2}}\,
	\frac{\Pch{\pn{x}}{m}\,\Pch{\pn{y}}{n}}{\Pch{\pn{x+y}}{m+n}}\,\BB(x,y).
\end{eqnarray}
In particular
\begin{eqnarray}
	\BB(x+1,y)\;&=&\;\frac{\pn{x}}{\phi_{1}\phi_{2}^{x}\,\pn{x+y}}\,\BB(x,y),\\
	\BB(x,y+1)\;&=&\;\frac{\pn{y}}{\phi_{1}^{y}\,\pn{x+y}}\,\BB(x,y).\label{eq:beta2-unit}
\end{eqnarray}
\end{theorem}
{\it Proof.}
By \eqref{eq:beta2-gamma2}, the quotient $\BB(x+m,y+n)/\BB(x,y)$ equals the Pochhammer quotient times
\begin{eqnarray}
	\phi_{1}^{-m-[(y+n)(y+n-1)-y(y-1)]/2}\,\phi_{2}^{-[(x+m)(x+m-1)-x(x-1)]/2},\nn
\end{eqnarray}
because $\gR(x+m)/\gR(x)=\Pch{\pn{x}}{m}$ by the periodicity of $\Omq$. One has $(y+n)(y+n-1)-y(y-1)=2ny+n(n-1)$ and $(x+m)(x+m-1)-x(x-1)=2mx+m(m-1)$. \hfill\(\Box\)

\begin{theorem}\label{thm:beta2-sym}
Assume $\Hphi$. The function
\begin{eqnarray}\label{eq:beta2-sym-def}
	\BBs(x,y)&:=&\frac{\phi_{1}^{x}}{\Kq(x)}\,\BB(x,y)\;=\;\frac{\phi_{1}^{x}\,\fq^{\frac{x(x-1)}{2}}}{\Omq(x)}\,\BB(x,y),\qquad (x,y)\in\Omega_{+},
\end{eqnarray}
satisfies
\begin{eqnarray}\label{eq:beta2-sym}
	\BBs(x,y)&=&\phi_{1}^{-\frac{x(x-1)+y(y-1)}{2}}\,\frac{\GR(x)\,\GR(y)}{\GR(x+y)}
	\;=\;\phi_{1}^{-\frac{x(x-1)+y(y-1)+(x+y-1)(x+y-2)}{2}}\,\BRs(x,y),
\end{eqnarray}
and is therefore commutative.
\end{theorem}
{\it Proof.}
Multiply \eqref{eq:beta2-gamma1} by $\phi_{1}^{x}/\Kq(x)$ and use $-(x^{2}+x+y^{2}-y)/2+x=-[x(x-1)+y(y-1)]/2$; the second expression follows from \eqref{eq:beta1-sym}. \hfill\(\Box\)

\begin{remark}\label{rem:beta2-limit}
In the undeformed limit $\phi_{1}\to1$ and $\fq\to1^{-}$ one has $\varepsilon=e^{-2\pi^{2}/\log(1/\fq)}\to0$, hence $\Omq\to1$ and $\Kq(x)=\fq^{-x(x-1)/2}\Omq(x)\to1$ by Proposition \ref{prop:Omega}, while $\Gq\to\Gamma$ \cite{Askey1978, GasperRahman}. Formula \eqref{eq:beta2-q} then reduces to the Euler integral $\int_{0}^{\infty}t^{x-1}(1+t)^{-x-y}dt=\Gamma(x)\Gamma(y)/\Gamma(x+y)$. For $\phi_{1}=1$ and $\fq<1$, \eqref{eq:beta2-q} is the $q-$beta integral on $(0,\infty)$ of Ramanujan type discussed in \cite{Askey1978, GasperRahman}.
\end{remark}

\subsection{The finite interval variant}\label{subsec:beta2-variant}

A second kind function on $[0,1]$ can also be built from the reversed deformed power $(\phi_{1}t\ominus1)^{n}_{\mathcal{R}(p,q)}=\prod_{k=0}^{n-1}(\phi_{1}^{k+1}t-\phi_{2}^{k})$. For integer second argument it is explicit, and it shows why the half line definition \eqref{eq:beta2-def} is the appropriate one.

\begin{proposition}\label{prop:beta2-variant}
Assume $\Hphi$, let $\Re x>0$ and $n\in\mathbb{N}\cup\{0\}$, and put $\mathcal{B}^{[0,1]}_{\mathcal{R}(p,q)}(x,n+1):=\int_{0}^{1}t^{x-1}(\phi_{1}t\ominus1)^{n}_{\mathcal{R}(p,q)}\,\dR t$. Then
\begin{eqnarray}\label{eq:beta2-variant}
	\mathcal{B}^{[0,1]}_{\mathcal{R}(p,q)}(x,n+1)&=&(-1)^{n}\,\phi_{2}^{\frac{n(2x+n-1)}{2}}\,\frac{\pf{n}}{\prod_{j=0}^{n}\pn{x+j}}\\
	\;&=&\;(-1)^{y-1}\phi_{2}^{\frac{(y-1)(2x+y-2)}{2}}\,\frac{\gR(x)\,\gR(y)}{\gR(x+y)}\Big|_{y=n+1},
\end{eqnarray}
and consequently, for $y\in\mathbb{N}$,
\begin{eqnarray}
	\mathcal{B}^{[0,1]}_{\mathcal{R}(p,q)}(x+1,y)&=&\frac{\phi_{2}^{\,y-1}\pn{x}}{\pn{x+y}}\,\mathcal{B}^{[0,1]}_{\mathcal{R}(p,q)}(x,y),\qquad
	\mathcal{B}^{[0,1]}_{\mathcal{R}(p,q)}(x,y+1)\;=\;-\frac{\phi_{2}^{\,x+y-1}\pn{y}}{\pn{x+y}}\,\mathcal{B}^{[0,1]}_{\mathcal{R}(p,q)}(x,y).\nn
\end{eqnarray}
\end{proposition}
{\it Proof.}
Since $(u\ominus1)^{n}_{\mathcal{R}(p,q)}=u^{n}(1\ominus u^{-1})^{n}_{\mathcal{R}(p,q)}$, the deformed Gauss formula of Lemma \ref{lem:powers}(iii) with $u=\phi_{1}t$, Proposition \ref{prop:integral}(ii) and the change $k\mapsto n-k$ give
\begin{eqnarray}
	\mathcal{B}^{[0,1]}_{\mathcal{R}(p,q)}(x,n+1)&=&(-1)^{n}\sum_{k=0}^{n}(-1)^{k}\bigg[\begin{array}{c} n \\ k\end{array}\bigg]_{\mathcal{R}(p,q)}\frac{\phi_{2}^{\binom{n-k}{2}}\phi_{1}^{\binom{k+1}{2}}}{\pn{x+k}} .\nn
\end{eqnarray}
This is the left hand side of \eqref{eq:partial-fractions} with $\phi_{1}$ and $\phi_{2}$ exchanged. Since the deformed numbers and binomial coefficients are symmetric in $(\phi_{1},\phi_{2})$, both sides of \eqref{eq:partial-fractions} are real analytic functions of $(\phi_{1},\phi_{2})\in(0,\infty)^{2}$ for real $x>0$, and the identity, valid for $\phi_{1}>\phi_{2}$, extends to $\phi_{1}<\phi_{2}$. Exchanging the roles of $\phi_{1}$ and $\phi_{2}$ gives the first expression in \eqref{eq:beta2-variant} for real $x>0$, and for $\Re x>0$ by analytic continuation. The second expression follows from $\gR(m)=\pf{m-1}$, $\gR(x+n+1)=\gR(x)\prod_{j=0}^{n}\pn{x+j}$ and $(y-1)(2x+y-2)/2=n(2x+n-1)/2$. The recurrences are quotients of consecutive values. \hfill\(\Box\)

Formula \eqref{eq:beta2-variant} agrees with the relation $\phi_{2}^{(y-1)(2x+y-2)/2}\gR(x)\gR(y)/\gR(x+y)$ proposed for this variant in preliminary versions only up to the sign $(-1)^{y-1}$: the integrand changes sign on the lattice, the function alternates in sign along $y\in\mathbb{N}$, and the second recurrence carries the factor $-1$. Moreover, for $t\in(0,1)$ the factors $\phi_{1}^{k+1}t-\phi_{2}^{k}$ of the infinite product $(\phi_{1}t\ominus1)^{\infty}_{\mathcal{R}(p,q)}$ do not tend to $1$, so that the quotient of infinite products used for complex $y$ in those versions is not defined. The half line function $\BB$ of Definition \ref{def:beta2} is positive, defined on all of $\Omega_{+}$ and carries the same valuation structure; it is the second kind function studied in this paper.

\section{Logarithmic convexity}\label{sec:convexity}

\begin{theorem}\label{thm:convex-beta}
Assume $\Hphi$. Then:
\begin{enumerate}
\item[\textnormal{(i)}] the function $(x,y)\mapsto\phi_{1}^{-(y-1)(y-2)/2}\BR(x,y)$ is jointly logarithmically convex on $(0,\infty)^{2}$; in particular, for every fixed $y>0$, $x\mapsto\BR(x,y)$ is logarithmically convex on $(0,\infty)$, and if $\phi_{1}\geq1$ then $\BR$ itself is jointly logarithmically convex;
\item[\textnormal{(ii)}] the function $(x,y)\mapsto\phi_{1}^{-[(x-1)(x-2)+(y-1)(y-2)]/2}\BRs(x,y)=\Bq(x,y)$ is jointly logarithmically convex on $(0,\infty)^{2}$;
\item[\textnormal{(iii)}] the function $(x,y)\mapsto\phi_{1}^{(x+y)(x+y-1)/2}\BB(x,y)$ is jointly logarithmically convex on $(0,\infty)^{2}$; in particular, if $\phi_{1}\leq1$, then $x\mapsto\BB(x,y)$ is logarithmically convex on $(0,\infty)$ for every fixed $y>0$.
\end{enumerate}
\end{theorem}
{\it Proof.}
Recall that a positive function is logarithmically convex if its logarithm is convex, that products and pointwise convergent sums of logarithmically convex functions are logarithmically convex \cite{Artin1964}, and that the sum of a convex function of $x$ and of a convex function of $y$ is jointly convex.

(i) By \eqref{eq:beta1-series}, $\phi_{1}^{-(y-1)(y-2)/2}\BR(x,y)=(1-\fq)\phi_{1}^{1-x}\sum_{j\geq0}u_{j}(x,y)$ with
\begin{eqnarray}
	\log u_{j}(x,y)&=&jx\log\fq+\log(\fq^{j+1};\fq)_{\infty}-\sum_{k\geq0}\log\big(1-\fq^{j+k}e^{-\lambda y}\big),\qquad \lambda=\log(1/\fq).\nn
\end{eqnarray}
The first two terms are affine in $x$. For $c\in(0,1)$, the function $\psi_{c}(y)=-\log(1-ce^{-\lambda y})$ satisfies $\psi_{c}''(y)=\lambda^{2}ce^{-\lambda y}(1-ce^{-\lambda y})^{-2}>0$, so that the last sum is convex in $y$ (the derivatives may be taken termwise by uniform convergence on compact sets). Hence each $u_{j}$ is jointly logarithmically convex, the factor $\phi_{1}^{1-x}$ is logarithmically affine, and the sum is logarithmically convex. For fixed $y$ the factor $\phi_{1}^{(y-1)(y-2)/2}$ is a constant, and for $\phi_{1}\geq1$ it is logarithmically convex in $y$.

(ii) is (i) for the Jackson case together with \eqref{eq:beta1-sym}; equivalently it follows from the logarithmic convexity of $\Gq$ and from the joint convexity of $(x,y)\mapsto-\log\Gq(x+y)$ composed with an affine map, but the series argument above gives it directly, since $\Bq(x,y)=(1-\fq)\sum_{j}u_{j}(x,y)$.

(iii) By \eqref{eq:beta2-series}, $\phi_{1}^{s(s-1)/2}\BB(x,y)=(1-\fq)\phi_{1}^{1-x}\sum_{k\in\mathbb{Z}}v_{k}(x,y)$ with $s=x+y$ and
\begin{eqnarray}
	\log v_{k}(x,y)&=&kx\log\fq-\log(-\fq^{k};\fq)_{\infty}+\sum_{j\geq0}\log\big(1+\fq^{k+j}e^{-\lambda(x+y)}\big).\nn
\end{eqnarray}
For $c>0$, $\chi_{c}(s)=\log(1+ce^{-\lambda s})$ satisfies $\chi_{c}''(s)=\lambda^{2}ce^{-\lambda s}(1+ce^{-\lambda s})^{-2}>0$, and the composition with the linear map $(x,y)\mapsto x+y$ is jointly convex. The rest of the argument is as in (i). Finally $\BB=\phi_{1}^{-s(s-1)/2}\cdot\big(\phi_{1}^{s(s-1)/2}\BB\big)$ and $s\mapsto-\tfrac{s(s-1)}{2}\log\phi_{1}$ is convex when $\phi_{1}\leq1$. \hfill\(\Box\)

\begin{corollary}[Deformed Tur\'an identity]\label{cor:turan}
Assume $\Hphi$. For all $x,y\in\mathbb{C}$,
\begin{eqnarray}\label{eq:turan-identity}
	\pn{x+1}\,\pn{x+y}-\pn{x}\,\pn{x+y+1}&=&(\phi_{1}\phi_{2})^{x}\,\pn{y}.
\end{eqnarray}
Consequently, for real $x,y>0$,
\begin{eqnarray}\label{eq:turan}
	\frac{\BR(x+1,y)^{2}}{\BR(x,y)\,\BR(x+2,y)}&=&1-\frac{(\phi_{1}\phi_{2})^{x}\,\pn{y}}{\pn{x+1}\,\pn{x+y}}\;<\;1,
\end{eqnarray}
which is the discrete form of the logarithmic convexity of $x\mapsto\BR(x,y)$ in Theorem \ref{thm:convex-beta}.
\end{corollary}
{\it Proof.}
Multiplying the left hand side of \eqref{eq:turan-identity} by $(\phi_{1}-\phi_{2})^{2}$ and expanding, the terms $\phi_{1}^{2x+y+1}$ and $\phi_{2}^{2x+y+1}$ cancel and the remaining four terms factor as $\phi_{1}^{x}\phi_{2}^{x}(\phi_{1}-\phi_{2})(\phi_{1}^{y}-\phi_{2}^{y})$; see Appendix \ref{app:exponents}. By \eqref{eq:shift-x} applied twice, the left hand side of \eqref{eq:turan} equals $\pn{x}\pn{x+y+1}/(\pn{x+1}\pn{x+y})$, and \eqref{eq:turan-identity} gives the right hand side. \hfill\(\Box\)

\begin{remark}\label{rem:convexity-noncomm}
Theorem \ref{thm:convex-beta} shows that the noncommutativity of the deformed beta functions does not destroy their convexity in the first argument, while the convexity in the second argument requires the normalisation by $\phi_{1}^{-(y-1)(y-2)/2}$ when $\phi_{1}<1$. This is consistent with Theorem \ref{thm:logconvex-gamma}: the valuation factors of the deformed beta functions are quadratic in the arguments, and the sign of $\log\phi_{1}$ decides whether they help or obstruct logarithmic convexity.
\end{remark}

\section{A deformed Euler integral and the deformed Gauss summation}\label{sec:euler}

The first kind deformed beta function is the kernel of a deformed Euler integral representation of a deformed Gauss hypergeometric function. We introduce this function entirely in terms of deformed Pochhammer symbols and deformed factorials.

\begin{definition}\label{def:hypergeometric}
Assume $\Hphi$. For $a,b,c\in\mathbb{C}$ with $\pn{c+j}\neq0$ for $j\in\mathbb{N}\cup\{0\}$, the deformed Gauss hypergeometric function is
\begin{eqnarray}\label{eq:hypergeometric}
	{}_{2}F_{1}^{\mathcal{R}(p,q)}(a,b;c;w)&:=&\sum_{n\geq0}\frac{\Pch{\pn{a}}{n}\,\Pch{\pn{b}}{n}}{\Pch{\pn{c}}{n}\,\pf{n}}\,w^{n},
\end{eqnarray}
for $|w|<\phi_{1}^{\Re(c-a-b)+1}$.
\end{definition}

\begin{lemma}\label{lem:hyp-basic}
Assume $\Hphi$. Then
\begin{eqnarray}\label{eq:hyp-basic}
	{}_{2}F_{1}^{\mathcal{R}(p,q)}(a,b;c;w)&=&{}_{2}\phi_{1}\big(\fq^{a},\fq^{b};\fq^{c};\fq,\phi_{1}^{a+b-c-1}w\big),
\end{eqnarray}
where ${}_{2}\phi_{1}(A,B;C;\fq,z)=\sum_{n}(A;\fq)_{n}(B;\fq)_{n}z^{n}/((C;\fq)_{n}(\fq;\fq)_{n})$ is the basic hypergeometric series \cite{GasperRahman}. In particular \eqref{eq:hypergeometric} converges absolutely for $|w|<\phi_{1}^{\Re(c-a-b)+1}$.
\end{lemma}
{\it Proof.}
By Lemma \ref{lem:numbers}(i), $\Pch{\pn{a}}{n}=\prod_{j=0}^{n-1}\phi_{1}^{a+j-1}[a+j]_{\fq}=\phi_{1}^{n(a-1)+n(n-1)/2}(\fq^{a};\fq)_{n}(1-\fq)^{-n}$ and $\pf{n}=\phi_{1}^{n(n-1)/2}(\fq;\fq)_{n}(1-\fq)^{-n}$. In the coefficient of \eqref{eq:hypergeometric} the factors $\phi_{1}^{n(n-1)/2}$ and $(1-\fq)^{-n}$ cancel, and $\phi_{1}^{n(a-1)+n(b-1)-n(c-1)}=\phi_{1}^{n(a+b-c-1)}$ remains. The ${}_{2}\phi_{1}$ series converges for $|z|<1$ \cite{GasperRahman}. \hfill\(\Box\)

\begin{theorem}[Deformed Euler integral]\label{thm:euler}
Assume $\Hphi$ and let $a,b,c\in\mathbb{C}$ with $\Re b>0$ and $\Re(c-b)>0$. Then, for $|w|<\phi_{1}^{\Re(c-a-b)+1}$,
\begin{eqnarray}\label{eq:euler}
	{}_{2}F_{1}^{\mathcal{R}(p,q)}(a,b;c;w)&=&\frac{\phi_{1}^{\frac{a(a-1)}{2}}}{\BR(b,c-b)}\int_{0}^{1}
	\frac{t^{b-1}\,(1\ominus\phi_{2}t)^{c-b-1}_{\mathcal{R}(p,q)}}{(1\ominus\phi_{1}^{a+b-c}w\,t)^{a}_{\mathcal{R}(p,q)}}\,\dR t ,
\end{eqnarray}
the deformed powers being those of Definition \ref{def:powers}.
\end{theorem}
{\it Proof.}
Put $z=\phi_{1}^{a+b-c-1}w$, so that $|z|<1$ and $\phi_{1}^{a+b-c}wt=\phi_{1}zt$. At the nodes $t_{j}=\fq^{j}/\phi_{1}$ one has $\phi_{1}zt_{j}=z\fq^{j}$, and by \eqref{eq:powers-complex} and the $\fq-$binomial theorem \eqref{eq:binomial},
\begin{eqnarray}
	\frac{1}{(1\ominus z\fq^{j})^{a}_{\mathcal{R}(p,q)}}&=&\phi_{1}^{-\frac{a(a-1)}{2}}\,\frac{(z\fq^{j+a};\fq)_{\infty}}{(z\fq^{j};\fq)_{\infty}}
	\;=\;\phi_{1}^{-\frac{a(a-1)}{2}}\sum_{n\geq0}\frac{(\fq^{a};\fq)_{n}}{(\fq;\fq)_{n}}\,(\phi_{1}zt_{j})^{n}.\nn
\end{eqnarray}
Since $|z|<1$ and the integrand of $\BR(b+n,c-b)$ is dominated by that of $\BR(\Re b,\Re(c-b))$, the double series converges absolutely and the summation may be interchanged with the deformed integral. Hence the integral in \eqref{eq:euler} equals
\begin{eqnarray}
	\phi_{1}^{-\frac{a(a-1)}{2}}\sum_{n\geq0}\frac{(\fq^{a};\fq)_{n}}{(\fq;\fq)_{n}}\,\phi_{1}^{n}z^{n}\,\BR(b+n,c-b).\nn
\end{eqnarray}
By \eqref{eq:finite-x} with $(x,y)=(b,c-b)$, $\BR(b+n,c-b)=\phi_{1}^{n(c-b-1)}\Pch{\pn{b}}{n}\Pch{\pn{c}}{n}^{-1}\BR(b,c-b)$, and by the proof of Lemma \ref{lem:hyp-basic}, $\Pch{\pn{b}}{n}/\Pch{\pn{c}}{n}=\phi_{1}^{n(b-c)}(\fq^{b};\fq)_{n}/(\fq^{c};\fq)_{n}$. The total power of $\phi_{1}$ is $n+n(c-b-1)+n(b-c)=0$, so that the integral equals $\phi_{1}^{-a(a-1)/2}\BR(b,c-b)\,{}_{2}\phi_{1}(\fq^{a},\fq^{b};\fq^{c};\fq,z)$, and Lemma \ref{lem:hyp-basic} concludes. \hfill\(\Box\)

\begin{corollary}[Deformed Gauss summation]\label{cor:gauss}
Assume $\Hphi$ and let $a,b,c\in\mathbb{C}$ with $\Re b>0$, $\Re(c-a-b)>0$ and $\pn{c+j}\neq0$ for $j\geq0$. Then
\begin{eqnarray}\label{eq:gauss}
	{}_{2}F_{1}^{\mathcal{R}(p,q)}\big(a,b;c;\phi_{1}\phi_{2}^{\,c-a-b}\big)&=&\phi_{1}^{-ab}\,\frac{\GR(c)\,\GR(c-a-b)}{\GR(c-a)\,\GR(c-b)} .
\end{eqnarray}
\end{corollary}
{\it Proof.}
Assume first $\Re(c-b)>0$. With $w=\phi_{1}\phi_{2}^{c-a-b}$ one has $z=\phi_{1}^{a+b-c-1}w=\fq^{c-a-b}$, $|z|<1$, and $\phi_{1}zt=\phi_{2}\fq^{c-a-b-1}t$. At the nodes, the numerator $(\phi_{1}z\fq^{a}t;\fq)_{\infty}=(\phi_{2}\fq^{c-b-1}t;\fq)_{\infty}$ of $1/(1\ominus\phi_{1}zt)^{a}_{\mathcal{R}(p,q)}$ cancels the denominator of $(1\ominus\phi_{2}t)^{c-b-1}_{\mathcal{R}(p,q)}$, and the integrand of \eqref{eq:euler} becomes
\begin{eqnarray}
	\phi_{1}^{\frac{(c-b-1)(c-b-2)}{2}-\frac{a(a-1)}{2}}\,t^{b-1}\,\frac{(\phi_{2}t;\fq)_{\infty}}{(\phi_{2}\fq^{c-a-b-1}t;\fq)_{\infty}}
	&=&\phi_{1}^{\,a(c-a-b-1)}\,t^{b-1}(1\ominus\phi_{2}t)^{c-a-b-1}_{\mathcal{R}(p,q)},\nn
\end{eqnarray}
by the exponent identity $(u-1)(u-2)-a(a-1)-(u-a-1)(u-a-2)=2a(u-a-1)$ with $u=c-b$. Therefore
\begin{eqnarray}
	{}_{2}F_{1}^{\mathcal{R}(p,q)}\big(a,b;c;\phi_{1}\phi_{2}^{\,c-a-b}\big)&=&\phi_{1}^{\frac{a(a-1)}{2}+a(c-a-b-1)}\,\frac{\BR(b,c-a-b)}{\BR(b,c-b)} .\nn
\end{eqnarray}
By Theorem \ref{thm:beta1} the quotient of beta functions equals $\phi_{1}^{E(b,c-a-b)-E(b,c-b)}\GR(c)\GR(c-a-b)/(\GR(c-a)\GR(c-b))$, and $E(b,c-a-b)-E(b,c-b)=-a(2c-a-3)/2$ (Appendix \ref{app:hyp}). The total exponent is $a(a-1)/2+a(c-a-b-1)-a(2c-a-3)/2=-ab$. Both sides of \eqref{eq:gauss} are analytic in $b$ in the region where the series converges, which removes the auxiliary condition $\Re(c-b)>0$. \hfill\(\Box\)

\begin{corollary}[Deformed Chu Vandermonde identity]\label{cor:chu}
Assume $\Hphi$. For $n\in\mathbb{N}$ and $b,c\in\mathbb{C}$ with $\pn{c+j}\neq0$ and $\pn{c-b+j}\neq0$ for $0\leq j<n$,
\begin{eqnarray}\label{eq:chu}
	\sum_{k=0}^{n}\frac{\Pch{\pn{-n}}{k}\,\Pch{\pn{b}}{k}}{\Pch{\pn{c}}{k}\,\pf{k}}\,\big(\phi_{1}\phi_{2}^{\,c-b+n}\big)^{k}
	&=&\phi_{1}^{\,nb}\,\frac{\Pch{\pn{c-b}}{n}}{\Pch{\pn{c}}{n}} .
\end{eqnarray}
\end{corollary}
{\it Proof.}
For $a=-n$ the series \eqref{eq:hypergeometric} terminates, since $\pn{0}=0$, and both sides of \eqref{eq:gauss} are rational in $\phi_{1}^{b},\phi_{2}^{b}$, so that the analytic continuation in $b$ is immediate. One uses $\GR(c)/\GR(c+n)=1/\Pch{\pn{c}}{n}$ and $\GR(c-b+n)/\GR(c-b)=\Pch{\pn{c-b}}{n}$. \hfill\(\Box\)

\begin{remark}\label{rem:euler}
For $\phi_{1}=1$, Theorem \ref{thm:euler} and Corollary \ref{cor:gauss} reduce to the Euler integral representation of ${}_{2}\phi_{1}$ and to the $q-$Gauss summation \cite[Sections 1.5 and 1.11]{GasperRahman}; letting further $\fq\to1$ one recovers the classical Euler integral and the Gauss summation formula. The argument $\phi_{1}\phi_{2}^{c-a-b}$ in \eqref{eq:gauss} plays the role of the classical point $w=1$; the valuation factor $\phi_{1}^{-ab}$ is symmetric in $a$ and $b$, as it must be, since ${}_{2}F_{1}^{\mathcal{R}(p,q)}$ is symmetric in its two numerator parameters. The proof of Corollary \ref{cor:gauss} derives the deformed Gauss summation from the first kind beta gamma relation alone, the only inputs being the $\fq-$binomial theorem and the finite shift identity \eqref{eq:finite-x}.
\end{remark}

\subsection{The incomplete deformed beta function}\label{subsec:incomplete}

\begin{definition}\label{def:incomplete}
Assume $\Hphi$. For $0<u\leq1$ and $(x,y)\in\Omega_{+}$, the incomplete first kind deformed beta function is
\begin{eqnarray}\label{eq:incomplete-def}
	\BR(u;x,y)&:=&\int_{0}^{u}t^{x-1}\,(1\ominus\phi_{2}t)^{y-1}_{\mathcal{R}(p,q)}\,\dR t ,
\end{eqnarray}
so that $\BR(1;x,y)=\BR(x,y)$.
\end{definition}

\begin{theorem}\label{thm:incomplete}
Assume $\Hphi$, let $0<u\leq1$ and $(x,y)\in\Omega_{+}$. Then
\begin{eqnarray}\label{eq:incomplete}
	\BR(u;x,y)&=&\phi_{1}^{\frac{(y-1)(y-2)}{2}}\,\frac{u^{x}}{\pn{x}}\;{}_{2}F_{1}^{\mathcal{R}(p,q)}\big(1-y,x;x+1;\phi_{1}\phi_{2}^{\,y}u\big),
\end{eqnarray}
and the integration by parts recurrence
\begin{eqnarray}\label{eq:incomplete-ibp}
	\BR(u;x+1,y)&=&\phi_{1}^{-x}\,\frac{\pn{x}}{\pn{y}}\,\BR(u;x,y+1)-\frac{u^{x}}{\phi_{1}^{x}\,\pn{y}}\,(1\ominus u)^{y}_{\mathcal{R}(p,q)}
\end{eqnarray}
holds. For $u=1$ the boundary term vanishes, \eqref{eq:incomplete-ibp} reduces to \eqref{eq:ibp-beta}, and \eqref{eq:incomplete} together with the deformed Gauss summation \eqref{eq:gauss} gives back the beta gamma relation \eqref{eq:beta1-gamma}.
\end{theorem}
{\it Proof.}
By \eqref{eq:powers-complex} and the $\fq-$binomial theorem with $b=\fq^{1-y}$ and $w=\phi_{2}\fq^{y-1}t$,
\begin{eqnarray}
	(1\ominus\phi_{2}t)^{y-1}_{\mathcal{R}(p,q)}&=&\phi_{1}^{\frac{(y-1)(y-2)}{2}}\,\frac{(\phi_{2}t;\fq)_{\infty}}{(\phi_{2}\fq^{y-1}t;\fq)_{\infty}}
	\;=\;\phi_{1}^{\frac{(y-1)(y-2)}{2}}\sum_{n\geq0}\frac{(\fq^{1-y};\fq)_{n}}{(\fq;\fq)_{n}}\,\big(\phi_{2}\fq^{y-1}t\big)^{n},\nn
\end{eqnarray}
for $|\phi_{2}\fq^{y-1}t|<1$. At the nodes $t=u\fq^{k}/\phi_{1}$ of \eqref{eq:incomplete-def} one has $|\phi_{2}\fq^{y-1}t|\leq|\fq^{y}|u<1$, so that the series may be integrated termwise. By Proposition \ref{prop:integral}(ii), $\int_{0}^{u}t^{x+n-1}\dR t=u^{x+n}/\pn{x+n}$, and by Lemma \ref{lem:numbers}(i),
$1/\pn{x+n}=\phi_{1}^{-n}\pn{x}^{-1}(\fq^{x};\fq)_{n}/(\fq^{x+1};\fq)_{n}$. Hence
\begin{eqnarray}
	\BR(u;x,y)&=&\phi_{1}^{\frac{(y-1)(y-2)}{2}}\,\frac{u^{x}}{\pn{x}}\sum_{n\geq0}\frac{(\fq^{1-y};\fq)_{n}(\fq^{x};\fq)_{n}}{(\fq;\fq)_{n}(\fq^{x+1};\fq)_{n}}\,\big(\fq^{y}u\big)^{n},\nn
\end{eqnarray}
since $\phi_{2}\fq^{y-1}/\phi_{1}=\fq^{y}$. By Lemma \ref{lem:hyp-basic} with $(a,b,c)=(1-y,x,x+1)$, the last series is ${}_{2}F_{1}^{\mathcal{R}(p,q)}(1-y,x;x+1;w)$ with $\phi_{1}^{a+b-c-1}w=\phi_{1}^{-y-1}w=\fq^{y}u$, that is $w=\phi_{1}\phi_{2}^{y}u$. For \eqref{eq:incomplete-ibp} we apply \eqref{eq:ibp} on $[0,u]$ with $f(t)=(t/\phi_{1})^{x}$ and $g(t)=-(1\ominus t)^{y}_{\mathcal{R}(p,q)}/\pn{y}$ exactly as in the proof of Lemma \ref{lem:ibp-pascal}; the boundary term $f(u)g(u)$ no longer vanishes. Finally, at $u=1$, $(a,b,c)=(1-y,x,x+1)$ satisfies $c-a-b=y$, so that $\phi_{1}\phi_{2}^{y}$ is the Gauss point of Corollary \ref{cor:gauss}, and
\begin{eqnarray}
	{}_{2}F_{1}^{\mathcal{R}(p,q)}\big(1-y,x;x+1;\phi_{1}\phi_{2}^{\,y}\big)&=&\phi_{1}^{x(y-1)}\,\frac{\GR(x+1)\GR(y)}{\GR(x+y)\GR(1)}\;=\;\phi_{1}^{x(y-1)}\pn{x}\,\frac{\GR(x)\GR(y)}{\GR(x+y)};\nn
\end{eqnarray}
since $(y-1)(y-2)/2+x(y-1)=E(x,y)$, \eqref{eq:incomplete} at $u=1$ is \eqref{eq:beta1-gamma}. \hfill\(\Box\)

Formula \eqref{eq:incomplete} is the deformed analogue of the classical representation $\mathrm{B}(u;x,y)=x^{-1}u^{x}\,{}_{2}F_{1}(1-y,x;x+1;u)$, to which it reduces when $\phi_{1}\to1$ and $\fq\to1$, the argument $\phi_{1}\phi_{2}^{y}u$ tending to $u$. The two derivations of the beta gamma relation, through the $\fq-$binomial theorem in Theorem \ref{thm:beta1} and through the deformed Gauss summation here, provide a consistency check of the valuation exponent $E(x,y)$.

\section{Deformed beta distributions}\label{sec:distributions}

The positivity of the two deformed beta integrands on the lattice of nodes provides two families of discrete probability distributions, in the spirit of the $\mathcal{R}(p,q)-$analogues of discrete distributions of \cite{Mahouton_Melong_Rpq_analogues_distribution, HounkonnouMelong2023}.

\begin{definition}\label{def:distributions}
Assume $\Hphi$ and let $x,y>0$.
\begin{enumerate}
\item[\textnormal{(i)}] The first kind deformed beta distribution $\mu^{(1)}_{x,y}$ is the probability measure on $\{t_{j}=\fq^{j}/\phi_{1}:j\geq0\}\subset(0,1/\phi_{1}]$ given by
\begin{eqnarray}
	\mu^{(1)}_{x,y}(\{t_{j}\})&:=&\frac{(1-\fq)\,\fq^{j}\,t_{j}^{x-1}\,(1\ominus\phi_{2}t_{j})^{y-1}_{\mathcal{R}(p,q)}}{\BR(x,y)} .\nn
\end{eqnarray}
\item[\textnormal{(ii)}] The second kind deformed beta distribution $\mu^{(2)}_{x,y}$ is the probability measure on $\{t_{k}=\fq^{k}/\phi_{1}:k\in\mathbb{Z}\}\subset(0,\infty)$ given by
\begin{eqnarray}
	\mu^{(2)}_{x,y}(\{t_{k}\})&:=&\frac{(1-\fq)\,\fq^{k}\,t_{k}^{x-1}}{(1\oplus\phi_{1}t_{k})^{x+y}_{\mathcal{R}(p,q)}\,\BB(x,y)} .\nn
\end{eqnarray}
\end{enumerate}
\end{definition}

By Lemma \ref{lem:powers}(v) and Lemmas \ref{lem:beta1-series} and \ref{lem:beta2-series}, both measures are well defined, their total mass being $1$ by the definition of the deformed integrals. For a random variable $T$ with law $\mu$ we write $\mathbb{E}_{\mu}[T^{k}]=\sum_{t}t^{k}\mu(\{t\})$.

\begin{theorem}\label{thm:moments1}
Assume $\Hphi$, let $x,y>0$ and let $T$ have law $\mu^{(1)}_{x,y}$. Then, for $k\in\mathbb{N}$,
\begin{eqnarray}\label{eq:moments1}
	\mathbb{E}[T^{k}]&=&\frac{\BR(x+k,y)}{\BR(x,y)}\;=\;\phi_{1}^{\,k(y-1)}\,\frac{\Pch{\pn{x}}{k}}{\Pch{\pn{x+y}}{k}},
\end{eqnarray}
and
\begin{eqnarray}\label{eq:mean-var1}
	\mathbb{E}[T]&=&\phi_{1}^{\,y-1}\,\frac{\pn{x}}{\pn{x+y}},\qquad
	\mathrm{Var}[T]\;=\;\phi_{1}^{\,2(y-1)}\,\frac{(\phi_{1}\phi_{2})^{x}\,\pn{x}\,\pn{y}}{\pn{x+y}^{2}\,\pn{x+y+1}} .
\end{eqnarray}
Moreover the moment sequence determines $\mu^{(1)}_{x,y}$, for $s\in\mathbb{C}$ with $\Re(x+s)>0$ one has $\mathbb{E}[T^{s}]=\BR(x+s,y)/\BR(x,y)$, and the distribution function at the nodes is
\begin{eqnarray}\label{eq:cdf1}
	\mathbb{P}\Big(T\leq\frac{\fq^{m}}{\phi_{1}}\Big)&=&\frac{\BR(\fq^{m};x,y)}{\BR(x,y)},\qquad m\in\mathbb{N}\cup\{0\},
\end{eqnarray}
with the incomplete deformed beta function of Theorem \ref{thm:incomplete}.
\end{theorem}
{\it Proof.}
Since $t^{k}\cdot t^{x-1}=t^{x+k-1}$, the moments are quotients of deformed beta functions, and \eqref{eq:moments1} is \eqref{eq:finite-x}. For the variance,
\begin{eqnarray}
	\mathrm{Var}[T]&=&\phi_{1}^{2(y-1)}\frac{\pn{x}}{\pn{x+y}}\Big(\frac{\pn{x+1}}{\pn{x+y+1}}-\frac{\pn{x}}{\pn{x+y}}\Big)\nn\\
	&=&\phi_{1}^{2(y-1)}\frac{\pn{x}\big(\pn{x+1}\pn{x+y}-\pn{x}\pn{x+y+1}\big)}{\pn{x+y}^{2}\pn{x+y+1}},\nn
\end{eqnarray}
and the deformed Tur\'an identity \eqref{eq:turan-identity} gives \eqref{eq:mean-var1}. The measure has bounded support, so that it is determined by its moments. The formula for $\mathbb{E}[T^{s}]$ is the definition of $\BR(x+s,y)$. Finally, the nodes $t_{j}$ with $t_{j}\leq\fq^{m}/\phi_{1}$ are those with $j\geq m$, and $\sum_{j\geq m}(1-\fq)\fq^{j}F(\fq^{j}/\phi_{1})=\int_{0}^{\fq^{m}}F\,\dR t$ by \eqref{eq:Rint} with $a=\fq^{m}$, which gives \eqref{eq:cdf1}. \hfill\(\Box\)

\begin{theorem}\label{thm:moments2}
Assume $\Hphi$, let $x,y>0$ and let $T$ have law $\mu^{(2)}_{x,y}$. Then $\mathbb{E}[T^{k}]<\infty$ if and only if $k<y$, and for $k\in\mathbb{N}$ with $k<y$,
\begin{eqnarray}\label{eq:moments2}
	\mathbb{E}[T^{k}]&=&\frac{\BB(x+k,y-k)}{\BB(x,y)}\;=\;\phi_{1}^{\frac{k(2y-k-3)}{2}}\,\phi_{2}^{-kx-\frac{k(k-1)}{2}}\,\frac{\Pch{\pn{x}}{k}}{\Pch{\pn{y-k}}{k}} .
\end{eqnarray}
In particular, for $y>1$, $\mathbb{E}[T]=\phi_{1}^{y-2}\phi_{2}^{-x}\pn{x}/\pn{y-1}$.
\end{theorem}
{\it Proof.}
The identity $t^{k}t^{x-1}/(1\oplus\phi_{1}t)^{x+y}=t^{x+k-1}/(1\oplus\phi_{1}t)^{(x+k)+(y-k)}$ shows that $\mathbb{E}[T^{k}]=\BB(x+k,y-k)/\BB(x,y)$, the series converging if and only if $y-k>0$ by the estimates of Lemma \ref{lem:beta2-series}, which are sharp for real arguments. By \eqref{eq:beta2-gamma1},
\begin{eqnarray}
	\frac{\BB(x+k,y-k)}{\BB(x,y)}&=&\phi_{1}^{-\frac{(x+k)^{2}+(x+k)+(y-k)^{2}-(y-k)-x^{2}-x-y^{2}+y}{2}}\,\frac{\Kq(x+k)}{\Kq(x)}\,\frac{\GR(x+k)\,\GR(y-k)}{\GR(x)\,\GR(y)},\nn
\end{eqnarray}
the denominators $\GR(x+y)$ cancelling. The exponent equals $-k(x-y)-k^{2}-k$, Proposition \ref{prop:Omega}(i) gives $\Kq(x+k)/\Kq(x)=\fq^{-kx-k(k-1)/2}$, and $\GR(x+k)/\GR(x)=\Pch{\pn{x}}{k}$, $\GR(y-k)/\GR(y)=1/\Pch{\pn{y-k}}{k}$. Writing $\fq=\phi_{2}/\phi_{1}$, the total power of $\phi_{1}$ is $-kx+ky-k^{2}-k+kx+k(k-1)/2=k(2y-k-3)/2$. \hfill\(\Box\)

\begin{figure}[ht]
\centering
\includegraphics[width=\linewidth]{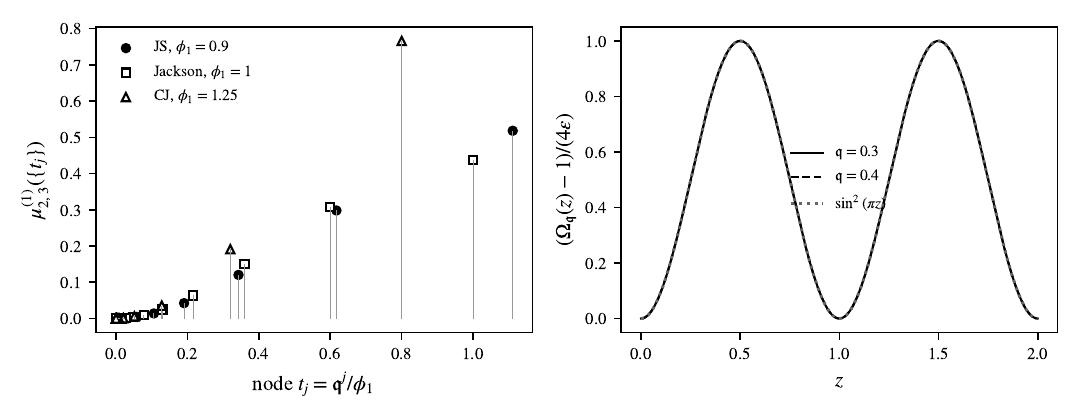}
\caption{Left: the first kind deformed beta distribution $\mu^{(1)}_{2,3}$ for the Jagannathan Srinivasa Rao deformation $(p,q)=(0.9,0.5)$, the Jackson case $q=0.6$ and the Chakrabarti Jagannathan deformation $(p,q)=(0.8,0.5)$; the atoms lie on the lattices $\fq^{j}/\phi_{1}$. Right: the normalised theta quotient $(\Omq(z)-1)/(4\varepsilon)$ of Proposition \ref{prop:Omega} for $\fq=0.3$ and $\fq=0.4$, indistinguishable at this scale from $\sin^{2}(\pi z)$.}\label{fig:distribution-theta}
\end{figure}

Figure \ref{fig:distribution-theta} shows how the deformation moves the atoms of $\mu^{(1)}_{x,y}$: for $\phi_{1}>1$ the support is contained in $(0,1/\phi_{1}]\subset(0,1)$ and the mass concentrates on the largest atom, while for $\phi_{1}<1$ the support extends beyond $1$.

\begin{remark}\label{rem:distributions}
In the undeformed limit, $\mu^{(1)}_{x,y}$ and $\mu^{(2)}_{x,y}$ converge to the beta distribution on $[0,1]$ and to the beta prime distribution on $(0,\infty)$, with means $x/(x+y)$ and $x/(y-1)$. The variance in \eqref{eq:mean-var1} is positive precisely because of the deformed Tur\'an identity, which thus appears as the deformed counterpart of the classical identity $(x+1)(x+y)-x(x+y+1)=y$.
\end{remark}

\section{Gauge invariance and particular cases}\label{sec:cases}

\subsection{Gauge invariance}\label{subsec:gauge}

We now consider structure functions satisfying the gauged hypothesis
\begin{eqnarray}\label{eq:Hphic}
	\textnormal{(H}^{c}_{\phi}\textnormal{)}\qquad\mathcal{R}(p^{n},q^{n})&=&c\,\frac{\phi_{1}^{n}-\phi_{2}^{n}}{\phi_{1}-\phi_{2}},\qquad n\in\mathbb{N}\cup\{0\},
\end{eqnarray}
with $c=c(p,q)>0$ and $\phi_{1}=\phi_{1}(p,q)>\phi_{2}=\phi_{2}(p,q)>0$, the gauge $c$ and the deformation functions being, like $\mathcal{R}(p,q)$, meromorphic in $(p,q)$. We write $[x]_{\phi}=(\phi_{1}^{x}-\phi_{2}^{x})/(\phi_{1}-\phi_{2})$, so that $\pn{x}:=c\,[x]_{\phi}$, and we denote with the index $\phi$ the objects of Sections \ref{sec:prelim} to \ref{sec:convexity} built from $[x]_{\phi}$. By \eqref{eq:diag}, $\partial_{\mathcal{R}(p,q)}=c\,\partial_{\phi}$, and the inverse of $\partial_{\mathcal{R}(p,q)}$ in the sense of Proposition \ref{prop:integral}(iv) is $\int\cdot\,\dR t:=c^{-1}\int\cdot\,d_{\phi}t$. The deformed powers \eqref{eq:powers-int} and \eqref{eq:powers-complex} only involve $\phi_{1},\phi_{2}$ and are unchanged, while $\ER(z)=E_{\phi}(z/c)$ and $\eR(z)=e_{\phi}(z/c)$. We define $\BR$ and $\BB$ by \eqref{eq:beta1-def} and \eqref{eq:beta2-def} with the integral $\dR t$, and
\begin{eqnarray}\label{eq:gamma-gauge}
	\GR(z)&:=&c^{z}\phi_{1}^{\frac{z(z-1)}{2}}\int_{0}^{a_{\mathcal{R}(p,q)}}t^{z-1}\ER(-c\,\phi_{2}t)\,\dR t,\nn\\
	\gR(z)&:=&c^{z}\phi_{2}^{\frac{z(z-1)}{2}}\int_{0}^{\infty;a_{\mathcal{R}(p,q)}}t^{z-1}\eR(-c\,\phi_{1}t)\,\dR t .
\end{eqnarray}

\begin{theorem}[Gauge invariance]\label{thm:gauge}
Assume \textnormal{(H}$^{c}_{\phi}$\textnormal{)}. Then $\GR=c^{z-1}\Gamma_{\phi}$, $\gR=c^{z-1}\gamma_{\phi}$, $\BR=c^{-1}\mathfrak{B}_{\phi}$ and $\BB=c^{-1}\mathcal{B}_{\phi}$. Consequently every identity of Theorems \ref{thm:gamma1}, \ref{thm:BM}, \ref{thm:logconvex-gamma}, \ref{thm:gamma2}, \ref{thm:wendel}, \ref{thm:beta1}, \ref{thm:exchange1}, \ref{thm:unit-shifts}, \ref{thm:shift-general}, \ref{thm:sym1}, \ref{thm:beta1-asymptotics}, \ref{thm:beta2}, \ref{thm:exchange2}, \ref{thm:beta2-shifts}, \ref{thm:beta2-sym}, \ref{thm:convex-beta}, \ref{thm:euler}, \ref{thm:incomplete}, \ref{thm:moments1} and \ref{thm:moments2} and of Corollaries \ref{cor:turan}, \ref{cor:gauss} and \ref{cor:chu} which is expressed through the deformed numbers $\pn{x}=c[x]_{\phi}$, the deformed factorials $\pf{n}=c^{n}[n]_{\phi}!$, the deformed Pochhammer symbols and the functions $\GR,\gR,\BR,\BB$ holds verbatim; the expressions through $\Gq$ and $\Bq$, such as \eqref{eq:gamma1-closed}, \eqref{eq:gamma2-closed}, \eqref{eq:beta2-q} and the first form of \eqref{eq:asym-x}, acquire the factor $c^{z-1}$ for the gamma functions and $c^{-1}$ for the beta functions. Theorem \ref{thm:multiplication} and Proposition \ref{prop:reflection} are stated for $c=1$.
\end{theorem}
{\it Proof.}
Since $\dR t=c^{-1}d_{\phi}t$ and $\ER(-c\phi_{2}t)=E_{\phi}(-\phi_{2}t)$, \eqref{eq:gamma-gauge} gives $\GR(z)=c^{z-1}\Gamma_{\phi}(z)$, and similarly for $\gR$; the beta integrals are multiplied by $c^{-1}$. In the beta gamma relations, the factors $c$ cancel because
\begin{eqnarray}
	\frac{\GR(x)\GR(y)}{\GR(x+y)}&=&c^{(x-1)+(y-1)-(x+y-1)}\,\frac{\Gamma_{\phi}(x)\Gamma_{\phi}(y)}{\Gamma_{\phi}(x+y)}\;=\;c^{-1}\,\frac{\Gamma_{\phi}(x)\Gamma_{\phi}(y)}{\Gamma_{\phi}(x+y)},\nn
\end{eqnarray}
and similarly with $\gR$. The same cancellation occurs in the Gauss summation \eqref{eq:gauss}, where the exponents of $c$ add up to $(c-1)+(c-a-b-1)-(c-a-1)-(c-b-1)=0$, and in the coefficients of ${}_{2}F_{1}^{\mathcal{R}(p,q)}$, which contain $2n$ deformed numbers in the numerator and in the denominator. In the shift identities and in the moments of the deformed beta distributions, the quotients $\Pch{\pn{x}}{m}\Pch{\pn{y}}{n}/\Pch{\pn{x+y}}{m+n}$ contain $m+n$ deformed numbers in the numerator and in the denominator, so that $c$ cancels. The functional equation becomes $\GR(z+1)=c^{z}\Gamma_{\phi}(z+1)=c[z]_{\phi}\,c^{z-1}\Gamma_{\phi}(z)=\pn{z}\GR(z)$. In Theorem \ref{thm:BM}, $f(x+1)=\pn{x}f(x)$ if and only if $g(x)=c^{1-x}f(x)$ satisfies $g(x+1)=[x]_{\phi}g(x)$, and $c^{1-x}$ is logarithmically affine, so that the convexity conditions are unchanged; for the same reason the second logarithmic derivatives of $\GR$ and $\Gamma_{\phi}$ coincide. \hfill\(\Box\)

\subsection{Particular cases}\label{subsec:cases}

We list the cases of Table \ref{tab:examples}; in each of them the deformed numbers and $(\phi_{1},\phi_{2},c)$ are read off from the structure function, following the classification of \cite{HounkBuk, Mahouton_Kyemba_Rpq_calculs}.

\paragraph{(i) Jagannathan Srinivasa Rao deformation.}
For $\mathcal{R}(u,v)=(u-v)/(p-q)$ with $0<q<p$ one has $\pn{n}=[n]_{p,q}=(p^{n}-q^{n})/(p-q)$, $\phi_{1}=p$, $\phi_{2}=q$, $c=1$ and $\fq=q/p$. Theorem \ref{thm:beta1} gives, for all $(x,y)\in\Omega_{+}$,
\begin{eqnarray}\label{eq:case-JS}
	\int_{0}^{1}t^{x-1}(1\ominus qt)^{y-1}_{p,q}\,d_{p,q}t&=&p^{\frac{(y-1)(2x+y-2)}{2}}\,\frac{\Gamma_{p,q}(x)\,\Gamma_{p,q}(y)}{\Gamma_{p,q}(x+y)},\qquad
	\Gamma_{p,q}(z)\;=\;p^{\frac{(z-1)(z-2)}{2}}\Gamma_{q/p}(z),
\end{eqnarray}
which extends to complex arguments the relation proved for positive integers in \cite{ref55}. The exchange defect is $p^{(y-x)(x+y-1)/2}$, and Theorem \ref{thm:logconvex-gamma} shows that $\Gamma_{p,q}$ is logarithmically convex on $(0,\infty)$ if and only if $p\geq1$.

For the second kind function, Theorem \ref{thm:beta2} gives
\begin{eqnarray}\label{eq:case-JS-2}
	\int_{0}^{\infty;1}\frac{t^{x-1}\,d_{p,q}t}{(1\oplus pt)^{x+y}_{p,q}}&=&p^{-x-\frac{y(y-1)}{2}}\,q^{-\frac{x(x-1)}{2}}\,\frac{\gamma_{p,q}(x)\,\Gamma_{p,q}(y)}{\Gamma_{p,q}(x+y)},\qquad
	\gamma_{p,q}\;=\;\Omega_{q/p}\,\Gamma_{p,q},
\end{eqnarray}
and its exchange defect is $p^{(y-x)(3-x-y)/2}q^{(y-x)(x+y-1)/2}\,\Omega_{q/p}(x)/\Omega_{q/p}(y)$.

\paragraph{(ii) Jackson case.}
For $p=1$ in (i), $\phi_{1}=1$ and $\phi_{2}=\fq=q$. All valuation factors of the first kind disappear, $\BR=\Bq$ is commutative and $\GR=\Gamma_{q}$, in accordance with \cite{GasperRahman, Askey1978}; the second kind function \eqref{eq:beta2-q} becomes $\Kq(x)\Bq(x,y)$ and remains noncommutative by Theorem \ref{thm:exchange2}. Letting further $q\to1$ restores the Euler beta function in both representations.

\paragraph{(iii) Chakrabarti Jagannathan deformation.}
For $\mathcal{R}(u,v)=(1-uv)/((p^{-1}-q)u)$ with $pq<1$, one has $\pn{n}=(p^{-n}-q^{n})/(p^{-1}-q)$, $\phi_{1}=p^{-1}$, $\phi_{2}=q$, $c=1$ and $\fq=pq$. Then
\begin{eqnarray}\label{eq:case-CJ}
	\BR(x,y)&=&p^{-\frac{(y-1)(2x+y-2)}{2}}\,\frac{\GR(x)\GR(y)}{\GR(x+y)},\qquad \GR(z)\;=\;p^{-\frac{(z-1)(z-2)}{2}}\Gamma_{pq}(z),
\end{eqnarray}
the exchange defect is $p^{-(y-x)(x+y-1)/2}$, and, if $p\leq1$, $\GR$ is logarithmically convex on $(0,\infty)$ by Theorem \ref{thm:logconvex-gamma}(i), while $x\mapsto\BR(x,y)$ is decreasing by Proposition \ref{prop:monotonicity}(i). The second kind function is $\BB(x,y)=p^{x+y(y-1)/2}q^{-x(x-1)/2}\gR(x)\GR(y)/\GR(x+y)$ with $\gR=\Omega_{pq}\GR$. Since $\phi_{1}>1$ when $p<1$, Theorem \ref{thm:beta1-asymptotics}(ii) shows that $\BR(x,y)$ grows like $p^{-y^{2}/2}$ as $y\to+\infty$.

\paragraph{(iv) Generalized $q-$Quesne deformation.}
For $\mathcal{R}(u,v)=(uv-1)/((q-p^{-1})v)$ with $pq\neq1$, one has
\begin{eqnarray}\label{eq:case-Quesne}
	\pn{n}&=&\frac{p^{n}-q^{-n}}{q-p^{-1}}\;=\;\frac{p}{q}\,\frac{\phi_{1}^{n}-\phi_{2}^{n}}{\phi_{1}-\phi_{2}},\qquad \{\phi_{1},\phi_{2}\}=\{p,q^{-1}\},
\end{eqnarray}
because $(p-q^{-1})/(q-p^{-1})=p/q$. Thus \textnormal{(H}$^{c}_{\phi}$\textnormal{)} holds with $c=p/q$, $\phi_{1}=\max\{p,q^{-1}\}$ and $\phi_{2}=\min\{p,q^{-1}\}$. If $pq<1$, then $\phi_{1}=q^{-1}$, $\phi_{2}=p$, $\fq=pq$, and by Theorem \ref{thm:gauge}
\begin{eqnarray}
	\BR(x,y)&=&q^{-\frac{(y-1)(2x+y-2)}{2}}\,\frac{\GR(x)\GR(y)}{\GR(x+y)},\qquad
	\GR(z)\;=\;\Big(\frac{p}{q}\Big)^{z-1}q^{-\frac{(z-1)(z-2)}{2}}\,\Gamma_{pq}(z);\nn
\end{eqnarray}
if $pq>1$, then $\phi_{1}=p$, $\phi_{2}=q^{-1}$, $\fq=(pq)^{-1}$ and the valuation factor is $p^{(y-1)(2x+y-2)/2}$.

\paragraph{(v) Multiparameter Hounkonnou Ngompe deformation.}
For $\mathcal{R}(u,v)=h(p,q)\,v^{\nu}u^{-\mu}(uv-1)/((q-p^{-1})v)$ \cite{HounkonnouNgompe, HmNe2}, put $\ell:=q^{\nu}p^{-\mu}$. Then
\begin{eqnarray}\label{eq:case-HN}
	\pn{n}&=&h\,\ell^{n}\,\frac{p^{n}-q^{-n}}{q-p^{-1}}\;=\;h\,\frac{p}{q}\,\ell\;\frac{(\ell p)^{n}-(\ell q^{-1})^{n}}{\ell p-\ell q^{-1}},
\end{eqnarray}
because $\ell^{n}[n]_{a,b}=\ell\,[n]_{\ell a,\ell b}$ with $[n]_{a,b}=(a^{n}-b^{n})/(a-b)$. Hence \textnormal{(H}$^{c}_{\phi}$\textnormal{)} holds with $c=h\,p^{1-\mu}q^{\nu-1}$ and $\{\phi_{1},\phi_{2}\}=\{p^{1-\mu}q^{\nu},\,p^{-\mu}q^{\nu-1}\}$, whose quotient is $pq$. In the usual range $0<pq<1$ of \cite{HounkonnouNgompe}, $\phi_{1}=p^{-\mu}q^{\nu-1}$, $\phi_{2}=p^{1-\mu}q^{\nu}$ and $\fq=pq$, so that
\begin{eqnarray}\label{eq:case-HN-beta}
	\BR(x,y)&=&\big(p^{-\mu}q^{\nu-1}\big)^{\frac{(y-1)(2x+y-2)}{2}}\,\frac{\GR(x)\GR(y)}{\GR(x+y)},\\
	\GR(z)\;&=&\;c^{z-1}\big(p^{-\mu}q^{\nu-1}\big)^{\frac{(z-1)(z-2)}{2}}\,\Gamma_{pq}(z),
\end{eqnarray}
with exchange defect $(p^{-\mu}q^{\nu-1})^{(y-x)(x+y-1)/2}$. The gauge constant $h$ and the factor $p^{1-\mu}q^{\nu-1}$ disappear from the beta gamma relation, in accordance with Theorem \ref{thm:gauge}, while the parameters $\mu$ and $\nu$ enter the valuation factor through $\phi_{1}$. Compared with the generalized $q-$Quesne case (iv) with the same $(p,q)$, for which $\phi_{1}=q^{-1}$, one has $\phi_{1}^{HN}=\ell\,\phi_{1}^{Q}$, $\phi_{2}^{HN}=\ell\,\phi_{2}^{Q}$ and $c^{HN}=h\ell\,c^{Q}$, hence
\begin{eqnarray}\label{eq:case-HN-vs-Q}
	\GR^{HN}(z)&=&h^{z-1}\ell^{\frac{z(z-1)}{2}}\,\GR^{Q}(z),\\
	\frac{\GR^{HN}(x)\GR^{HN}(y)}{\GR^{HN}(x+y)}\;&=&\;\frac{1}{h\,\ell^{xy}}\,\frac{\GR^{Q}(x)\GR^{Q}(y)}{\GR^{Q}(x+y)},\\
	\BR^{HN}(x,y)&=&h^{-1}\ell^{\frac{(y-1)(y-2)}{2}-x}\,\BR^{Q}(x,y).
\end{eqnarray}
The gamma quotients thus differ by the factor $h^{-1}(p^{\mu}/q^{\nu})^{xy}$, which is symmetric in $(x,y)$, while the beta functions themselves differ by a nonsymmetric factor, since the valuation factor $\phi_{1}^{(y-1)(2x+y-2)/2}$ is rescaled by $\ell^{(y-1)(2x+y-2)/2}$.

\begin{table}[ht]
\centering
\caption{Valuation and exchange factors of the first kind deformed beta function.}\label{tab:valuation}
\begin{tabular}{lllll}
\hline
Deformation & $\phi_{1}$ & $\fq$ & valuation factor & exchange defect\\
\hline
Jagannathan Srinivasa Rao & $p$ & $q/p$ & $p^{(y-1)(2x+y-2)/2}$ & $p^{(y-x)(x+y-1)/2}$\\
Jackson & $1$ & $q$ & $1$ & $1$\\
Chakrabarti Jagannathan & $p^{-1}$ & $pq$ & $p^{-(y-1)(2x+y-2)/2}$ & $p^{-(y-x)(x+y-1)/2}$\\
$q-$Quesne, $pq<1$ & $q^{-1}$ & $pq$ & $q^{-(y-1)(2x+y-2)/2}$ & $q^{-(y-x)(x+y-1)/2}$\\
Hounkonnou Ngompe, $pq<1$ & $p^{-\mu}q^{\nu-1}$ & $pq$ & $\phi_{1}^{(y-1)(2x+y-2)/2}$ & $\phi_{1}^{(y-x)(x+y-1)/2}$\\
\hline
\end{tabular}
\end{table}

\section{Numerical verification}\label{sec:numerics}

All identities of the paper were checked with the multiple precision library \texttt{mpmath}, working with forty significant digits. The deformed integrals were computed from their defining Jackson type series \eqref{eq:Rint} and \eqref{eq:Rint-bilateral}, truncated at $4000$ terms on $[0,1]$ and at $|k|\leq1500$ for the bilateral integrals, the deformed powers from \eqref{eq:powers-int} and \eqref{eq:powers-complex}, and the closed forms from $\Gq$ and the infinite products. The truncation errors are bounded by geometric tails of ratio $\fq^{\Re x}$ and $\fq^{\Re y}$, far below the working precision. Three deformations are used: Jagannathan Srinivasa Rao with $(p,q)=(0.9,0.5)$, so that $\phi_{1}=0.9$ and $\fq=5/9$; Jackson with $q=0.6$; Chakrabarti Jagannathan with $(p,q)=(0.8,0.5)$, so that $\phi_{1}=1.25$ and $\fq=0.4$. A gauged case $(\phi_{1},\phi_{2},c)=(2,0.8,0.4)$ of Quesne type was also tested.

Table \ref{tab:num1} compares the Jackson series of $\BR(x,y)$ with the closed form \eqref{eq:beta1-gamma} and the computed exchange quotient with the factor of Theorem \ref{thm:exchange1}. In all cases the relative discrepancy is below $10^{-30}$ and the exchange quotient agrees with $\phi_{1}^{(y-x)(x+y-1)/2}$ to all displayed digits.

\begin{table}[ht]
\centering
\caption{First kind deformed beta function: Jackson series versus closed form, and exchange quotient.}\label{tab:num1}
\begin{tabular}{llllll}
\hline
Case & $(x,y)$ & $\BR(x,y)$ & rel. error & $\BR(x,y)/\BR(y,x)$ & $\phi_{1}^{(y-x)(x+y-1)/2}$\\
\hline
JS & $(0.5,1.5)$ & $1.408035512089$ & $3\times10^{-31}$ & $0.9486832980505$ & $0.9486832980505$\\
JS & $(1.3,2.7)$ & $0.4097476007829$ & $7\times10^{-31}$ & $0.801510599147$ & $0.801510599147$\\
JS & $(2.7,1.3)$ & $0.5112191918847$ & $<10^{-31}$ & $1.247644137288$ & $1.247644137288$\\
Jackson & $(0.5,1.5)$ & $1.477639940991$ & $9\times10^{-31}$ & $1$ & $1$\\
Jackson & $(1.3,2.7)$ & $0.3962833511377$ & $4\times10^{-31}$ & $1$ & $1$\\
CJ & $(0.5,1.5)$ & $1.536420967390$ & $<10^{-31}$ & $1.118033988750$ & $1.118033988750$\\
CJ & $(1.3,2.7)$ & $0.5693285480355$ & $2\times10^{-31}$ & $1.597758097756$ & $1.597758097756$\\
CJ & $(2.7,1.3)$ & $0.3563296276422$ & $4\times10^{-31}$ & $0.6258769718675$ & $0.6258769718675$\\
\hline
\end{tabular}
\end{table}

For the Jagannathan Srinivasa Rao case at $(x,y)=(1.3,2.7)$, the relative residuals of \eqref{eq:shift-x}, \eqref{eq:shift-y}, \eqref{eq:shift-xy}, of \eqref{eq:finite-x}, \eqref{eq:finite-y}, \eqref{eq:finite-mixed} for $n=2,3$, of the general identity \eqref{eq:shift-general} for $(m,n)=(2,3)$ at $(x,y)=(1.3,0.7)$, of the integration by parts and Pascal relations \eqref{eq:ibp-beta}, \eqref{eq:pascal-beta} and of the symmetry of $\BRs$ are all between $3\times10^{-31}$ and $2\times10^{-30}$. The partial fraction identity \eqref{eq:partial-fractions} was checked for $n\in\{0,1,3,5\}$ and $x\in\{0.3,2.2\}$, the special value \eqref{eq:special2} and the Tur\'an identity \eqref{eq:turan-identity} with the same accuracy. The shift sums \eqref{eq:shift-sum1} and \eqref{eq:shift-sumn} hold for $n\in\{1,2,3\}$ with relative residuals below $5\times10^{-41}$ (forty digits), and the finite interval variant \eqref{eq:beta2-variant} for $n\in\{0,\dots,4\}$ and $x\in\{0.7,2.3\}$ below $2\times10^{-40}$; it is negative for odd $n$, for instance $\mathcal{B}^{[0,1]}_{\mathcal{R}(p,q)}(0.7,4)=-0.0265388545478134\ldots$. For the generalized $q-$Quesne structure function with $(p,q)=(0.8,0.9)$, which is in Regime II and satisfies \textnormal{(H}$^{c}_{\phi}$\textnormal{)} with $(\phi_{1},\phi_{2},c)=(1/0.9,0.8,8/9)$, the power series $I_{\mathcal{R}(p,q)}f(a)$ of Theorem \ref{thm:HB-Jackson} for $f(t)=(1-t/2)^{-1}$ and $a=0.9$, the Regime II series followed by the multiplier \eqref{eq:multiplier}, and $c^{-1}$ times the series \eqref{eq:Rint} in the variables $(\phi_{1},\phi_{2})$ all equal $1.3661999273378131404\ldots$, in accordance with Corollary \ref{cor:HB-Hphi}.

The results of Sections \ref{subsec:gamma-further}, \ref{subsec:beta1-asymptotics}, \ref{sec:euler} and \ref{sec:distributions} were tested in the Jagannathan Srinivasa Rao case. The multiplication formula \eqref{eq:multiplication} holds for $n\in\{2,3,4\}$ at $z=0.83$ and the reflection formula \eqref{eq:reflection} at $x\in\{0.3,0.71\}$ with relative errors below $4\times10^{-31}$. The Wendel type estimate \eqref{eq:wendel} holds at $x\in\{0.4,2.3,8,20\}$, $s\in\{0.25,0.6\}$, the quotient of $\GR(x+s)/\GR(x)$ by its upper bound being $0.806$, $0.980$, $0.99947$ and $0.9999995$ for $s=0.25$, in agreement with \eqref{eq:wendel-limit}. The quotient of $\BR(x,1.7)$ by $\phi_{1}^{(x-1)(y-1)}\GR(y)\pn{x}^{-y}$ equals $0.98473$, $0.99922$ and $0.9999978$ at $x=5,10,20$, as predicted by \eqref{eq:asym-x}. The deformed Euler integral \eqref{eq:euler} was checked for $(a,b,c)=(0.7,1.3,3.1)$ and $w\in\{0.2,-0.5\}$, the deformed Gauss summation \eqref{eq:gauss} for the same parameters, the deformed Chu Vandermonde identity \eqref{eq:chu} for $n\in\{2,4\}$, the moments \eqref{eq:moments2} for $(x,y)=(1.2,3.4)$ and $k\in\{1,2\}$, and the incomplete beta formulas \eqref{eq:incomplete} and \eqref{eq:incomplete-ibp} for $(x,y)=(1.3,2.6)$ and $u\in\{1,\fq^{2},0.4\}$ (for instance $\BR(0.4;1.3,2.6)=0.196105082621922\ldots$), all with relative errors below $10^{-30}$. The product identities \eqref{eq:product-swap} and \eqref{eq:product-assoc} hold at $(x,y,z)=(0.7,1.6,2.3)$ to $5\times10^{-31}$, the entries of Table \ref{tab:integer-values} agree with the Jackson series to the working precision (for instance $\BR(3,3)=0.2137354293\ldots$), and the minimum of Proposition \ref{prop:monotonicity}(ii) for $y=1.5$ is attained at $x_{\min}=3.59357\ldots$ with $\BR(x_{\min},1.5)=0.428001\ldots$.

Table \ref{tab:num2} concerns the second kind function. The bilateral series \eqref{eq:beta2-series} agrees with \eqref{eq:beta2-q} with relative error below $10^{-40}$ at $(x,y)\in\{(0.5,0.5),(1.3,2.7),(2,3)\}$ in all cases, including the gauged one; the table shows the exchange quotient, which is not equal to $1$ in the Jackson case, and the unit shifts \eqref{eq:beta2-unit}.

\begin{table}[ht]
\centering
\caption{Second kind deformed beta function at $(x,y)=(1.3,2.7)$.}\label{tab:num2}
\begin{tabular}{llll}
\hline
Quantity & JS $(0.9,0.5)$ & Jackson $q=0.6$ & prediction\\
\hline
$\BB(x,y)$ & $0.9205877849376$ & $0.4377904998887$ & \eqref{eq:beta2-q}\\
$\BB(y,x)$ & $3.666046384441$ & $1.279819453615$ & \eqref{eq:beta2-q}\\
$\BB(x,y)/\BB(y,x)$ & $0.2511118759557$ & $0.3420720779420$ & \eqref{eq:exchange2}, same digits\\
$\BB(x,y+1)/\BB(x,y)$ & $1.340071766733$ & $0.8596363183029$ & \eqref{eq:beta2-unit}, same digits\\
$\BB(x+1,y)/\BB(x,y)$ & $2.147180690599$ & $1.083053099515$ & \eqref{eq:beta2-unit}, same digits\\
\hline
\end{tabular}
\end{table}

Table \ref{tab:num3} illustrates Proposition \ref{prop:Omega} and Theorem \ref{thm:logconvex-gamma}. The integral representations \eqref{eq:gamma1-def} and \eqref{eq:gamma2-def} were also compared with the closed forms \eqref{eq:gamma1-closed} and \eqref{eq:gamma2-closed} at $z\in\{0.5,1,2.5,4\}$, with relative errors below $2\times10^{-40}$; for the Jagannathan Srinivasa Rao case, $\GR(1/2)=1.53698378365855918$ and $\gR(1/2)=1.53698378365857518$, whose quotient is $\Omq(1/2)=1+1.04\times10^{-14}$.

\begin{table}[ht]
\centering
\caption{Left: the theta quotient $\Omq(1/2)-1$ compared with $4\varepsilon$, $\varepsilon=e^{-2\pi^{2}/\log(1/\fq)}$. Right: $\frac{d^{2}}{dx^{2}}\log\GR(x)$.}\label{tab:num3}
\begin{tabular}{lll}
\hline
$\fq$ & $\Omq(1/2)-1$ & $4\varepsilon$\\
\hline
$0.3$ & $3.03232\times10^{-7}$ & $3.03232\times10^{-7}$\\
$0.5$ & $1.71542\times10^{-12}$ & $1.71542\times10^{-12}$\\
$0.6$ & $6.60922\times10^{-17}$ & $6.60924\times10^{-17}$\\
\hline
\end{tabular}
\qquad
\begin{tabular}{lllll}
\hline
Case & $x=0.5$ & $x=2$ & $x=5$ & $x=10$\\
\hline
JS, $\phi_{1}=0.9$ & $4.5355$ & $0.28838$ & $-0.061224$ & $-0.10318$\\
Jackson & $4.6794$ & $0.42186$ & $0.056176$ & $0.0039746$\\
CJ, $\phi_{1}=1.25$ & $4.6998$ & $0.51204$ & $0.23769$ & $0.22329$\\
\hline
\end{tabular}
\end{table}

In the Jagannathan Srinivasa Rao case, the second logarithmic derivative of $\GR$ vanishes at $x_{*}=3.672228654\ldots$ and tends to $\log0.9=-0.10536\ldots$, while in the Chakrabarti Jagannathan case it tends to $\log1.25=0.22314\ldots$, in agreement with Theorem \ref{thm:logconvex-gamma}.

\section{Concluding remarks}\label{sec:conclusion}

We have developed the theory of the noncommutative beta functions attached to the $\mathcal{R}(p,q)-$deformed quantum algebras under the two root hypothesis $\Hphi$ and its gauged version, working throughout with the Jackson type series representation of the deformed integral. Two deformed gamma functions have been identified: the first kind function $\GR=\phi_{1}^{(z-1)(z-2)/2}\Gq$, characterised by a Bohr Mollerup type property and logarithmically convex on the positive half line exactly when $\phi_{1}\geq1$, and the second kind function $\gR=\Omq\GR$, obtained from a bilateral deformed integral through Ramanujan's ${}_{1}\psi_{1}$ summation, which differs from the first one by a one periodic theta quotient exponentially close to $1$. The first kind deformed beta function has been evaluated for complex arguments, with the valuation factor $\phi_{1}^{(y-1)(2x+y-2)/2}$; its general shift identity, its exchange defect $\phi_{1}^{(y-x)(x+y-1)/2}$ and its symmetrisation, unique up to a constant, have been established, together with a deformed partial fraction identity. The second kind deformed beta function on the half line has been evaluated in closed form; it pairs the two deformed gamma functions, and its exchange defect involves both deformation functions and the theta quotient, so that it remains noncommutative in the Jackson case. The first kind deformed beta function has been shown to be the kernel of a deformed Euler integral for a deformed Gauss hypergeometric function written with deformed Pochhammer symbols, which yields a deformed Gauss summation and a deformed Chu Vandermonde identity, and the two beta integrands have produced deformed beta and beta prime distributions with explicit moments. The deformed gamma function has been shown to satisfy a multiplication formula whose base is the deformed number $\pn{n}$ and a Wendel type estimate, while the first kind beta function grows or decays in its second argument like a Gaussian determined by $\phi_{1}$. The logarithmic convexity of normalised beta functions has been obtained from the series representation, with a deformed Tur\'an identity as discrete counterpart, and the gauge invariance of all relations has been used to treat the generalized $q-$Quesne and Hounkonnou Ngompe deformations. The numerical checks confirm every identity with thirty significant digits.

Several questions remain open. The deformed beta functions should allow the construction of deformed hypergeometric functions and of the associated Euler type integral representations, as well as deformed fractional operators of Riemann Liouville and Erd\'elyi Kober type, whose kernels are deformed beta kernels \cite{SoontharanonSitthiwirattham2020, Kiryoka}. Structure functions which do not satisfy the two root hypothesis, for instance those whose deformed numbers are combinations of more than two geometric sequences, require a different reduction, and the existence of a deformed Bohr Mollerup theorem in that setting is not clear. Finally, the second kind deformed gamma function suggests a deformed analogue of the theory of quasi constants of Askey type, in which the theta quotient $\Omq$ would play the role of a deformed modular factor.

\appendix

\section{Detailed computations}\label{app:computations}

\subsection{Exponent identities}\label{app:exponents}

Throughout, $E(x,y):=\tfrac{(y-1)(2x+y-2)}{2}$, so that $2E(x,y)=2xy+y^{2}-3y-2x+2$.

\paragraph{Beta gamma relation of first kind.}
From $\GR(z)=\phi_{1}^{(z-1)(z-2)/2}\Gq(z)$ and $(z-1)(z-2)=z^{2}-3z+2$,
\begin{eqnarray}
	(x-1)(x-2)+(y-1)(y-2)-(x+y-1)(x+y-2)&=&x^{2}+y^{2}-3x-3y+4-(x+y)^{2}\nn\\&+&3(x+y)-2
	=2-2xy,
\end{eqnarray}
so that $\GR(x)\GR(y)/\GR(x+y)=\phi_{1}^{1-xy}\Bq(x,y)$. Since $\BR=\phi_{1}^{1-x+(y-1)(y-2)/2}\Bq$ by Lemma \ref{lem:beta1-series} and Theorem \ref{thm:beta1}, the valuation exponent is
\begin{eqnarray}
	1-x+\frac{(y-1)(y-2)}{2}-1+xy&=&\frac{y^{2}-3y+2-2x+2xy}{2}\;=\;E(x,y).\nn
\end{eqnarray}

\paragraph{Exchange and symmetrisation.}
\begin{eqnarray}
	2E(x,y)-2E(y,x)&=&(y^{2}-3y-2x)-(x^{2}-3x-2y)\;=\;y^{2}-x^{2}-(y-x)\;=\;(y-x)(x+y-1),\nn\\
	x(x-1)+2E(x,y)&=&x^{2}+y^{2}+2xy-3x-3y+2\;=\;(x+y-1)(x+y-2).\nn
\end{eqnarray}

\paragraph{General shift.}
With $X=x+m$ and $Y=y+n$,
\begin{eqnarray}
	(Y-1)(2X+Y-2)&=&\big[(y-1)+n\big]\big[(2x+y-2)+(2m+n)\big]\nn\\
	&=&(y-1)(2x+y-2)+(y-1)(2m+n)+n(2x+y-2)+n(2m+n),\nn
\end{eqnarray}
so that $2[E(X,Y)-E(x,y)]=2m(y-1)+2n(x+y)+2mn+n^{2}-3n$, which is twice the exponent of \eqref{eq:shift-general}. For $(m,n)=(1,0)$, $(0,1)$ and $(1,1)$ one recovers $y-1$, $x+y-1$ and $x+2y-1$.

\paragraph{Second kind.}
From \eqref{eq:beta2-q} and $\Bq=\phi_{1}^{xy-1}\GR(x)\GR(y)/\GR(x+y)$,
\begin{eqnarray}
	1-x-\frac{(x+y)(x+y-1)}{2}+xy-1&=&\frac{-2x-x^{2}-2xy-y^{2}+x+y+2xy}{2}\;=\;-\frac{x^{2}+x+y^{2}-y}{2}.\nn
\end{eqnarray}
Using $\Kq(x)=\fq^{-x(x-1)/2}\Omq(x)$ and $\fq=\phi_{2}/\phi_{1}$, $\phi_{1}^{-(x^{2}+x+y^{2}-y)/2}\fq^{-x(x-1)/2}=\phi_{1}^{-x-y(y-1)/2}\phi_{2}^{-x(x-1)/2}$, which gives \eqref{eq:beta2-gamma2}.

\paragraph{Tur\'an identity.}
With $a=\phi_{1}$ and $b=\phi_{2}$,
\begin{eqnarray}
	(a^{x+1}-b^{x+1})(a^{x+y}-b^{x+y})-(a^{x}-b^{x})(a^{x+y+1}-b^{x+y+1})
	&=&-a^{x+1}b^{x+y}-b^{x+1}a^{x+y}+a^{x}b^{x+y+1}+b^{x}a^{x+y+1}\nn\\
	&=&a^{x}b^{x}\big(a^{y+1}-ab^{y}-ba^{y}+b^{y+1}\big)\nn\\
	&=&a^{x}b^{x}(a-b)(a^{y}-b^{y}),\nn
\end{eqnarray}
and division by $(a-b)^{2}$ gives \eqref{eq:turan-identity}.

\paragraph{Deformed Gauss summation.}\label{app:hyp}
With $E(b,y)=(y-1)(2b+y-2)/2$, $y_{1}=c-a-b$ and $y_{2}=y_{1}+a=c-b$,
\begin{eqnarray}
	(y_{1}-1)(2b+y_{1}-2)-(y_{1}+a-1)(2b+y_{1}+a-2)&=&-a(2b+y_{1}-2)-a(y_{1}+a-1)\nn\\
	&=&-a(2b+2y_{1}+a-3),\nn
\end{eqnarray}
and $2b+2y_{1}+a-3=2c-a-3$, so that $E(b,c-a-b)-E(b,c-b)=-a(2c-a-3)/2$. Then
$\tfrac{a(a-1)}{2}+a(c-a-b-1)-\tfrac{a(2c-a-3)}{2}=\tfrac{a}{2}\big[(a-1)+2c-2a-2b-2-2c+a+3\big]=-ab$. As a consistency check, with $f(z)=(z-1)(z-2)$ one has $f(c)+f(c-a-b)-f(c-a)-f(c-b)=2ab$, so that $\GR(c)\GR(c-a-b)/(\GR(c-a)\GR(c-b))=\phi_{1}^{ab}\,\Gq(c)\Gq(c-a-b)/(\Gq(c-a)\Gq(c-b))$, and \eqref{eq:gauss} is equivalent to the $\fq-$Gauss summation \cite[Section 1.5]{GasperRahman}.

\subsection{Second kind deformed beta function: the bilateral sum}\label{app:bilateral}

For completeness we detail the use of \eqref{eq:psi11} in Theorem \ref{thm:beta2}. With $a=-1$, $b=-\fq^{x+y}$ and $w=\fq^{x}$ one has $b/a=\fq^{x+y}$, $aw=-\fq^{x}$, $\fq/(aw)=-\fq^{1-x}$, $\fq/a=-\fq$ and $b/(aw)=\fq^{y}$, whence
\begin{eqnarray}
	\frac{(-\fq^{x+y};\fq)_{\infty}}{(-1;\fq)_{\infty}}\,
	\frac{(\fq,\fq^{x+y},-\fq^{x},-\fq^{1-x};\fq)_{\infty}}{(-\fq^{x+y},-\fq,\fq^{x},\fq^{y};\fq)_{\infty}}
	&=&\frac{(-\fq^{x},-\fq^{1-x};\fq)_{\infty}}{(-1,-\fq;\fq)_{\infty}}\,\frac{(\fq,\fq^{x+y};\fq)_{\infty}}{(\fq^{x},\fq^{y};\fq)_{\infty}}
	\;=\;\Kq(x)\,\frac{\Bq(x,y)}{1-\fq}.\nn
\end{eqnarray}
The same computation with $b=0$ gives the second kind deformed gamma function of Theorem \ref{thm:gamma2}, the factor $(-\fq^{x+y};\fq)_{\infty}$ being then absent.

\subsection{Integer second argument: derivation by iteration}\label{app:integer}

For $y=n+1$, $n\in\mathbb{N}\cup\{0\}$, the beta gamma relation can be obtained without the $\fq-$binomial theorem. Starting from $\BR(x,1)=1/\pn{x}$ (Proposition \ref{prop:integral}(ii)) and iterating the unit shift \eqref{eq:shift-y}, which was derived in Theorem \ref{thm:unit-shifts} from the integration by parts and Pascal relations only, one obtains
\begin{eqnarray}
	\BR(x,n+1)&=&\prod_{j=1}^{n}\phi_{1}^{\,x+j-1}\frac{\pn{j}}{\pn{x+j}}\;\BR(x,1)
	\;=\;\phi_{1}^{\,nx+\frac{n(n-1)}{2}}\,\frac{\pf{n}}{\prod_{j=0}^{n}\pn{x+j}} .\nn
\end{eqnarray}
Since $nx+n(n-1)/2=n(2x+n-1)/2=E(x,n+1)$ and $\GR(x)/\GR(x+n+1)=1/\prod_{j=0}^{n}\pn{x+j}$, this is \eqref{eq:beta1-gamma} at $y=n+1$. Expanding the deformed power by Lemma \ref{lem:powers}(iii) gives in addition the deformed partial fraction identity \eqref{eq:partial-fractions}, so that the beta gamma relation for integer second argument is equivalent to that identity. The extension to complex $y$ requires an input beyond the recursions, here the $\fq-$binomial theorem; the recursions determine $\BR(x,\cdot)$ only up to a one periodic factor, exactly as the functional equation determines $\GR$ only up to such a factor (compare Theorems \ref{thm:BM} and \ref{thm:gamma2}).

\subsection{An independent check of the unit shifts by the integral}\label{app:direct}

For $y=2$, $(1\ominus\phi_{2}t)^{1}_{\mathcal{R}(p,q)}=1-\phi_{2}t$ and Proposition \ref{prop:integral}(ii) give directly
\begin{eqnarray}
	\BR(x,2)&=&\frac{1}{\pn{x}}-\frac{\phi_{2}}{\pn{x+1}}\;=\;\frac{\pn{x+1}-\phi_{2}\pn{x}}{\pn{x}\pn{x+1}}\;=\;\frac{\phi_{1}^{x}}{\pn{x}\pn{x+1}},\nn
\end{eqnarray}
by the addition formula $\pn{x+1}=\phi_{1}^{x}+\phi_{2}\pn{x}$. This agrees with \eqref{eq:shift-y} at $y=1$, namely $\BR(x,2)=\phi_{1}^{x}\pn{1}\pn{x+1}^{-1}\BR(x,1)$, and with Corollary \ref{cor:special-values}, since $E(x,2)=x$ and $\pf{1}=1$.

\end{document}